\documentclass[11pt]{article}

\usepackage[a4paper]{geometry}
\usepackage{authblk}

\usepackage{amsmath, amssymb, amsthm, mathtools, thmtools}
\usepackage{bm}
\usepackage{comment}
\usepackage{commath}
\usepackage[shortlabels]{enumitem}
\usepackage{etoolbox}
\usepackage{hyperref}
\usepackage{mathtools}
\usepackage[authoryear]{natbib}
\usepackage{xcolor}

\declaretheoremstyle[
    spaceabove = \topsep,
    spacebelow = \topsep,
    headfont = \bfseries,
    headpunct = {},
    notefont = \normalfont,
    notebraces = {(}{)}, 
    bodyfont = \normalfont,
    postheadspace = 1em
]{theoremstyle}

\declaretheorem[numberwithin = section, style = theoremstyle]{theorem}
\declaretheorem[sibling = theorem, style = theoremstyle]{corollary}
\declaretheorem[sibling = theorem, style = theoremstyle]{definition}
\declaretheorem[sibling = theorem, style = theoremstyle]{example}
\declaretheorem[sibling = theorem, style = theoremstyle]{lemma}
\declaretheorem[sibling = theorem, style = theoremstyle]{remark}
\declaretheorem[sibling = theorem, style = theoremstyle]{proposition}
\declaretheorem[sibling = theorem, style = theoremstyle]{setting}

\newcommand\declaresymbol[2]{\newcommand{#1}{\TextOrMath{$#2$\xspace}{#2}}}

\declaresymbol\banachspace F
\declaresymbol\banachspaceelement f

\declaresymbol\semigroup T
\declaresymbol\generator A

\declaresymbol\probabilityspace\Omega
\declaresymbol\eventsystem{\mathcal A}
\declaresymbol\probabilitymeasure{\operatorname P}
\declaresymbol\expectation{\operatorname E}
\declaresymbol\shift\theta
\declaresymbol\outcome\omega
\declaresymbol\event A
\declaresymbol\secondevent B
\declaresymbol\probabilityfamily{\mathcal P}

\declaresymbol\dominatingvariable D

\declaresymbol\filtration{\mathcal F}
\declaresymbol\secondfiltration{\mathcal G}
\declaresymbol\thirdfiltration{\mathcal H}
\declaresymbol\filtrationset F

\declaresymbol\measurablespace E 
\declaresymbol\point x
\declaresymbol\secondpoint y
\declaresymbol\measurableset B
\declaresymbol\measurablesystem{\mathcal E}
\declaresymbol\generatorsystem{\mathcal G}
\declaresymbol\generatorset G

\declaresymbol\measure\nu

\declaresymbol\topology\eta

\declaresymbol\deadstate\Delta

\declaresymbol\expparameter\alpha
\declaresymbol\expvariable\xi

\declaresymbol\integrand f
\declaresymbol\secondintegrand g

\declaresymbol\indexset I
\declaresymbol\indexsetelement i
\declaresymbol\secondindexsetelement j
\declaresymbol\thirdindexsetelement k
\declaresymbol\fourthindexsetelement l

\declaresymbol\additivefunctional A
\declaresymbol\multiplicativefunctional M
\declaresymbol\lebesguestieltjes\alpha

\declaresymbol\closedsubspace C
\declaresymbol\core D

\declaresymbol\process X
\declaresymbol\localprocess Y
\declaresymbol\killedlocalprocess Z
\declaresymbol\concatenatedprocess X

\declaresymbol\kernel\kappa
\declaresymbol\localsemigroup\kappa
\declaresymbol\killedlocalsemigroup Q
\declaresymbol\concatenatedsemigroup Q

\declaresymbol\localgenerator L
\declaresymbol\killedlocalgenerator K
\declaresymbol\globalgenerator G
\declaresymbol\concatenatedgenerator A

\declaresymbol\killingrate c

\declaresymbol\probabilitykernel\varpi
\declaresymbol\transferkernel T
\declaresymbol\regenerationdistribution\mu

\declaresymbol\projectionvariable V

\declaresymbol\probabilityprojection p

\declaresymbol\timedomain T
\declaresymbol\discretetimedomain{\mathbb N_0}
\declaresymbol\timepoint t
\declaresymbol\prevtimepoint s
\declaresymbol\prevprevtimepoint r
\declaresymbol\discretetimepoint n

\declaresymbol\lifetime\tau
\declaresymbol\lifetimesum\sigma

\setlistdepth{9}
\renewlist{itemize}{itemize}{9}

\newenvironment{listclaim}{%
    
    \stepcounter{listclaimlevel}%
    \refstepcounter{listclaim\roman{listclaimlevel}}%
    \normalfont\trivlist
    \item[\hskip\labelsep{\bfseries Claim \csuse{thelistclaim\roman{listclaimlevel}}:}]\itshape\ignorespaces\normalfont 
}{%
    \endtrivlist\addtocounter{listclaimlevel}{-1}%
}

\newcounter{listclaimlevel}
\newcounter{listclaimi}[theorem] 
\newcounter{listclaimii}[listclaimi]
\newcounter{listclaimiii}[listclaimii]
\newcounter{listclaimiv}[listclaimiii]
\newcounter{listclaimv}[listclaimiv]
\newcounter{listclaimvi}[listclaimv]

\let\originalleft\left
\let\originalright\right
\renewcommand{\left}{\mathopen{}\mathclose\bgroup\originalleft}
\renewcommand{\right}{\aftergroup\egroup\originalright}

\newcommand{\underbraceinsidebracket}[6]{%
    \left.\vphantom{#2 #3 #5}\right#1 %
    #2\underbrace{#3}_{#4}#5 %
    \left.\vphantom{#2 #3 #5}\right#6
}

\newcommand\expectationmid{\;\middle|\;}

\newcommand\largebracket[1]{%
    \scalebox{1.4}{#1}
}

\newif\ifincludeproofs
\includeproofsfalse

\ifincludeproofs
    
\else
    \excludecomment{conditionalproof}
\fi

\makeatletter
\newcommand{\setword}[2]{%
  \phantomsection
  #1\def\@currentlabel{\unexpanded{#1}}\label{#2}%
}
\makeatother

\newcommand*\sq{\mathbin{\vcenter{\hbox{\rule{.6ex}{.6ex}}}}}

\usepackage{nomencl}
\makenomenclature

\title{Concatenation of Markov processes}

\author{
    \href{https://people.mpi-inf.mpg.de/~sholl/}{Sascha~Holl}\\
	Max Planck Institute for Informatics\\
	\texttt{sholl@mpi-inf.mpg.de}
}

\date{\today}

\begin{document}

\maketitle

In this paper, we investigate the concatenation of Markov processes. Our primary concern is to utilize processes constructed in this manner for Monte Carlo integration. To enable this using conventional methods, it is essential to demonstrate the Markov property and invariance with respect to a given target distribution. We provide mild sufficient conditions for this.

In \autoref{sec:stochastic-base}, we first introduce the theoretical concepts necessary for our elaboration. The central notions of additive and multiplicative functionals are defined similarly to \citep[Section~IV.1 and Section~III.1]{sharpe1988markov}.

The concatenation of processes is clearly only meaningful if the processes to be concatenated are executed for finite times. From the perspective of an individual process, it is "killed" before the execution of the next process. The aforementioned multiplicative functionals arise naturally from the rule describing the desired type of killing. We present the corresponding construction in \autoref{sec:killing-by-a-multiplicative-functional}. This approach has already been described in the literature in \citep[Section~III.2 and Section~III.3]{blumenthal1968markov}. In contrast to their presentation, we avoid unnecessary restrictions, allowing for a less obscure representation.

In \autoref{sec:concatenation}, we finally describe how the concatenation of Markov processes can be formally constructed. Similar to the previous section, this construction has also been considered in the literature in \citep[Section~II.14]{sharpe1988markov}. Once again, our presentation differs by making fewer restrictive assumptions and by being tailored towards the identification of an invariant measure.

Our main result is the identification of the generator of the concatenation of Markov processes. This result provides the theoretical foundation for Monte Carlo methods based on this construction. An initial instance of such methods has already been considered in \citep{wang2020thesis}, \citep{wang2021regeneration}, and \citep{mckimm2022restore}. The theory presented in this work significantly generalizes and rigorously establishes the arguments on which these works are based.

\section{Stochastic base, multiplicative and additive functional, realization of a Markov semigroup}\label{sec:stochastic-base}

Following, for example, \citep[Section~IX.4]{cinlar2011probability}, we deal with Markov processes in a modern setting. Fundamental to this is the notion of a \emph{stochastic base}:

\begin{definition}\label{def:stochastic-base}
    $(\probabilityspace,\eventsystem,(\filtration_\timepoint)_{\timepoint\ge0},(\shift_\timepoint)_{\timepoint\ge0})$ is called \textbf{stochastic base} $:\Leftrightarrow$ \begin{enumerate}[(i)]
        \item $(\probabilityspace,\eventsystem)$ is a measurable space.
        \item $(\filtration_\timepoint)_{\timepoint\ge0}$ is a filtration on $(\probabilityspace,\eventsystem)$.
        \item $\shift_0=\operatorname{id}_\probabilityspace$ and $\shift_\timepoint:\probabilityspace\to\probabilityspace$ is $(\eventsystem,\eventsystem)$-measurable for all $\timepoint>0.$
    \end{enumerate}
    \begin{flushright}
        $\square$
    \end{flushright}
\end{definition}

$\shift_\prevtimepoint$ will function as a time shift operator. Its existence is a condition of richness on $\probabilityspace$. If, for example, $\probabilityspace=\measurablespace^{[0,\;\infty)}$ for some vector space $\measurablespace$, we may define \begin{equation}
    \left(\left(\shift_\prevtimepoint\right)(\omega)\right)(\timepoint):=\omega(\prevtimepoint+\timepoint)-\omega(\prevtimepoint)\;\;\;\text{for }\timepoint\ge0\text{ and }\omega\in\probabilityspace.
\end{equation} In probability theory, we usually have no particular interest in the measurable space $(\probabilityspace,\eventsystem)$ itself and hence do not hesitate to assume that it is rich enough to admit any object we desire. However, the interested reader can find a more detailed elaboration on the existence of $\left(\shift_\prevtimepoint\right)_{\prevtimepoint\ge0}$ in \citep[Section~VII-3]{cinlar2011probability}.

We now fix a stochastic base $(\probabilityspace,\eventsystem,(\filtration_\timepoint)_{\timepoint\ge0},(\shift_\timepoint)_{\timepoint\ge0})$ and $\probabilityfamily\subseteq\mathcal M_1(\probabilityspace,\eventsystem)$. By $\mathcal M_1(\probabilityspace,\eventsystem)$ and $\mathcal M_{\le1}(\probabilityspace,\eventsystem)$, we denote the set of \setword{probability}{inline:probability-measures} and \setword{sub-probability}{inline:sub-probability-measures} measures on $(\probabilityspace,\eventsystem)$, respectively.

Central to our construction in \autoref{sec:killing-by-a-multiplicative-functional} is the notion of a \emph{multiplicative functional}.

\begin{definition}[multiplicative functional]\label{def:multiplicative-functional}
    $(\multiplicativefunctional_\timepoint)_{\timepoint\ge0}$ is called \textbf{multiplicative functional on} $\bm{(\probabilityspace,\eventsystem,(\filtration_\timepoint)_{\timepoint\ge0},(\shift_\timepoint)_{\timepoint\ge0},\probabilityfamily)}$ $:\Leftrightarrow$ \begin{enumerate}[(i)]
        \item $(\multiplicativefunctional_\timepoint)_{\timepoint\ge0}$ is a $[0,1]$-valued $(\filtration_\timepoint)_{\timepoint\ge0}$-adapted process on $(\probabilityspace,\eventsystem)$.
        \item Let $\prevtimepoint,\timepoint\ge0$ and $\probabilitymeasure\in\probabilityfamily$ $\Rightarrow$ \begin{equation}\label{eq:multiplicative-functional}
            \multiplicativefunctional_{\prevtimepoint+t}=\multiplicativefunctional_\prevtimepoint(\multiplicativefunctional_\timepoint\circ\shift_\prevtimepoint)\;\;\;\probabilitymeasure\text{-almost surely}.
        \end{equation}
    \end{enumerate}
    \textit{Remark}: \begin{equation}
        \multiplicativefunctional_0=\multiplicativefunctional_{0+0}=\multiplicativefunctional_0(\multiplicativefunctional_0\circ\shift_0)=\multiplicativefunctional_0^2\;\;\;\probabilitymeasure\text{-almost surely}
    \end{equation} and hence \begin{equation}\label{eq:multiplicative-functional-inistial-state}
        \multiplicativefunctional_0\in\{0,1\}\;\;\;\probabilitymeasure\text{-almost surely}
    \end{equation} for all $\probabilitymeasure\in\probabilityfamily$.
    \begin{flushright}
        $\square$
    \end{flushright}
\end{definition}

Multiplicative functionals naturally emerge when we \emph{kill} a process according to a given rule as it is done in \autoref{sec:killing-by-a-multiplicative-functional}. For further elaboration on them, we refer to \citep[Chapter~III]{blumenthal1968markov} and \citep[Chapter~VII]{sharpe1988markov}.

One important instance of a multiplicative functional is induced by the analogously defined concept of an \emph{additive functional}:

\begin{definition}[additive functional]\label{def:additive-functional}
    $(\additivefunctional_\timepoint)_{\timepoint\ge0}$ is called \textbf{additive functional on} $\bm{(\probabilityspace,\eventsystem,(\filtration_\timepoint)_{\timepoint\ge0},(\shift_\timepoint)_{\timepoint\ge0},\probabilityfamily)}$ $:\Leftrightarrow$ \begin{enumerate}[(i)]
        \item $(\additivefunctional_\timepoint)_{\timepoint\ge0}$ is a $[0,\infty]$-valued $(\filtration_\timepoint)_{\timepoint\ge0}$-adapted process on $(\probabilityspace,\eventsystem)$.
        \item Let $\prevtimepoint,\timepoint\ge0$ and $\probabilitymeasure\in\probabilityfamily$ $\Rightarrow$ \begin{equation}\label{eq:additive-functional}
            \additivefunctional_{\prevtimepoint+\timepoint}=\additivefunctional_\prevtimepoint+\additivefunctional_\timepoint\circ\shift_\prevtimepoint\;\;\;\probabilitymeasure\text{-almost surely}.
        \end{equation}
    \end{enumerate}
    \textit{Remark}: \begin{equation}
        A_0=A_{0+0}=A_0+A_0\circ\shift_0=A_0+A_0\;\;\;\probabilitymeasure\text{-almost surely}
    \end{equation} and hence \begin{equation}
        A_0=0\;\;\;\probabilitymeasure\text{-almost surely}
    \end{equation} for all $\probabilitymeasure\in\probabilityfamily$.
    \begin{flushright}
        $\square$
    \end{flushright}
\end{definition}

For an extended discussion on additive functionals, we refer to \citep[Chapter~IV]{blumenthal1968markov} and \citep[Chapter~VIII]{sharpe1988markov}. The mentioned induced multiplicative functional is constructed in the following way:

\begin{proposition}[multiplicative functional induced by additive functional]\label{prop:multiplicative-functional-induced-by-additive-functional}
    Let $(\additivefunctional_\timepoint)_{\timepoint\ge0}$ be an additive functional on $(\probabilityspace,\eventsystem,(\filtration_\timepoint)_{\timepoint\ge0},(\shift_\timepoint)_{\timepoint\ge0},\probabilityfamily)$ $\Rightarrow$ \begin{equation}
        \multiplicativefunctional_\timepoint:=e^{-\additivefunctional_\timepoint}\;\;\;\text{for }\timepoint\ge0
    \end{equation} is a multiplicative functional on $(\probabilityspace,\eventsystem,(\filtration_\timepoint)_{\timepoint\ge0},(\shift_\timepoint)_{\timepoint\ge0},\probabilityfamily)$. If $(\additivefunctional_\timepoint)_{\timepoint\ge0}$ is nondecreasing, then $(\multiplicativefunctional_\timepoint)_{\timepoint\ge0}$ is nonincreasing.
    \begin{proof}[Proof\textup:\nopunct]
        \leavevmode
        \begin{itemize}[$\circ$]
            \item $e^{-x}\in[0,1]$ for all $x\in[0,\infty]$ $\Rightarrow (\multiplicativefunctional_\timepoint)_{\timepoint\ge0}$ is a $[0,1]$-valued $(\filtration_\timepoint)_{\timepoint\ge0}$-adapted process on $(\probabilityspace,\eventsystem)$.
            \item Let $\prevtimepoint,\timepoint\ge0$ and $\probabilitymeasure\in\probabilityfamily$ $\Rightarrow$ \begin{equation}
                \multiplicativefunctional_{\prevtimepoint+\timepoint}\stackrel{\text{def}}=e^{-\additivefunctional_{\prevtimepoint+\timepoint}}\stackrel{\eqref{eq:additive-functional}}=e^{-(\additivefunctional_\prevtimepoint+\additivefunctional_\timepoint\:\circ\:\shift_\prevtimepoint)}=\underbrace{e^{-\additivefunctional_\prevtimepoint}}_{\stackrel{\text{def}}=\:\multiplicativefunctional_\prevtimepoint}\underbraceinsidebracket{(}{}{e^{-\additivefunctional_\timepoint}}{\stackrel{\text{def}}=\;\multiplicativefunctional_\timepoint}{\circ\;\shift_\prevtimepoint}{)}
            \end{equation} $\probabilitymeasure$-almost surely.
        \end{itemize}
    \end{proof}
\end{proposition}

While a stopping time $\lifetime$ cannot "look in the future", it might have a certain "memory of the past". If we concatenate two Markov processes by terminating the execution of the first one at time $\lifetime$ and then starting the second process in its place, this memory destroys the Markov property of the combined process. For that reason, we need to come up with a suitable "memorylessness" property:

\begin{definition}[terminal time]\label{def:terminal-time}
    $\lifetime$ is called \textbf{terminal time on} $\bm{(\probabilityspace,\eventsystem,(\filtration_\timepoint)_{\timepoint\ge0},(\shift_\timepoint)_{\timepoint\ge0},\probabilityfamily)}$ $:\Leftrightarrow$ \begin{enumerate}[(i)]
        \item $\lifetime$ is an $(\filtration_\timepoint)_{\timepoint\ge0}$ stopping time on $(\probabilityspace,\eventsystem)$;
        \item Let $s\ge0$ and $\probabilitymeasure\in\mathcal P$ $\Rightarrow$ \begin{equation}\label{eq-terminal-time}
            \prevtimepoint+\lifetime\circ\shift_\prevtimepoint=\lifetime\;\;\;\probabilitymeasure\text{-almost surely on }\{\:\prevtimepoint<\lifetime\:\}.
        \end{equation}
    \end{enumerate}
\end{definition}

A terminal time gives rise to a multiplicative functional:

\begin{proposition}[multiplicative functional induced by stopping time]\label{ex:multiplicative-functional-induced-by-stopping-time}
	If $\lifetime$ be a terminal time on $(\probabilityspace,\eventsystem,(\filtration_\timepoint)_{\timepoint\ge0},(\shift_t)_{\timepoint\ge0},\mathcal P)$, then \begin{equation}
		\multiplicativefunctional_\timepoint:=1_{[0,\:\lifetime)}(\timepoint)\;\;\;\text{for }\timepoint\ge0
	\end{equation} is a multiplicative functional on $(\probabilityspace,\eventsystem,(\filtration_\timepoint)_{\timepoint\ge0},(\shift_t)_{\timepoint\ge0},\mathcal P)$.
	\begin{proof}
        \leavevmode
		\begin{itemize}[$\circ$]
			\item Let $\prevtimepoint,\timepoint\ge0$ $\Rightarrow$ \begin{equation}
				\begin{split}
					&\multiplicativefunctional_{\prevtimepoint+\timepoint}(\omega)=1\\\Leftrightarrow&\:(\omega,\prevtimepoint+\timepoint)\in[0,\lifetime)\\\Leftrightarrow&\:\prevtimepoint<\lifetime(\omega)\text{ and }\prevtimepoint+\timepoint<\lifetime(\omega)
				\end{split}
			\end{equation} and \begin{equation}
				\begin{split}
					&\prevtimepoint<\lifetime(\omega)\text{ and }\timepoint<\lifetime\left(\shift_\prevtimepoint(\omega)\right)\\\Leftrightarrow&\:(\omega,\prevtimepoint)\in[0,\lifetime)\text{ and }(\shift_\prevtimepoint(\omega),\timepoint)\in[0,\lifetime)\\\Leftrightarrow&\:\multiplicativefunctional_\prevtimepoint(\omega)\multiplicativefunctional\left(\shift_\prevtimepoint(\omega)\right)=1
				\end{split}
			\end{equation} for all $\omega\in\omega$.
			\item $\lifetime$ is a terminal time on $(\probabilityspace,\eventsystem,(\filtration_\timepoint)_{\timepoint\ge0},(\shift_\timepoint)_{\timepoint\ge0},\probabilitymeasure)$ $\Rightarrow$ Claim.
		\end{itemize}
	\end{proof}
\end{proposition}

We now go over to consider general processes living on a stochastic basis. For that, we fix a measurable space $(\measurablespace,\measurablesystem)$ on which the process will evolve.

\begin{definition}[process on stochastic base]\label{def:process-on-stochastic-process}
    $(\process_\timepoint)_{\timepoint\ge0}$ is called $\bm{(\measurablespace,\measurablesystem)}$\textbf{-valued process on} $\bm{(\probabilityspace,\eventsystem,(\filtration_\timepoint)_{\timepoint\ge0},(\shift_\timepoint)_{\timepoint\ge0},\probabilityfamily)}$ $:\Leftrightarrow$ \begin{enumerate}[(i)]
        \item $(\process_\timepoint)_{\timepoint\ge0}$ is an $(\measurablespace,\measurablesystem)$-valued $(\filtration_\timepoint)_{\timepoint\ge0}$-adapted process on $(\probabilityspace,\eventsystem)$.
        \item Let $\prevtimepoint,\timepoint\ge0$ and $\probabilitymeasure\in\probabilityfamily$ $\Rightarrow$ \begin{equation}\label{eq:process-on-stochastic-base}
            \process_{\prevtimepoint+\timepoint}=\process_\timepoint\circ\shift_\prevtimepoint\;\;\;\probabilitymeasure\text{-almost surely}.
        \end{equation}
    \end{enumerate} In that case, $(\process_\timepoint)_{\timepoint\ge0}$ is called \textbf{progressive} $:\Leftrightarrow$ $(\process_\timepoint)_{\timepoint\ge0}$ is $(\filtration_\timepoint)_{\timepoint\ge0}$-progressive.
    \begin{flushright}
        $\square$
    \end{flushright}
\end{definition}

Each such process does induce an additive functional:

\begin{proposition}[additive functional induced by process]\label{prop:additive-functional-induced-by-process}
    \normalfont Let $(\process_\timepoint)_{\timepoint\ge0}$ be an $(\measurablespace,\measurablesystem)$-valued progressive process on $(\probabilityspace,\eventsystem,(\filtration_\timepoint)_{\timepoint\ge0},(\shift_\timepoint)_{\timepoint\ge0},\probabilityfamily)$ and $c:E\to[0,\infty)$ be $\measurablespace$-measurable. Then, \begin{equation}
        \additivefunctional_\timepoint:=\int_0^\timepoint\killingrate(\process_\prevtimepoint)\dif\prevtimepoint\;\;\;\text{for }\timepoint\ge0
    \end{equation} is an additive functional on $(\probabilityspace,\eventsystem,(\filtration_\timepoint)_{\timepoint\ge0},(\shift_\timepoint)_{\timepoint\ge0},\probabilityfamily)$. Moreover, $\additivefunctional_0=0$ and $(\additivefunctional)_{\timepoint\ge0}$ is nondecreasing.
    \begin{proof}[Proof\textup:\nopunct]
        \leavevmode
        \begin{itemize}[$\circ$]
            \item $c\ge0$ $\Rightarrow$ $(\additivefunctional)_{\timepoint\ge0}$ is $[0,\infty]$-valued.
            \item $(\process_\timepoint)_{\timepoint\ge0}$ is $(\filtration_\timepoint)_{\timepoint\ge0}$-progressive $\Rightarrow$ $(\additivefunctional)_{\timepoint\ge0}$ is $(\filtration_\timepoint)_{\timepoint\ge0}$-adapted.
            \item Let $\prevtimepoint,\timepoint\ge0$ and $\probabilitymeasure\in\probabilityfamily$ $\Rightarrow$ \begin{equation}
                \additivefunctional_{\prevtimepoint+\timepoint}\stackrel{\text{def}}=\int_0^{\prevtimepoint+\timepoint}\killingrate(\process_\prevprevtimepoint)\dif\prevprevtimepoint\stackrel{\eqref{eq:process-on-stochastic-base}}=\underbrace{\int_0^\prevtimepoint\killingrate(\process_\prevprevtimepoint)\dif\prevprevtimepoint}_{\stackrel{\text{def}}=\:\multiplicativefunctional_\prevtimepoint}+\underbrace{\int_0^\timepoint\killingrate(\process_\prevprevtimepoint)\dif\prevprevtimepoint}_{\stackrel{\text{def}}=\:\multiplicativefunctional}\circ\:\shift_\prevtimepoint
            \end{equation} $\probabilitymeasure$-almost surely.
        \end{itemize}
    \end{proof}
\end{proposition}

We now consider processes $(\process_\timepoint)_{\timepoint\ge0}$ together with a family of $(\probabilitymeasure_\point)_{\point\in\measurablespace}$ such that $(\process_\timepoint)_{\timepoint\ge0}$ is almost surely started at $\point\in\measurablespace$ under $\probabilitymeasure_\point$. For a purely technical reason, we need to assume that \begin{equation}
    \{\point\}\in\measurablesystem\;\;\;\text{for all }\point\in\measurablesystem.
\end{equation} Otherwise, the following definition would not be well-posed.


\begin{definition}[normal process on stochastic base]~
    \normalfont$((\process_\timepoint)_{\timepoint\ge0},(\probabilitymeasure_\point)_{\point\in\measurablespace})$ is called \textbf{normal} $\bm{(\measurablespace,\measurablesystem)}$\textbf{-valued process on} $\bm{(\probabilityspace,\eventsystem,(\filtration_\timepoint)_{\timepoint\ge0},(\shift_\timepoint)_{\timepoint\ge0})}$ $:\Leftrightarrow$ \begin{enumerate}[(i)]
        \item $(\probabilitymeasure_\point)_{\point\in\measurablespace}\subseteq\mathcal M_1(\probabilityspace,\eventsystem)$.
        \item $(\process_\timepoint)_{\timepoint\ge0}$ is an $(\measurablespace,\measurablesystem)$-valued process on $(\probabilityspace,\eventsystem,(\filtration_\timepoint)_{\timepoint\ge0},(\shift_\timepoint)_{\timepoint\ge0},(\probabilitymeasure_\point)_{\point\in\measurablespace})$ with \begin{equation}
            \process_0=\point\;\;\;\probabilitymeasure_\point\text{-almost surely for all }\point\in\measurablespace.
        \end{equation}
    \end{enumerate}
    \begin{flushright}
        $\square$
    \end{flushright}
\end{definition}

Since it plays a crucial role in our constructions made in \autoref{sec:killing-by-a-multiplicative-functional}, it is worth elaborating on the additive functional defined in \autoref{prop:additive-functional-induced-by-process} in more detail:

\begin{remark}\label{rem:additive-functional-induced-by-process}
    Let $((\process_\timepoint)_{\timepoint\ge0},(\probabilitymeasure_\point)_{\point\in\measurablespace})$ be an $(\measurablespace,\measurablesystem)$-valued progressive process on $(\probabilityspace,\eventsystem,(\filtration_\timepoint)_{\timepoint\ge0},(\shift_\timepoint)_{\timepoint\ge0})$ and \begin{equation}
        \additivefunctional_\timepoint:=\int_0^\timepoint\killingrate(X_\prevtimepoint)\dif\prevtimepoint\;\;\;\text{for }\timepoint\ge0
    \end{equation} for some $\measurablespace$-measurable $\killingrate:\measurablespace\to[0,\infty)$ denote the additive functional on $(\probabilityspace,\eventsystem,(\filtration_\timepoint)_{\timepoint\ge0},(\shift_\timepoint)_{\timepoint\ge0},(\probabilitymeasure_\point)_{\point\in\measurablespace})$ induced by $(\process_\timepoint)_{\timepoint\ge0}$ and $\killingrate$ introduced in \autoref{prop:additive-functional-induced-by-process}. Let $\point\in\measurablespace$. Assume $(\killingrate\circ\process)(\outcome)$ is (right-)continuous at $0$ for $\probabilitymeasure_\point$-almost all $\outcome\in\probabilityspace$ $\Rightarrow$ \begin{equation}
        \left.\frac{\dif}{\dif t}\additivefunctional_\timepoint(\outcome)\right|_{\timepoint=0+}=\killingrate(\point)
    \end{equation} for $\probabilitymeasure_\point$-almost all $\outcome\in\probabilityspace$. Let \begin{equation}
        \multiplicativefunctional:=e^{-\additivefunctional}
    \end{equation} denote the multiplicative functional on $(\probabilityspace,\eventsystem,(\filtration_\timepoint)_{\timepoint\ge0},(\shift_\timepoint)_{\timepoint\ge0},(\probabilitymeasure_\point)_{\point\in\measurablespace})$ induced by $(\additivefunctional)_{\timepoint\ge0}$ introduced in \autoref{prop:multiplicative-functional-induced-by-additive-functional} $\Rightarrow$ $\multiplicativefunctional(\outcome)$ is (right-)differentiable at $0$ with \begin{equation}
        \left.\frac{\dif}{\dif´\timepoint}\multiplicativefunctional_\timepoint(\outcome)\right|_{\timepoint=0+}=-\killingrate(\point)
    \end{equation} for $\probabilitymeasure_\point$-almost all $\outcome\in\probabilityspace$. Assume $\killingrate$ is bounded $\Rightarrow$ \begin{equation}
        0\le\frac{\multiplicativefunctional_0-\multiplicativefunctional_\timepoint}\timepoint=\frac{1-e^{-\additivefunctional}}\timepoint\le\frac{\additivefunctional}\timepoint=\frac1\timepoint\int_0^\timepoint\killingrate(\process_\prevtimepoint)\dif\prevtimepoint\le\left\|c\right\|_\infty<\infty
    \end{equation} for all $\timepoint>0$ and hence \begin{equation}
        \left.\frac{\dif}{\dif\timepoint}\expectation_\point\left[\multiplicativefunctional\right]\right|_{\timepoint=0+}=-\killingrate(\point)
    \end{equation} by Lebesgue’s dominated convergence theorem.
    \begin{flushright}
        $\square$
    \end{flushright}
\end{remark}

To the end of this section, we fix a Markov semigroup $(\kernel_\timepoint)_{\timepoint\ge0}$ on $(\measurablespace,\measurablesystem)$. As shown in \citep[Theorem~11.4]{kallenberg2021probability}, as long as $(\measurablespace,\measurablesystem)$ is nice enough, we can always construct a Markov process $(\process_\timepoint)_{\timepoint\ge0}$ with transition semigroup $(\kernel_\timepoint)_{\timepoint\ge0}$. Consequentially, the concrete shape of $(\process_\timepoint)_{\timepoint\ge0}$ is of minor interest; it is only the distribution determined by $(\kernel_\timepoint)_{\timepoint\ge0}$ which is of importance. From the proof of \citep[Theorem~11.4]{kallenberg2021probability}, it is apparent that we can select a family of probability measures $(\probabilitymeasure_\point)_{\point\in\measurablespace}$ such that $(\process_\timepoint)_{\timepoint\ge0}$ is almost surely started at $\point\in\measurablespace$ under $\probabilitymeasure_\point$. Since this will be important in the constructions made in \autoref{sec:killing-by-a-multiplicative-functional} and \autoref{sec:concatenation}, we will assume that the dependence of $\probabilitymeasure_\point$ on $\point$ is $\measurablesystem$-measurable. altogether, this gives us the following definition of a \textit{realization} of a Markov semigroup:

\begin{definition}\label{def:realization-of-markov-semigroup}
    $((\process_\timepoint)_{\timepoint\ge0},(\probabilitymeasure_\point)_{\point\in\measurablespace})$ is called \textbf{realization of} $\bm{(\kernel_\timepoint)_{\timepoint\ge0}}$ \textbf{on} $\bm{(\probabilityspace,\eventsystem,(\filtration_\timepoint)_{\timepoint\ge0},(\shift_\timepoint)_{\timepoint\ge0})}$ $:\Leftrightarrow$ \begin{enumerate}[(i)]
        \item$((\process_\timepoint)_{\timepoint\ge0},(\probabilitymeasure_\point)_{\point\in\measurablespace})$ is a normal $(\measurablespace,\measurablesystem)$-valued process on $(\probabilityspace,\eventsystem,(\filtration_\timepoint)_{\timepoint\ge0},(\shift_\timepoint)_{\timepoint\ge0})$.
        \item Let $\prevtimepoint,\timepoint\ge0$ and $(\point,\integrand)\in\measurablespace\times\measurablespace_b$ $\Rightarrow$ \begin{equation}\label{eq:realization-of-markov-semigroup}
            \expectation_\point\left[\integrand(\process_{\prevtimepoint+\timepoint})\mid\filtration_\prevtimepoint\right]=(\kernel\integrand)(\process_\prevtimepoint).
        \end{equation}
        \item \begin{equation}\label{eq:markov-kernel-induced-by-realization-of-markov-semigroup}
            \measurablespace\to[0,1]\;,\;\;\;\point\mapsto\probabilitymeasure_\point[\eventsystem]
        \end{equation} is $\measurablespace$-measurable for all $\eventsystem\in\eventsystem$.\footnote{\eqref{eq:realization-of-markov-semigroup} $\Rightarrow$ \begin{equation}
        \probabilitymeasure_\point[\process_\timepoint\in\measurableset]=\kernel(\point,\measurableset)\;\;\;\text{for all }(\point,\measurableset)\in\measurablespace\times\measurablespace\text{ and }\timepoint\ge0
    \end{equation} and hence \eqref{eq:markov-kernel-induced-by-realization-of-markov-semigroup} is $\measurablespace$-measurable for all $\event\in\bigcup_{\timepoint\ge0}\sigma(\process_\timepoint)$.}\footnote{i.e. \begin{equation}
            (\measurablespace,\eventsystem)\ni(\point,\event)\mapsto\probabilitymeasure_\point\left[\event\right]
        \end{equation} is a Markov kernel with source $(\measurablespace,\measurablesystem)$ and target $(\probabilityspace,\eventsystem)$.}
    \end{enumerate}
    \begin{flushright}
        $\square$
    \end{flushright}
\end{definition}

\section{Killing of processes}\label{sec:killing-by-a-multiplicative-functional}

In this section, we describe how to \emph{kill} a given Markov process according to a rule defined by a multiplicative functional, thereby limiting its lifetime. The resulting "subordinated" processes were also considered in a modified form in \citep[SectionIII.2 and SectionIII.3]{blumenthal1968markov}. We omit certain restrictions given there, which are not necessary for our purposes, and focus on the results of interest to us. These include verification of a Markov property and identification of the generator.

Let \begin{itemize}[$\circ$]
    \item $(\probabilityspace,\eventsystem,(\filtration_\timepoint)_{\timepoint\ge0},(\shift_\timepoint)_{\timepoint\ge0})$ be a stochastic base\footnote{see \autoref{def:stochastic-base}.};
    \item $(\measurablespace,\measurablesystem)$ be a measurable space with \begin{equation}
        \{\point\}\in\measurablesystem\;\;\;\text{for all }\point\in\measurablespace;
    \end{equation}
    \item $(\localsemigroup_\timepoint)_{\timepoint\ge0}$ be a Markov semigroup on $(\measurablespace,\measurablesystem)$;
    \item $((\localprocess_\timepoint)_{\timepoint\ge0},(\probabilitymeasure_\point)_{\point\in\measurablespace})$ be a realization\footnote{see \autoref{def:realization-of-markov-semigroup}} of $(\localsemigroup_\timepoint)_{\timepoint\ge0}$ on $(\probabilityspace,\eventsystem,(\filtration_\timepoint)_{\timepoint\ge0},(\shift_\timepoint)_{\timepoint\ge0})$;
    \item $(\multiplicativefunctional_\timepoint)_{\timepoint\ge0}$ be a multiplicative functional\footnote{see \autoref{def:multiplicative-functional}} on $\largebracket(\probabilityspace,\eventsystem,\largebracket(\filtration^\localprocess_\timepoint\largebracket)_{\timepoint\ge0},(\shift_\timepoint)_{\timepoint\ge0},(\probabilitymeasure_\point)_{\point\in\measurablespace}\largebracket)$.
\end{itemize}

In the following, we will restrict the \emph{lifetime} of the process $(\localprocess_\timepoint)_{\timepoint\ge0}$ according to the rule defined by the multiplicative functional $(\multiplicativefunctional_\timepoint)_{\timepoint\ge0}$.

We are interested in various special cases of this scenario. In the simplest and most conventional case, we restrict the lifetime of $(\localprocess_\timepoint)_{\timepoint\ge0}$ by directly specifying the lifetime:

\begin{setting}[killing at a terminal time]\label{set:killing-at-terminal-time}
    Assume \begin{equation}
        \multiplicativefunctional=1_{[0,\;\lifetime)}
    \end{equation} for some terminal time $\lifetime$ on $\largebracket(\probabilityspace,\eventsystem,\largebracket(\filtration^\localprocess_\timepoint\largebracket)_{\timepoint\ge0},(\shift_\timepoint)_{\timepoint\ge0},(\probabilitymeasure_\point)_{\point\in\measurablespace}\largebracket)$.
    \begin{flushright}
        $\square$
    \end{flushright}
\end{setting}

An important generic instance of this setting is given when the multiplicative functional is the one induced by an additive functional:

\begin{setting}\label{set:multiplicative-functional-induced-by-additive-functional}
    Assume \begin{equation}
        \multiplicativefunctional=e^{-\additivefunctional}
    \end{equation} for some additive functional $(\additivefunctional_{\timepoint})_{\timepoint\ge0}$ on $(\probabilityspace,\eventsystem,(\filtration_{\timepoint})_{\timepoint\ge0},(\shift_{\timepoint})_{\timepoint\ge0},(\probabilitymeasure_{\point})_{\point\in \measurablespace})$.\bigbreak
    \noindent\textit{Remark}: If $(\additivefunctional_\timepoint)_{\timepoint\ge0}$ is nondecreasing, then $(\multiplicativefunctional_\timepoint)_{\timepoint\ge0}$ is nonincreasing.
    \begin{flushright}
        $\square$
    \end{flushright}
\end{setting}

As it will turn out, the integrals induced by the process are those additive functionals that will conveniently enable the adjustment towards a given invariant measure:

\begin{setting}[killing at an exponential rate]\label{set:additive-functional-induced-by-integral}
    Assume \autoref{set:multiplicative-functional-induced-by-additive-functional}, $(\localprocess_{\timepoint})_{\timepoint\ge0}$ is $(\filtration_{\timepoint})_{\timepoint\ge0}$-progressive and \begin{equation}
        \additivefunctional_{\timepoint}=\int_0^{\timepoint}\killingrate(\localprocess_{\prevtimepoint})\dif \prevtimepoint\;\;\;\text{for all }\timepoint\ge0
    \end{equation} for some $\measurablesystem$-measurable $\killingrate:\measurablespace\to[0,\infty)$.\bigbreak
    \noindent\textit{Remark}: $(\additivefunctional_{\timepoint})_{\timepoint\ge0}$ is $(\filtration_{\timepoint})_{\timepoint\ge0}$-adapted and continuous $\Rightarrow$ $(\additivefunctional_{\timepoint})_{\timepoint\ge0}$ is $(\filtration_{\timepoint})_{\timepoint\ge0}$-progressive.
    \begin{flushright}
        $\square$
    \end{flushright}
\end{setting}

Orthogonal to the above settings, mild regularity properties are required for specific conclusions:

\begin{setting}\label{set:left-regular-local-process}
    Assume $\measurablesystem=\sigma(\topology)$ for some topology $\topology$ on $\measurablespace$ and $(\localprocess_\timepoint)_{\timepoint\ge0}$ is left-regular with respect to $\topology$
\end{setting}

By \eqref{eq:multiplicative-functional-inistial-state}, \begin{equation}\label{eq:multiplicative-functional-inistial-state-repeated}
    \multiplicativefunctional_0\in\{0,1\}\;\;\;\probabilitymeasure_\point\text{-almost surely}\text{ for all }\point\in\measurablespace.
\end{equation} The points $\point\in\measurablespace$ for which $\multiplicativefunctional_0$ starts $\probabilitymeasure_\point$-almost surely at $1$ play an important role:

\begin{definition}\label{def:set-em}
    Let \begin{equation}
        \measurablespace_\multiplicativefunctional:=\left\{\point\in\measurablespace:\probabilitymeasure_\point\left[\multiplicativefunctional_0=1\right]=1\right\}.
    \end{equation}
    \begin{flushright}
        $\square$
    \end{flushright}
\end{definition}

\begin{example}[killing at a random time (cont.)]
    Assume \autoref{set:killing-at-terminal-time} $\Rightarrow$ \begin{equation}
        \measurablespace_\multiplicativefunctional=\left\{\point\in\measurablespace:\probabilitymeasure_\point\left[\lifetime>0\right]=1\right\}.
    \end{equation}
    \begin{flushright}
        $\square$
    \end{flushright}
\end{example}

\begin{example}
    Assume \autoref{set:multiplicative-functional-induced-by-additive-functional} $\Rightarrow$ \begin{equation}\label{eq:set-em-1}
        \multiplicativefunctional_0=e^{-\overbrace{\additivefunctional_0}^{=\;0}}=1\;\;\;\probabilitymeasure_\point\text{-almost surely for all }\point\in\measurablespace
    \end{equation} and hence \begin{equation}
        \measurablespace_\multiplicativefunctional=\measurablespace.
    \end{equation}
    \begin{flushright}
        $\square$
    \end{flushright}
\end{example}

\begin{lemma}\label{lem:set-em}
    \leavevmode
    \begin{enumerate}[(i)]
        \item $\measurablespace_{\multiplicativefunctional}\in\measurablesystem$.
        \item\label{lem:set-em-ii}Let $\point\in\measurablespace$ $\Rightarrow$ \begin{equation}
            \probabilitymeasure_\point\left[\multiplicativefunctional_0=1\right]\in\{0,1\}.
        \end{equation}
        \item Let $\point\in\measurablespace\setminus \measurablespace_{\multiplicativefunctional}$ $\Rightarrow$ \begin{equation}\label{eq:vanishing-set-for-multiplicative-functional}
            \multiplicativefunctional_0=0\;\;\;\probabilitymeasure_\point\text{-almost surely}.
        \end{equation}
    \end{enumerate}
    \begin{proof}[Proof\textup:\nopunct]
        \leavevmode
        \begin{enumerate}[(i)]
            \item \begin{equation}
                1_{\measurablespace_{\multiplicativefunctional}}(\point)=\probabilitymeasure_\point\left[\multiplicativefunctional_0=1\right]\;\;\;\text{for all }\point\in\measurablespace
            \end{equation} and hence $1_{\measurablespace_{\multiplicativefunctional}}$ is $\measurablesystem$-measurable.
            \item \leavevmode\begin{itemize}[$\circ$]
                \item $\probabilitymeasure_\point\left[\localprocess_0=\point\right]=1$ $\Rightarrow$ $\sigma(\localprocess_0)$ is $\probabilitymeasure_\point$-trivial.
                \item $\multiplicativefunctional_0$ is $\sigma(\localprocess_0)$-measurable $\Rightarrow$ $\sigma(\multiplicativefunctional_0)\subseteq\sigma(\localprocess_0)$ is $\probabilitymeasure_\point$-trivial.
            \end{itemize}
            \item \leavevmode\begin{itemize}[$\circ$]
                \item\eqref{eq:multiplicative-functional-inistial-state-repeated} $\Rightarrow$ \begin{equation}
                    \multiplicativefunctional_0\in\{0,1\}\;\;\;\probabilitymeasure_\point\text{-almost surely}.
                \end{equation}
                \item\ref{lem:set-em-ii} and $\point\not\in\measurablespace_{\multiplicativefunctional}$ $\Rightarrow$ \begin{equation}
                    \probabilitymeasure_\point\left[\multiplicativefunctional_0=1\right]=0.
                \end{equation}
            \end{itemize}
        \end{enumerate}
    \end{proof}
\end{lemma}

We now can define the \emph{subordinate} semigroup\footnote{cf. \citep[Section~III.2]{blumenthal1968markov}}, which will turn out to be the transition semigroup of the subsequently constructed \emph{killed process}:

\begin{definition}[transition semigroup of killed process]\label{def:killed-local-semigroup}
    \begin{equation}
        (\killedlocalsemigroup_\timepoint\integrand)(\point):=\expectation_\point\left[\multiplicativefunctional_\timepoint\integrand(\localprocess_\timepoint)\right]\;\;\;\text{for }(\point,\integrand)\in \measurablespace\times\measurablesystem_b\text{ and }\timepoint\ge0
    \end{equation} is called \textbf{killing of} $\bm{(\localsemigroup_\timepoint)_{\timepoint\ge0}}$ \textbf{by} $\bm{(\multiplicativefunctional_\timepoint)_{\timepoint\ge0}}$. Assume \autoref{set:additive-functional-induced-by-integral} $\Rightarrow$ $(\killedlocalsemigroup_\timepoint)_{\timepoint\ge0}$ is called \textbf{killing of} $\bm{(\localsemigroup_\timepoint)_{\timepoint\ge0}}$ \textbf{at rate} $\bm\killingrate$.
    \begin{flushright}
        $\square$
    \end{flushright}
\end{definition}

\begin{example}[killing at a terminal time (cont.)]
    Assume \autoref{set:killing-at-terminal-time} $\Rightarrow$ \begin{equation}
        (Q_tf)(x)=\operatorname E_x\left[f(X_t);t<\tau\right]\;\;\;\text{for all }(x,f)\in E\times\mathcal E_b\text{ and }t\ge0.
    \end{equation}
    \begin{flushright}
        $\square$
    \end{flushright}
\end{example}

\begin{proposition}
    $(\killedlocalsemigroup_\timepoint)_{\timepoint\ge0}$
    is a sub-Markov semigroup on $(\measurablespace,\measurablesystem)$. If $(\point,\integrand)\in(\measurablespace_\multiplicativefunctional,\measurablesystem_b)$, then \begin{equation}
        (\killedlocalsemigroup_\timepoint\integrand)(\point)=(\localsemigroup_\timepoint\integrand)(\point)-\expectation_\point\left[(\multiplicativefunctional_0-\multiplicativefunctional_\timepoint)\integrand(\localprocess_\timepoint)\right].
    \end{equation}
    \begin{proof}[Proof\textup:\nopunct]
        \leavevmode
        \begin{itemize}[$\circ$]
            \item \begin{listclaim}
                Let $\timepoint\ge0$ $\Rightarrow$ $\killedlocalsemigroup_\timepoint$ is a sub-Markov kernel on $(\measurablespace,\measurablesystem)$.
                \begin{proof}[Proof\textup:\nopunct]
                    Trivial.
                \end{proof}
            \end{listclaim}
            \item \begin{listclaim}
                Let $\prevtimepoint,\timepoint\ge0$ and $(\point,\integrand)\in \measurablespace\times\measurablesystem_b$ $\Rightarrow$ \begin{equation}
                    (\killedlocalsemigroup_{\prevtimepoint+\timepoint}\integrand)(\point)=\left(\killedlocalsemigroup_{\prevtimepoint}(\killedlocalsemigroup_\timepoint\integrand)\right)(\point).
                \end{equation}
                \begin{proof}[Proof\textup:\nopunct]
                    \begin{equation}
                        \begin{split}
                            (\killedlocalsemigroup_{\prevtimepoint+\timepoint}\integrand)(\point)&\stackrel{\text{def}}=\expectation_\point\left[\multiplicativefunctional_{\prevtimepoint+\timepoint}\integrand(\localprocess_{\prevtimepoint+\timepoint})\right]\\&=\expectation_\point\left[\multiplicativefunctional_{\prevtimepoint}(\multiplicativefunctional_\timepoint\integrand(\localprocess_\timepoint)\circ\shift_{\prevtimepoint})\right]\\&=\expectation_\point\underbraceinsidebracket{[}{\multiplicativefunctional_{\prevtimepoint}}{\expectation_\point\left[\multiplicativefunctional_\timepoint\integrand(\localprocess_\timepoint)\circ\shift_{\prevtimepoint}\mid\filtration_\prevtimepoint\right]}{=\;\expectation_{\localprocess_{\prevtimepoint}}\left[\multiplicativefunctional_\timepoint\integrand(\localprocess_\timepoint)\right]\;\stackrel{\text{def}}=\;(\killedlocalsemigroup_\timepoint\integrand)(\localprocess_{\prevtimepoint})}{}{]}\\&\stackrel{\text{def}}=\left(\killedlocalsemigroup_{\prevtimepoint}(\killedlocalsemigroup_\timepoint\integrand)\right)(\point).
                        \end{split}
                    \end{equation}
                \end{proof}
            \end{listclaim}
        \end{itemize}
    \end{proof}
\end{proposition}

We will construct the killed process on a suitable augmentation of the given probability space:

\begin{definition}\label{def:killed-process-augmented-space}
    Let \begin{equation}
        \begin{split}
            \tilde\probabilityspace&:=\probabilityspace\times[0,\infty];\\
            \tilde\eventsystem&:=\eventsystem\otimes\mathcal B([0,\infty])
        \end{split}
    \end{equation} and $\probabilityprojection_\indexsetelement$ denote the projection from $\tilde\probabilityspace$ onto the $\indexsetelement$th coordinate $\Rightarrow$ \begin{equation}
        \begin{split}
            \tilde\localprocess&:=\localprocess\circ\probabilityprojection_1;\\
            \tilde\multiplicativefunctional&:=\multiplicativefunctional\circ\probabilityprojection_1
        \end{split}
    \end{equation} and \begin{equation}
        \tilde\lifetime:=\probabilityprojection_2.
    \end{equation}
    \begin{flushright}
        $\square$
    \end{flushright}
\end{definition}

Note that \begin{equation}
    \tilde\probabilityspace_\timepoint:=\probabilityspace\times(\timepoint,\infty]\in\tilde\eventsystem\;\;\;\text{for all }\timepoint\ge0.
\end{equation}

The information on the augmented system will be characterized by the following filtration:

\begin{proposition}
    \begin{equation}
        \begin{split}
            \tilde\filtration_\timepoint&:=\left\{\tilde\event\in\tilde\eventsystem:\exists\eventsystem\in\filtration_\timepoint:\tilde\event\cap\tilde\probabilityspace_\timepoint=\event\times(\timepoint,\infty]\right\}\\&=\left\{\tilde\event\in\tilde\eventsystem:\exists\event\in\filtration_\timepoint:\tilde\event\cap\{\;\timepoint<\tilde\lifetime\;\}=(\event\times[0,\infty])\cap\{\;\timepoint<\tilde\lifetime\;\}\right\}
        \end{split}
    \end{equation} for $\timepoint\ge0$ is a filtration on $\largebracket(\tilde\probabilityspace,\tilde\eventsystem\largebracket)$.
    \begin{proof}[Proof\textup:\nopunct]
        This is easily been checked.
    \end{proof}
\end{proposition}

\begin{proposition}\label{prop:lifetime-on-augmented-space-is-stopping-time}
    $\tilde\lifetime$ is a $\largebracket(\tilde\filtration_\timepoint\largebracket)_{\timepoint\ge0}$-stopping time.
    \begin{proof}[Proof\textup:\nopunct]
        Let $\timepoint\ge0$ $\Rightarrow$ \begin{equation}
            \largebracket\{\tilde\lifetime\le\timepoint\largebracket\}\cap\left\{\timepoint<\tilde\lifetime\right\}=\emptyset=(\emptyset\times[0,\infty])\cap\largebracket\{\timepoint<\tilde\lifetime\largebracket\}.
        \end{equation}
    \end{proof}
\end{proposition}

Let \begin{equation}
    \tilde\shift_\timepoint(\outcome,\lifetime):=(\shift_\timepoint(\outcome),\lifetime-\timepoint\wedge\lifetime)=(\shift_\timepoint(\outcome),(\lifetime-\timepoint)\vee0)\;\;\;\text{for }(\outcome,\lifetime)\in\tilde\probabilityspace.
\end{equation}

\begin{remark}
    $\largebracket(\tilde\probabilityspace,\tilde{\eventsystem},\largebracket(\tilde\filtration_\timepoint\largebracket)_{\timepoint\ge0},\largebracket(\tilde\shift_\timepoint\largebracket)_{\timepoint\ge0}\largebracket)$ is a stochastic base.
    \begin{flushright}
        $\square$
    \end{flushright}
\end{remark}

The stopping time $\tilde\lifetime$ trivially has the terminal time property:

\begin{proposition}\label{prop:lifetime-on-augmented-space-has-terminal-time-property}
    Let $\timepoint\ge0$ $\Rightarrow$ \begin{equation}
        \timepoint+\tilde\lifetime\circ\tilde\shift_\timepoint=\tilde\lifetime\;\;\;\text{on }\left\{\timepoint<\tilde\lifetime\right\}.
    \end{equation}
    \begin{proof}[Proof\textup:\nopunct]
        Let $\timepoint\ge0$ and $(\outcome,\lifetime)\in\left\{\timepoint<\tilde\lifetime\right\}$ $\Rightarrow$ \begin{equation}
            \timepoint<\lifetime
        \end{equation} and hence \begin{equation}
            \tilde\lifetime\left(\tilde\shift_\timepoint(\outcome,\lifetime)\right)\stackrel{\text{def}}=\lifetime-\timepoint\wedge\lifetime=\lifetime-\timepoint\stackrel{\text{def}}=\tilde\lifetime(\outcome,\lifetime)-\timepoint.
        \end{equation}
    \end{proof}
\end{proposition}

Let $\deadstate\not\in\measurablespace$ and $(\measurablespace^\ast,\measurablesystem^\ast)$ denote the one-point extension of $(\measurablespace,\measurablesystem)$ by $\deadstate$.

\begin{definition}[killed process]\label{def:killed-process}
    \begin{equation}
        \tilde\killedlocalprocess_\timepoint:=\left.\begin{cases}
            \tilde\localprocess_\timepoint&\text{, if }\timepoint<\tilde\lifetime;\\
            \deadstate&\text{, otherwise}
        \end{cases}\right\}\;\;\;\text{for }\timepoint\in[0,\infty]
    \end{equation} is called \textbf{killing of} $\bm{(\localprocess_\timepoint)_{\timepoint\ge0}}$ $\textbf{by}$ $\bm{(\multiplicativefunctional_\timepoint)_{\timepoint\ge0}}$. Assume \autoref{set:additive-functional-induced-by-integral} $\Rightarrow$ $\largebracket(\tilde\killedlocalprocess_\timepoint\largebracket)_{\timepoint\ge0}$ is called \textbf{killing of} $\bm{(\localprocess_\timepoint)_{\timepoint\ge0}}$ $\textbf{at rate}$ $\bm\killingrate$.
    \begin{flushright}
        $\square$
    \end{flushright}
\end{definition}

Intuitively, we set the killed process $\tilde\killedlocalprocess$ to the abstract \emph{dead state} $\deadstate$ after its \emph{lifetime} $\tilde\lifetime$ has elapsed. By construction, $\deadstate$ is a \emph{trap} for $(\killedlocalprocess_\timepoint)_{\timepoint\ge0}$: After $\deadstate$ has been reached, it is never left. Moreover, \begin{equation}
    \tilde\lifetime=\inf\left\{\timepoint\ge0:\tilde\killedlocalprocess_\timepoint=\deadstate\right\}
\end{equation} with the common convention $\inf\emptyset:=\infty$.

We easily observe that this process is adapted to the previously constructed filtration:

\begin{proposition}\label{prop:killed-process-is-adapted}
    $\largebracket(\tilde \killedlocalprocess_\timepoint\largebracket)_{\timepoint\ge0}$ is an $(\measurablespace^\ast,\measurablesystem^\ast)$-valued $\largebracket(\tilde\filtration_\timepoint\largebracket)_{\timepoint\ge0}$-adapted process on $\largebracket(\tilde\probabilityspace,\tilde\eventsystem\largebracket)$.
    \begin{proof}[Proof\textup:\nopunct]
        This is obvious by construction.
    \end{proof}
\end{proposition}

We now strengthen our assumptions on $(\multiplicativefunctional_\timepoint)_{\timepoint\ge0}$. First of all, we assume that $(\multiplicativefunctional_\timepoint)_{\timepoint\ge0}$ is nonincreasing.

\begin{remark}
    $(\multiplicativefunctional_\timepoint)_{\timepoint\ge0}$ is nonincreasing $\Rightarrow$ \begin{equation}
        0\le\multiplicativefunctional_\timepoint\le \multiplicativefunctional_0\stackrel{\eqref{eq:vanishing-set-for-multiplicative-functional}}=0\;\;\;\probabilitymeasure_\point\text{-almost surely}
    \end{equation} and hence \begin{equation}
        (\killedlocalsemigroup_\timepoint\integrand)(\point)\stackrel{\text{def}}=\expectation_\point\left[\multiplicativefunctional_\timepoint\integrand(\localprocess_\timepoint)\right]=0\;\;\;\text{for all }\integrand\in\measurablesystem_b
    \end{equation} for all $\point\in \measurablespace\setminus\measurablespace_\multiplicativefunctional$ and $\timepoint\ge0$.
    \begin{flushright}
        $\square$
    \end{flushright}
\end{remark}

Furthermore, we now assume that $(\multiplicativefunctional_\timepoint)_{\timepoint\ge0}$ is right-continuous. Let \begin{equation}
    \multiplicativefunctional_\infty:=0.
\end{equation}

\begin{proposition}[Lebesgue-Stieltjes measure induced by $M$]
    Let $\outcome\in\probabilityspace$ $\Rightarrow$ $\exists!\lebesguestieltjes_\outcome\in\mathcal M_{\le1}([0,\infty])$ with \begin{equation}
        \lebesguestieltjes_\outcome(\{0\})=0
    \end{equation} and \begin{equation}\label{eq:killing-with-mf-ls-identity}
        \lebesguestieltjes_\outcome\left((\prevtimepoint,\timepoint]\right)=\multiplicativefunctional_\prevtimepoint(\outcome)-\multiplicativefunctional_\timepoint(\outcome)\;\;\;\text{for all }0\le \prevtimepoint<\timepoint\le\infty.
    \end{equation} Moreover, \begin{equation}\label{eq:killing-with-mf-ls-identity-2}
        \lebesguestieltjes_\outcome([0,\infty])=\lebesguestieltjes_\outcome\left([0,\infty)\right)=\multiplicativefunctional_0(\outcome)
    \end{equation} and hence \begin{equation}\label{eq:lebesgue-stieltjes-measure-induced-by-multiplicative-functional-eq4}
        \lebesguestieltjes_\outcome\in\mathcal M_1([0,\infty])\Leftrightarrow \multiplicativefunctional_0(\outcome)=1.
    \end{equation}
    \textit{Remark}: $\lebesguestieltjes_\outcome$ is the Lebesgue-Stieltjes measure on $\mathcal B([0,\infty])$ induced by\newline $\left(\multiplicativefunctional_0(\outcome)-\multiplicativefunctional_\timepoint(\outcome)\right)_{\timepoint\in[0,\;\infty]}$.
    \begin{proof}[Proof\textup:\nopunct]
        \leavevmode
        \begin{itemize}[$\circ$]
            \item Let $$\secondintegrand(\timepoint):=\multiplicativefunctional_0(\outcome)-\multiplicativefunctional_\timepoint(\outcome)\in[0,1]\;\;\;\text{for }\timepoint\in[0,\infty].$$
            \item $\multiplicativefunctional(\outcome)$ is right-continuous $\Rightarrow$ $\secondintegrand$ is right-continuous.
            \item $\multiplicativefunctional(\outcome)$ is bounded and nonincreasing $\Rightarrow$\newline$\secondintegrand$ is bounded and nondecreasing $\Rightarrow$\newline$\secondintegrand$ has bounded variation $\Rightarrow$\newline There is a unique measure $\measure_\secondintegrand$ on $\mathcal B([0,\infty))$ with \begin{equation}
                \measure_\secondintegrand\left((\prevtimepoint,\timepoint]\right)=\secondintegrand(\timepoint)-\secondintegrand(\prevtimepoint)\ge0\;\;\;\text{for all }0\le \prevtimepoint\le \timepoint.
            \end{equation}
            \item \begin{equation}
                \measure_\secondintegrand(\{0\})=\secondintegrand(0)=0.
            \end{equation}
            \item $[0,\infty)=\{0\}\uplus\bigcup_{n\in\mathbb N}(0,n]$ and $\measure_\secondintegrand$ is continuous from below $\Rightarrow$ \begin{equation}
                \begin{split}
                    \measure_\secondintegrand\left([0,\infty)\right)&=\underbrace{\measure_\secondintegrand(\{0\})}_{=\;0}+\lim_{n\to\infty}\measure_\secondintegrand\left((0,n]\right)\\&=\lim_{n\to\infty}\secondintegrand(n)=\secondintegrand(\infty)=\multiplicativefunctional_0(\outcome)\le1
                \end{split}
            \end{equation} and hence $\measure_\secondintegrand\in\mathcal M_{\le1}\left([0,\infty)\right)$.
            \item Note that $\measure_\secondintegrand$ has a unique extension to a sub-probability measure $\tilde\measure_\secondintegrand$ on $[0,\infty]$ with \begin{equation}
                \begin{split}
                    \tilde\measure_\secondintegrand\left((\prevtimepoint,\infty]\right)&=\lim_{\timepoint\to\infty}\measure_\secondintegrand\left((\prevtimepoint,\timepoint]\right)=\lim_{\timepoint\to\infty}\left(\secondintegrand(\timepoint)-\secondintegrand(\prevtimepoint)\right)\\&=\secondintegrand(\infty)-\secondintegrand(\prevtimepoint)
                \end{split}
            \end{equation} for all $\prevtimepoint\ge0$.
            \item \begin{equation}
                \begin{split}
                    \tilde\measure_\secondintegrand([0,\infty])&=\underbrace{\measure_\secondintegrand(\{0\})}_{=\;0}+\tilde\measure_\secondintegrand\left((0,\infty]\right)\\&=\secondintegrand(\infty)-\underbrace{\secondintegrand(0)}_{=\;0}=\multiplicativefunctional_0(\outcome)
                \end{split}
            \end{equation} and hence $\tilde\measure_\secondintegrand$ is a probability measure iff \begin{equation}
                \multiplicativefunctional_0(\outcome)=1.
            \end{equation}
            \item Moreover, \begin{equation}
                \tilde\measure_\secondintegrand(\{\infty\})=\tilde\measure_\secondintegrand([0,\infty])-\measure_\secondintegrand\left([0,\infty)\right)=0.
            \end{equation}
        \end{itemize}
    \end{proof}
\end{proposition}

\begin{remark}
	\begin{equation}
		\lebesguestieltjes(\outcome,\measurableset):=\lebesguestieltjes_\outcome(\measurableset)\;\;\;\text{for }(\outcome,\measurableset)\in\probabilityspace\times\mathcal B([0,\infty])
	\end{equation} is a sub-Markov kernel with source $(\probabilityspace,\eventsystem)$ and target $([0,\infty],\mathcal B([0,\infty]))$.
    \begin{flushright}
        $\square$
    \end{flushright}
\end{remark}

\begin{example}\label{ex:lebesgue-stieltjes-measure-in-set-additive-functional-induced-by-integral}
    Assume \autoref{set:additive-functional-induced-by-integral} $\Rightarrow$ \begin{equation}
        \multiplicativefunctional_\prevtimepoint-\multiplicativefunctional_\timepoint=\int_\prevtimepoint^\timepoint\multiplicativefunctional_\prevprevtimepoint\killingrate\left(\localprocess_\prevprevtimepoint\right)\dif\prevprevtimepoint\;\;\;\text{for all }0\le\prevtimepoint\le\timepoint
    \end{equation} and hence \begin{equation}
        \lebesguestieltjes(\outcome,\timedomain)=\int_{\timedomain\cap[0,\;\infty)}\multiplicativefunctional_\timepoint(\outcome)\killingrate\left(\localprocess_\timepoint(\outcome)\right)\dif\timepoint\;\;\;\text{for all }(\outcome,\timedomain)\in\probabilityspace\times\mathcal B([0,\infty]).
    \end{equation}
    \begin{flushright}
        $\square$
    \end{flushright}
\end{example}

We can now construct the probability measure with respect to which we will verify the Markov property of the killed process:

\begin{definition}\label{def:killed-process-probability-measure}
	Let $\point\in\measurablespace$ $\Rightarrow$ \begin{equation}
		\tilde\probabilitymeasure_\point:=\probabilitymeasure_\point\lebesguestieltjes;
	\end{equation} i.e. \begin{equation}
		\begin{split}
			\tilde\expectation_\point\left[\tilde\secondintegrand\right]&=\expectation_\point\left[\lebesguestieltjes \tilde\secondintegrand\right]\stackrel{\text{def}}=\int\probabilitymeasure_\point\left[\dif\outcome\right]\int\lebesguestieltjes_\outcome(\dif\lifetime)\tilde\secondintegrand(\outcome,\lifetime)\\&=\int\probabilitymeasure_\point\left[\dif\outcome\right]\int\tilde\secondintegrand(\outcome,\timepoint)\dif\left(\multiplicativefunctional_0(\outcome)-\multiplicativefunctional_\lifetime(\outcome)\right)
		\end{split}
	\end{equation} for all $\tilde\secondintegrand\in\tilde\eventsystem_b$.
    \begin{flushright}
        $\square$
    \end{flushright}
\end{definition}

\begin{remark}
    By \eqref{eq:lebesgue-stieltjes-measure-induced-by-multiplicative-functional-eq4}, \begin{equation}
        \tilde\probabilitymeasure_\point\in\mathcal M_1(\tilde\probabilityspace,\tilde\eventsystem)
    \end{equation} for all $\point\in\measurablespace_\multiplicativefunctional$. In contrast, \begin{equation}
        \multiplicativefunctional_0=0\;\;\;\probabilitymeasure_\point\text{-almost surely}
    \end{equation} and hence \begin{equation}
        \lebesguestieltjes_\outcome=0\;\;\;\text{for }\probabilitymeasure_\point\text{-almost surely all }\outcome\in\probabilityspace
    \end{equation} for all $\point\in\measurablespace\setminus\measurablespace_\multiplicativefunctional$. Thus, $\tilde\probabilitymeasure_\point$ is the trivial measure on $(\tilde\probabilityspace,\tilde\eventsystem)$ for all $\point\in\measurablespace\setminus\measurablespace_\multiplicativefunctional$.
    \begin{flushright}
        $\square$
    \end{flushright}
\end{remark}

In \autoref{set:multiplicative-functional-induced-by-additive-functional} an equivalent, more natural, description of the killed process is available:

\begin{example}\label{ex:killedprocess-distribution-multiplicative-functional-induced-by-additive-functional}
    Assume \autoref{set:multiplicative-functional-induced-by-additive-functional}. Let \begin{itemize}[$\circ$]
        \item $\point\in\measurablespace$;
        \item $\expvariable$ be a real-valued $\operatorname{Exp}(1)$-distributed random variable on $(\probabilityspace,\eventsystem,\probabilitymeasure_\point)$;
        \item $\lifetime:=\inf\left\{\timepoint\ge0:\additivefunctional_\timepoint\ge\expvariable\right\}$ and \begin{equation}
            \killedlocalprocess_\timepoint:=\left.\begin{cases}
                \localprocess_\timepoint&\text{, if }\timepoint<\lifetime;\\
                \deadstate&\text{, otherwise}
        \end{cases}\right\}\;\;\;\text{for }\timepoint\in[0,\infty].
        \end{equation}
    \end{itemize} If $\expvariable$ and $\localprocess$ are $\probabilitymeasure_\point$-independent, then \begin{equation}
        \tilde\probabilitymeasure_\point\circ\left(\tilde\localprocess,\tilde\lifetime\right)^{-1}=\probabilitymeasure_\point\circ\;\left(\localprocess,\lifetime\right)^{-1}.
    \end{equation}
    \begin{proof}[Proof\textup:\nopunct]
        \leavevmode
        \begin{itemize}[$\circ$]
            \item $(\additivefunctional_\timepoint)_{\timepoint\ge0}$ is $\largebracket(\filtration^\localprocess_\timepoint\largebracket)_{\timepoint\ge0}$-adapted $\Rightarrow$ $\sigma(\additivefunctional)\subseteq\sigma(\localprocess)$ $\Rightarrow$ $\expvariable$ and $\additivefunctional$ are $\probabilitymeasure_\point$-independent $\Rightarrow$ \begin{equation}
                \probabilitymeasure_\point\left[\additivefunctional_\timepoint<\expvariable\mid\localprocess\right]=\left.\probabilitymeasure_\point\left[\MakeLowercase\additivefunctional<\expvariable\right]\right|_{\MakeLowercase\additivefunctional\;=\;\additivefunctional_\timepoint}=e^{-\additivefunctional_\timepoint}\stackrel{\text{def}}=\multiplicativefunctional_\timepoint
            \end{equation} for all $\timepoint\ge0$.
            \item Let $\timepoint\ge0$ $\Rightarrow$ \begin{equation}
                \timepoint<\lifetime\Leftrightarrow\additivefunctional_\timepoint<\expvariable
            \end{equation} and hence \begin{equation}
                \begin{array}{rcl}
                    \tilde\expectation_\point\left[\integrand\left(\tilde\localprocess\right);\timepoint<\tilde\lifetime\right]&=&\displaystyle\int\probabilitymeasure_\point\left[\dif\omega\right]\integrand\left(\localprocess(\outcome)\right)\lebesguestieltjes_\outcome\left((\timepoint,\infty]\right)\\
                    &=&\expectation_\point\left[\multiplicativefunctional_\timepoint\integrand(\localprocess)\right]\\
                    &=&\expectation_\point\left[\probabilitymeasure_\point\left[\timepoint<\lifetime\mid\localprocess\right]\integrand(\localprocess)\right]=\expectation_\point\left[\integrand
                    (\localprocess);\timepoint<\lifetime\right]
                \end{array}
            \end{equation} for all $\integrand\in\measurablesystem^{\otimes[0,\;\infty)}_b$.       
        \end{itemize}
    \end{proof}
\end{example}

\begin{example}\label{ex:killedprocess-distribution-additive-functional-induced-by-integral}
    Reconsider \autoref{ex:killedprocess-distribution-multiplicative-functional-induced-by-additive-functional} under the specialized assumption in \autoref{set:additive-functional-induced-by-integral}. If $\killingrate$ is constant, then \begin{equation}
        \lifetime=\frac\expvariable\killingrate
    \end{equation} and hence \begin{equation}
        \probabilitymeasure_\point\circ\;\lifetime^{-1}=\operatorname{Exp}(\killingrate).
    \end{equation} That is, our construction is reduced to the killing of $(\localprocess_\timepoint)_{\timepoint\ge0}$ after an exponentially distributed amount of time. This has been studied extensively in the literature; for example in \citep{avrachenkov2013restart}.
    \begin{flushright}
        $\square$
    \end{flushright}
\end{example}


\begin{proposition}
	\begin{equation}\label{eq:propability-measure-of-killing-with-mf-has-measurable-dependence}
		\measurablespace\to[0,1]\;,\;\;\;\point\mapsto\tilde\probabilitymeasure_\point\largebracket[\tilde\event\largebracket]
	\end{equation} is $\measurablesystem$-measurable for all $\tilde\event\in\tilde\eventsystem$.
	\begin{proof}[Proof\textup:\nopunct]
		\leavevmode
		\begin{itemize}[$\circ$]
			\item Let \begin{equation}
				\mathcal M:=\eventsystem\times\left\{(\prevtimepoint,\timepoint]:0\le\prevtimepoint<\timepoint\le\infty\right\}
			\end{equation} and \begin{equation}
				\mathcal H:=\left\{\tilde\event\in\tilde\eventsystem:\eqref{eq:propability-measure-of-killing-with-mf-has-measurable-dependence}\text{ is }\measurablesystem\text{-measurable}\right\}.
			\end{equation}
			\item \begin{listclaim}\label{prop:propability-measure-of-killing-with-mf-has-measurable-dependence-claim1}
				$\mathcal M\subseteq\mathcal H$.
				\begin{proof}[Proof\textup:\nopunct]
					\leavevmode
					\begin{itemize}[$\circ$]
						\item Let $\tilde\event\in\mathcal M$ $\Rightarrow$ \begin{equation}
							\tilde\event:=\event\times(\prevtimepoint,\timepoint]
						\end{equation} for some $\event\in\eventsystem$ and $0\le\prevtimepoint<\timepoint\le\infty$ $\Rightarrow$ \begin{equation}
							\tilde\event_\outcome=\left.\begin{cases}(\prevtimepoint,\timepoint]&\text{, if }\outcome\in\event;\\\emptyset&\text{, otherwise}\end{cases}\right\}\;\;\;\text{for all }\outcome\in\probabilityspace
						\end{equation} and hence \begin{equation}
    							\tilde\probabilitymeasure_\point\largebracket[\tilde\event\largebracket]=\int_\event\probabilitymeasure_\point\left[\dif\outcome\right]\lebesguestieltjes_\outcome\left((\prevtimepoint,\timepoint]\right)=\expectation_\point\left[\multiplicativefunctional_\prevtimepoint-\multiplicativefunctional_\timepoint; \event \right]\;\;\;\text{for all }\point\in\measurablespace_\multiplicativefunctional.
						\end{equation}
						\item \begin{equation}
							\measurablespace\to[0,1]\;,\;\;\;\point\mapsto\expectation_\point[\secondintegrand]
						\end{equation} is $\measurablesystem$-measurable for all $\secondintegrand\in\eventsystem_b$ $\Rightarrow$ Claim.
					\end{itemize}
				\end{proof}
			\end{listclaim}
			\item\autoref{prop:propability-measure-of-killing-with-mf-has-measurable-dependence-claim1} $\Rightarrow$ Claim.
		\end{itemize}
	\end{proof}
\end{proposition}

The embedding of objects, which live on the initial space $(\probabilitymeasure,\eventsystem)$, into the augmented space does not affect their distribution:

\begin{lemma}\label{lem:killed-process-space-dependence-on-first-argument}
	Let \begin{itemize}[$\circ$]
		\item $(\probabilityspace',\eventsystem')$ be a measurable space;
		\item $V:\probabilityspace\to\probabilityspace'$ be bounded and $(\eventsystem,\eventsystem')$-measurable;
		\item $\tilde V:=V\circ\probabilityprojection_1$;
		\item $\point\in\measurablespace_\multiplicativefunctional$.
	\end{itemize} Then, \begin{equation}\label{eq:killed-process-space-dependence-on-first-argument}
		\tilde\probabilitymeasure_\point\left[\tilde V\in\event'\right]=\probabilitymeasure_\point\left[V\in\event'\right]\;\;\;\text{for all }\event'\in\eventsystem'.
	\end{equation}
	\begin{proof}[Proof\textup:\nopunct]
		Let $\event'\in\eventsystem'$ and $\tilde\event:=\left\{\tilde V\in\event'\right\}$ $\Rightarrow$ \begin{equation}
		    \begin{split}
		        \tilde\event_\outcome&=\left\{\lifetime\in[0,\infty]:(\outcome,\lifetime)\in\tilde\event\right\}=\left\{\lifetime\in[0,\infty]:\tilde V(\outcome)\in\measurableset\right\}\\&=\begin{cases}
		            [0,\infty]&\text{, if }V(\outcome)\in\measurableset;\\\emptyset&\text{, otherwise}
		        \end{cases}
		    \end{split}
		\end{equation} for all $\outcome\in\probabilitymeasure$ $\Rightarrow$ \begin{equation}
		    \begin{split}
		        \tilde\probabilitymeasure_\point\left[\tilde V\in \measurableset\right]&\stackrel{\text{def}}=\int \probabilitymeasure_\point \left[\dif \outcome \right]\lebesguestieltjes_\outcome \left(\tilde \event_\outcome \right)\\&=\int_{\{\;V\;\in\;\event'\;\}}\probabilitymeasure_\point\left[\dif\outcome \right]\underbrace{\lebesguestieltjes_\outcome([0,\infty])}_{=\;\multiplicativefunctional_0(\outcome)}=\probabilitymeasure_\point\left[V\in\measurableset\right].
		    \end{split}
		\end{equation}
	\end{proof}
\end{lemma}

In particular, $\largebracket(\tilde\localprocess_\timepoint\largebracket)_{\timepoint\ge0}$ is still a well-posed process on the augmented space:

\begin{lemma}\label{prop:localprocess-on-augmented-space}
    $\largebracket(\tilde\localprocess_\timepoint\largebracket)_{\timepoint\ge0}$ is an $(\measurablespace,\measurablesystem)$-valued process on $\largebracket(\tilde\probabilityspace,\tilde\eventsystem,(\filtration_\timepoint\times[0,\infty])_{\timepoint\ge0},\largebracket(\tilde\shift_\timepoint\largebracket)_{\timepoint\ge0},\largebracket(\tilde\probabilitymeasure_\point\largebracket)_{\point\in\measurablespace_\multiplicativefunctional}\largebracket)$.
	\begin{proof}[Proof\textup:\nopunct]
        \leavevmode
        \begin{itemize}[$\circ$]
            \item Let $\point\in\measurablespace_\multiplicativefunctional$.
            \item\autoref{lem:killed-process-space-dependence-on-first-argument} $\Rightarrow$ \begin{equation}
                \tilde\probabilitymeasure_\point\left[\tilde\localprocess_0=x\right]=1.
            \end{equation}
            \item Let $\prevtimepoint,\timepoint\ge0$ $\Rightarrow$ \begin{equation}
                \localprocess_{\prevtimepoint+\timepoint}(\outcome)=(\localprocess_\timepoint\circ\shift_\prevtimepoint)(\outcome)\;\;\;\text{for all }\outcome\in\probabilityspace\setminus N
            \end{equation} for some $\probabilitymeasure_\point$-null set $N\in\eventsystem$.
            \item \autoref{lem:killed-process-space-dependence-on-first-argument} $\Rightarrow$ \begin{equation}
                \tilde N:=N\times[0,\infty]\in\tilde\eventsystem
            \end{equation} is a $\tilde\probabilitymeasure_\point$-null set.
            \item Let $(\outcome,\lifetime)\in\tilde\probabilityspace\setminus\tilde N=(\probabilityspace\setminus N)\times[0,\infty]$ $\Rightarrow$ \begin{equation}
                \left(\tilde\localprocess_\timepoint\circ\tilde\shift_\prevtimepoint\right)(\outcome,\lifetime)=(\localprocess_\timepoint\circ\shift_\prevtimepoint)(\outcome)=\localprocess_{\prevtimepoint+\timepoint}(\outcome)=\tilde\localprocess_{\prevtimepoint+\timepoint}(\outcome,\lifetime).
            \end{equation}
        \end{itemize}
	\end{proof}
\end{lemma}

We are also able to define the transition semigroup of the killed process in terms of objects on the augmented space only:

\begin{remark}
    \autoref{lem:killed-process-space-dependence-on-first-argument} $\Rightarrow$ \begin{equation}
        (\killedlocalsemigroup_\timepoint\integrand)(\point)=\tilde\expectation_\point\left[\tilde\multiplicativefunctional_\timepoint\integrand\left(\tilde\localprocess_\timepoint\right)\right]\;\;\;\text{for all }(\point,\integrand)\in\measurablespace\times\eventsystem_b\text{ and }\timepoint\ge0.
    \end{equation}
    \begin{flushright}
        $\square$
    \end{flushright}
\end{remark}

We now verify our intuition that $\largebracket(\tilde\killedlocalprocess_\timepoint\largebracket)_{\timepoint\ge0}$ is obtained by killing $\largebracket(\localprocess_\timepoint)_{\timepoint\ge0}$ at rate \begin{equation}
    -\frac{\dif\multiplicativefunctional_\timepoint}{\multiplicativefunctional_\timepoint}.
\end{equation} To do so, we let \begin{equation}
	\tilde\eventsystem_0:=\eventsystem\times[0,\infty],
\end{equation} which captures the information on the initial space.

\begin{lemma}
    Let $\point\in\measurablespace$ and $\timepoint\ge0$ $\Rightarrow$ \begin{equation}\label{eq:killed-process-lifetime-distribution}
        \tilde\probabilitymeasure_\point\left[\timepoint<\tilde\lifetime\expectationmid\tilde\eventsystem_0\right]=\tilde{\multiplicativefunctional}_\timepoint.
    \end{equation}
    \begin{proof}[Proof\textup:\nopunct]
		\leavevmode
		\begin{itemize}[$\circ$]
			\item Let $\tilde\event\in\tilde\eventsystem_0$ and \begin{equation}
				\tilde{\secondevent}:=\tilde\event\cap\{\timepoint<\tilde\lifetime\}.
			\end{equation}
			\item $\tilde\event\in\tilde\eventsystem_0$ $\Rightarrow$ $\tilde\event=\event\times[0,\infty]$ for some $\event\in\eventsystem$ $\Rightarrow$ \begin{equation}
				\tilde\event_\outcome=\begin{cases}[0,\infty]&\text{, if }\outcome\in\event;\\\emptyset&\text{, otherwise}\end{cases}
			\end{equation} and \begin{equation}
				\tilde{\secondevent}_\outcome=\begin{cases}(\timepoint,\infty]&\text{, if }\outcome\in\event;\\\emptyset&\text{, otherwise}\end{cases}
			\end{equation} for all $\outcome\in\probabilityspace$ $\Rightarrow$ \begin{equation}
			    \begin{array}{rcl}
                    \tilde\probabilitymeasure_\point\left[\timepoint<\tilde\lifetime;\tilde\event\right]&\stackrel{\text{def}}=&\tilde\probabilitymeasure_\point\left[\tilde{\secondevent}\right]\stackrel{\text{def}}=\displaystyle\int\probabilitymeasure_\point\left[\dif\outcome\right]\lebesguestieltjes_\outcome\left(\tilde{\secondevent}_\outcome\right)\\&=&\displaystyle\int_\event\probabilitymeasure_\point\left[\dif\outcome\right]\lebesguestieltjes_\outcome\left((\timepoint,\infty]\right)\\&\stackrel{\eqref{eq:killing-with-mf-ls-identity}}=&\displaystyle\int_\event\probabilitymeasure\left[\dif\outcome\right]\multiplicativefunctional_\timepoint(\outcome)\stackrel{\text{def}}=\expectation\left[\multiplicativefunctional_\timepoint;\event\right]\\&=&\displaystyle\int\probabilitymeasure_\point\left[\dif\outcome\right]\multiplicativefunctional_\timepoint(\outcome)\underbrace{\lebesguestieltjes_\outcome\left(\tilde\event_\outcome\right)}_{=\;\multiplicativefunctional_0(\outcome)}\\&=&\displaystyle\int\probabilitymeasure_\point\left[\dif\outcome\right]\int\lebesguestieltjes_\outcome(\dif\lifetime)\tilde{\multiplicativefunctional}_\timepoint(\outcome,\lifetime)1_{\tilde\event}(\outcome,\lifetime)\\&\stackrel{\text{def}}=&\tilde\expectation_\point\left[\tilde{\multiplicativefunctional}_\timepoint;\tilde\event\right].
			    \end{array}
			\end{equation}
		\end{itemize}
    \end{proof}
\end{lemma}

\begin{corollary}
    Let $(\point,\integrand)\in\measurablespace\times\measurablesystem_b$ and $\timepoint\ge0$ $\Rightarrow$ \begin{equation}
        (\killedlocalsemigroup_\timepoint\integrand)(\point)=\tilde\expectation_\point\left[\integrand\left(\tilde\localprocess_\timepoint\right);\timepoint<\tilde\lifetime\right].
    \end{equation}
    \begin{proof}[Proof\textup:\nopunct]
        $\tilde\localprocess_\timepoint$ is $\tilde\eventsystem_0$-measurable $\Rightarrow$ \begin{equation}
            \begin{split}
                (\killedlocalsemigroup_\timepoint\integrand)(\point)&=\tilde\expectation_\point\left[\tilde\multiplicativefunctional_\timepoint\integrand\left(\tilde\localprocess_\timepoint\right)\right]=\tilde\expectation_\point\left[\tilde\probabilitymeasure_\point\left[\timepoint<\tilde\lifetime\mid\tilde\eventsystem_0\right]\integrand\left(\tilde\localprocess_\timepoint\right)\right]\\&=\tilde\expectation_\point\left[\integrand\left(\tilde\localprocess_\timepoint\right);\timepoint<\tilde\lifetime\right].
            \end{split}
        \end{equation}
    \end{proof}
\end{corollary}

For subsequent usage in \autoref{sec:concatenation}, the following result on the joint distribution of the lifetime $\tilde\lifetime$ and the \emph{exit point} $\tilde\killedlocalprocess_{\tilde\lifetime-}$ is useful. However, it will not be needed in the remainder of this section.

\begin{lemma}\label{lem:killed-process-lifefime-joint-distribution}
    Assume \autoref{set:left-regular-local-process} $\Rightarrow$ \begin{enumerate}[(i)]
        \item $\largebracket(\tilde\killedlocalprocess_\timepoint\largebracket)_{\timepoint\ge0}$ is left-regular with respect to $\topology$ at $\tilde\lifetime$ on $\left\{\;\tilde\lifetime<\infty\;\right\}$ and \begin{equation}
            \tilde\killedlocalprocess_{\tilde\lifetime-}(\outcome,\lifetime)=\left.\begin{cases}
                \localprocess_{\lifetime-}(\outcome)&\text{, if }\timepoint<\lifetime;\\
                \deadstate&\text{, otherwise}
            \end{cases}\right\}\;\;\;\text{for all }(\outcome,\lifetime)\in\probabilityspace\times[0,\infty).
        \end{equation}
        \item Let $\point\in\measurablespace_\multiplicativefunctional$ and $(\timedomain,\measurableset)\in\mathcal B([0,\infty))\times\measurablesystem$ $\Rightarrow$ \begin{equation}
            \begin{split}
                \tilde\probabilitymeasure_\point\left[\largebracket(\tilde\lifetime,\tilde\killedlocalprocess_{\tilde\lifetime-}\largebracket)\in\timedomain\times\measurableset\right]&=\expectation_\point\left[\int_\timedomain1_\measurableset\left(\localprocess_{\lifetime-}\right)\dif\left(\multiplicativefunctional_0-\multiplicativefunctional_\lifetime\right)\right]\\
                &=\tilde\expectation_\point\left[\int_\timedomain1_\measurableset\largebracket(\tilde\killedlocalprocess_{\lifetime-}\largebracket)\dif\largebracket(\tilde\multiplicativefunctional_0-\tilde\multiplicativefunctional_\lifetime\largebracket)\right].
            \end{split}
        \end{equation}
    \end{enumerate}
    \begin{proof}[Proof\textup:\nopunct]
        \leavevmode
        \begin{enumerate}[(i)]
            \item \leavevmode\begin{itemize}[$\circ$]
                \item Let $\tilde\omega:=(\omega,\lifetime)\in\tilde\Omega$ $\Rightarrow$ \begin{equation}
                    \tilde\killedlocalprocess_t(\tilde\omega)\stackrel{\text{def}}=\left.\begin{cases}
                        \localprocess_\timepoint(\omega)&\text{, if }t<\lifetime;\\\Delta&\text{, otherwise}
                    \end{cases}\right\}\;\;\;\text{for all }t\in[0,\infty]
                \end{equation} and \begin{equation}
                    \tilde\lifetime(\tilde\omega)=\lifetime.
                \end{equation}
                \item If $\lifetime=0$, then \begin{equation}
                    \tilde\killedlocalprocess_t(\tilde\omega)=\Delta\;\;\;\text{for all }t\in[0,\infty]
                \end{equation} is trivially left-regular with respect to $\tau$ at $\tilde\lifetime(\tilde\omega)\stackrel{\text{def}}=\lifetime=0$.
                \item Assume $\lifetime\in(0,\infty)$.
                \item Let \begin{equation}
                    \point:=\localprocess(\omega).
                \end{equation}
                \item $x$ is left-regular with respect to $\tau$ at $\lifetime$ $\Rightarrow$ \begin{equation}
                    \forall N\in\mathcal N_\tau\left(x(\lifetime-)\right):\exists\delta>0:x\left((\lifetime-\delta,\delta)\right)\subseteq N.
                \end{equation}
                \item Let $\tau^\ast$ denote the topology on $E^\ast$ with $\mathcal E^\ast=\sigma(\tau^\ast)$ and $N^\ast\in\mathcal N_{\tau^\ast}\left(x(\lifetime-)\right)$ $\Rightarrow$ $\exists O^\ast\in\tau^\ast$ with \begin{equation}
                    x(\lifetime-)\in O^\ast\subseteq N^\ast.
                \end{equation}
                \item \autoref{prop:subspace-topology-in-one-point-extension-is-original-topology} $\Rightarrow$ \begin{equation}
                    \left.\tau^\ast\right|_E=\tau
                \end{equation} and hence \begin{equation}
                    N:=O^\ast\cap E\in\tau.
                \end{equation}
                \item $x(\lifetime-)\in E$ $\Rightarrow$ $N\in\mathcal N_\tau\left(x(\lifetime-)\right)$ $\Rightarrow$ $\exists\delta>0$ with \begin{equation}
                    \tilde\killedlocalprocess_{\left(\tilde\lifetime(\tilde\omega)-\delta,\:\tilde\lifetime(\tilde\omega)\right)}(\tilde\omega)=x\left((\lifetime-\delta,\lifetime)\right)\subseteq N\subseteq O^\ast\subseteq N^\ast.
                \end{equation}
            \end{itemize}
            \item \leavevmode\begin{itemize}[$\circ$]
                \item Let \begin{equation}
					\tilde A:=\left\{\tilde\lifetime\in I\right\}\cap\left\{\tilde\killedlocalprocess_{\tilde\lifetime-}\in B\right\}.
				\end{equation}
				\item Note that \begin{equation}
					\begin{split}
						\tilde A_\omega&=\left\{\lifetime\in I:\tilde\killedlocalprocess_{\tilde\lifetime-}(\omega,\lifetime)\in B\right\}\\&=\left\{\lifetime\in I\cap(0,\infty):X_{\lifetime-}(\omega)\in B\right\}
					\end{split}
				\end{equation} for all $\omega\in\Omega$ $\Rightarrow$ \begin{equation}
					\begin{array}{rcl}
						\tilde{\operatorname P}_x\left[\tilde A\right]&=&\displaystyle\int\operatorname P_x\left[\dif\omega\right]\int_{I\cap(0,\:\infty)}\alpha_\omega(\dif\lifetime)1_B\left(X_{\lifetime-}(\omega)\right)\\&\stackrel{\alpha_\omega\{0\}\:=\:0}=&\displaystyle\int_I\alpha_\omega(\dif\lifetime)1_B\left(X_{\lifetime-}(\omega)\right)\\&=&\displaystyle\int_I1_B\left(X_{\lifetime-}(\omega)\right)\dif\left(M_0(\omega)-M_\lifetime(\omega)\right)\\&=&\displaystyle\operatorname E_x\left[\int_I1_B\left(X_{\lifetime-}\right)\dif\left(M_0-M_\lifetime)\right)\right].
					\end{array}
				\end{equation}
            \end{itemize}
        \end{enumerate}
    \end{proof}
\end{lemma}

\begin{example}\label{ex:killed-process-lifefime-joint-distribution-in-set-additive-functional-induced-by-integral}
    Assume \autoref{set:additive-functional-induced-by-integral} and \autoref{set:left-regular-local-process} $\Rightarrow$ \begin{equation}
        \begin{split}
            \tilde\probabilitymeasure_\point\left[\largebracket(\tilde\lifetime,\tilde\killedlocalprocess_{\tilde\lifetime-}\largebracket)\in\timedomain\times\measurableset\right]&=\int_\timedomain\expectation_\point\left[\multiplicativefunctional_\timepoint\killingrate\left(\localprocess_\timepoint\right);\localprocess_{\timepoint-}\in\measurableset\right]\dif\timepoint\\
            &=\int_\timedomain\tilde\expectation_\point\left[\tilde\multiplicativefunctional_\timepoint\killingrate\largebracket(\tilde\killedlocalprocess_\timepoint\largebracket);\tilde\localprocess_{\timepoint-}\in\measurableset\right]\dif\timepoint
        \end{split}
    \end{equation} for all $\point\in\measurablespace$ and $(\timedomain,\measurableset)\in\mathcal B([0,\infty))\times\measurablesystem$.
    \begin{proof}[Proof\textup:\nopunct]
        \autoref{lem:killed-process-lifefime-joint-distribution} and \autoref{ex:lebesgue-stieltjes-measure-in-set-additive-functional-induced-by-integral} $\Rightarrow$ Claim.
    \end{proof}
\end{example}

$\largebracket(\tilde\killedlocalprocess_\timepoint\largebracket)_{\timepoint\ge0}$ takes values in the space $(\measurablespace^\ast,\measurablesystem^\ast)$. For that reason, we need to verify the Markov property on this larger space. We start by showing that it almost surely starts in the right initial point:

\begin{lemma}
    Let $\point\in\measurablespace^\ast$ $\Rightarrow$ \begin{equation}
        \left\{\tilde\killedlocalprocess_0=\point\right\}=\begin{cases}
            \left\{\localprocess_0=\point\right\}\times(0,\infty]&\text{, if }\point\in\measurablespace;\\\probabilityspace\times\{0\}&\text{, otherwise}.
        \end{cases}
    \end{equation}
    \begin{proof}[Proof\textup:\nopunct]
        \begin{equation}
            \tilde\killedlocalprocess(\outcome,\lifetime)=\begin{cases}
                \localprocess_0(\outcome)&\text{, if }0<\lifetime;\\
                \deadstate&\text{, otherwise}
            \end{cases}
        \end{equation} and hence \begin{equation}
            \begin{split}
                \left\{\tilde\killedlocalprocess_0=\point\right\}&=\underbraceinsidebracket{\{}{(\outcome,\lifetime)\in\probabilityspace\times\{0\}:}{\tilde\killedlocalprocess_0(\outcome,\lifetime)}{=\;\deadstate}{=\point}{\}}\;\uplus\\&\hphantom=\;\;\underbraceinsidebracket{\{}{(\outcome,\lifetime)\in\probabilityspace\times(0,\infty]:}{\tilde\killedlocalprocess_0(\outcome,\lifetime)}{=\;\localprocess_0(\outcome)}{=\point}{\}}\\
                &=\left.\begin{cases}
                    \emptyset&\text{, if }\point\in\measurablespace;\\
                    \probabilityspace\times\{0\}&\text{, otherwise}
                \end{cases}\right\}\;\uplus\\&\hphantom=\;\;\left.\begin{cases}
                    \left\{\localprocess_0=\point\right\}\times(0,\infty]&\text{, if }\point\in\measurablespace;\\
                    \emptyset&\text{, otherwise}
                \end{cases}\right\}
            \end{split}
        \end{equation} for all $(\omega,\lifetime)\in\tilde\probabilityspace$.
    \end{proof}
\end{lemma}

At this point, we have defined $\tilde\probabilitymeasure_\point$ only for those $\point$ which belong to $\measurablespace$. Since $\deadstate$ is a trap for $\largebracket(\tilde\killedlocalprocess_\timepoint\largebracket)_{\timepoint\ge0}$, it is unimportant how we choose $\tilde\probabilitymeasure_\point$ for $\point=\deadstate$. Formally, we let \begin{itemize}[$\circ$]
    \item $\measure\in\mathcal M_1(\probabilityspace,\eventsystem)$;
    \item $\delta$ denote the Dirac kernel on $([0,\infty],\mathcal B([0,\infty])$;
    \item $\tilde\probabilitymeasure_\deadstate:=\measure\otimes\delta_0$.
\end{itemize}



\begin{lemma}\label{lem:lifetime-in-deadstate}
    $\tilde\probabilitymeasure_\deadstate\left[\tilde\lifetime>0\right]=0$.
	\begin{proof}[Proof\textup:\nopunct]
        \begin{equation}
                \{\;\timepoint<\tilde\lifetime\;\}=\probabilityspace\times(\timepoint,\infty]
            \end{equation} and hence \begin{equation}
                \tilde\probabilitymeasure_\deadstate\left[\timepoint<\tilde\lifetime\right]=\underbrace{\measure(\probabilityspace)}_{=\;1}\underbrace{\delta_0\left((\timepoint,\infty]\right)}_{=\;0}=0.
        \end{equation} for all $\timepoint\ge0$.
    \end{proof}
\end{lemma}

Given all these definitions, we see that the killed process is a normal process on the augmented space:

\begin{lemma}\label{lem:killed-process-on-stochastic-base}
    $\largebracket(\tilde\killedlocalprocess_\timepoint\largebracket)_{\timepoint\ge0}$ is a normal $(\measurablespace^\ast,\measurablesystem^\ast)$-valued process on $\largebracket(\tilde\probabilityspace,\tilde{\eventsystem},\largebracket(\tilde\filtration_\timepoint\largebracket)_{\timepoint\ge0},\largebracket(\tilde\shift_\timepoint\largebracket)_{\timepoint\ge0},\largebracket(\tilde\probabilitymeasure_\point\largebracket)_{\point\in\measurablespace_\multiplicativefunctional\uplus\{\deadstate\}}\largebracket)$.
	\begin{proof}[Proof\textup:\nopunct]
        \leavevmode
        \begin{itemize}[$\circ$]
            \item\autoref{prop:killed-process-is-adapted} $\Rightarrow$ $\largebracket(\tilde \killedlocalprocess_\timepoint\largebracket)_{\timepoint\ge0}$ is an $(\measurablespace^\ast,\measurablesystem^\ast)$-valued $\largebracket(\tilde\filtration_\timepoint\largebracket)_{\timepoint\ge0}$-adapted process on $\largebracket(\tilde\probabilityspace,\tilde\eventsystem\largebracket)$.
            \item If $\point\in\measurablespace_\multiplicativefunctional\uplus\{\deadstate\}$, then \begin{equation}\label{eq:killed-process-initial-state}
                \tilde\probabilitymeasure_\point\left[\tilde\killedlocalprocess_0=\point\right]=\left.\begin{cases}
                    \expectation_\point\left[\multiplicativefunctional_0;\localprocess_0=\point\right]&\text{, if }\point\in\measurablespace;\\
                    \measure(\probabilityspace)\delta_0(\{0\})&\text{, otherwise}
                \end{cases}\right\}=1.
            \end{equation}
            \item If $\point=\deadstate$, then \begin{equation}
                \tilde\killedlocalprocess_\timepoint\circ\tilde\shift_\prevtimepoint=\deadstate=\tilde\killedlocalprocess_{\prevtimepoint+\timepoint}\;\;\;\tilde\probabilitymeasure_\point\text{-almost surely}
            \end{equation} by \autoref{lem:lifetime-in-deadstate}.
            \item Let $\point\in\measurablespace_\multiplicativefunctional$ and $\tilde N$ denote the set from the proof of \autoref{prop:localprocess-on-augmented-space} $\Rightarrow$ \begin{equation}
                \begin{split}
                    \left(\tilde\killedlocalprocess_\timepoint\circ\tilde\shift_\prevtimepoint\right)(\outcome,\lifetime)&\stackrel{\text{def}}=\begin{cases}
                        \underbrace{(\localprocess_\timepoint\circ\shift_\prevtimepoint)(\outcome)}_{=\;\localprocess_{\prevtimepoint+\timepoint}(\outcome)}&\text{, if }\timepoint<\lifetime-\prevtimepoint\wedge\lifetime\\\deadstate&\text{, otherwise}
                    \end{cases}\\&=\begin{cases}
                        \localprocess_{\prevtimepoint+\timepoint}(\outcome)&\text{, if }\prevtimepoint+\timepoint<\lifetime;\\\deadstate&\text{, otherwise}
                    \end{cases}\\&\stackrel{\text{def}}=\tilde\killedlocalprocess_{\prevtimepoint+\timepoint}(\outcome,\lifetime)
                \end{split}
            \end{equation} for all $(\outcome,\lifetime)\in\tilde\probabilityspace\setminus\tilde N\times[0,\infty]$.
        \end{itemize}
	\end{proof}
\end{lemma}

We are finally able to state and verify the Markov property of the killed process:

\begin{theorem}[killed process is Markov]\label{thm:killed-process-is-markov}
    Let $\prevtimepoint,\timepoint\ge0$ and $(\point,\integrand)\in(\measurablespace_\multiplicativefunctional\uplus\{\deadstate\})\times\measurablesystem_b$ $\Rightarrow$ \begin{equation}\label{eq:killed-process-is-markov-1}
        \tilde\expectation_\point\left[\integrand\left(\tilde\killedlocalprocess_{\prevtimepoint+\timepoint}\right)\expectationmid\tilde\filtration_\prevtimepoint\right]=(\killedlocalsemigroup_\timepoint\integrand)\left(\tilde\killedlocalprocess_\prevtimepoint\right).
    \end{equation}
    \begin{proof}[Proof\textup:\nopunct]
        \leavevmode
        \begin{itemize}[$\circ$]
            \item If $\point=\deadstate$, then \begin{equation}
                \begin{array}{rcl}
                    \tilde\expectation_\point\left[\integrand\left(\tilde\killedlocalprocess_{\prevtimepoint+\timepoint}\right)\expectationmid\tilde\filtration_\prevtimepoint\right]&=&\tilde\expectation_\deadstate\left[\integrand\left(\tilde\localprocess_{\prevtimepoint+\timepoint}\right);\prevtimepoint+\timepoint<\tilde\lifetime\expectationmid\tilde\filtration_\prevtimepoint\right]\\
                    &=&0\\
                    &=&\tilde\expectation_\deadstate\left[(\killedlocalsemigroup_\timepoint\integrand)\left(\tilde\localprocess_\prevtimepoint\right);\prevtimepoint<\tilde\lifetime\expectationmid\tilde\filtration_\prevtimepoint\right]\\
                    &=&\tilde\expectation_\point\left[(\killedlocalsemigroup_\timepoint\integrand)\left(\tilde\killedlocalprocess_\prevtimepoint\right)\expectationmid\tilde\filtration_\prevtimepoint\right]
                \end{array}
            \end{equation}
            \item Assume $\point\in\measurablespace_\multiplicativefunctional$ $\Rightarrow$ \begin{equation}
                \multiplicativefunctional_{\prevtimepoint+\timepoint}=\multiplicativefunctional_\prevtimepoint(\multiplicativefunctional_\timepoint\circ\shift_\prevtimepoint)\;\;\;\probabilitymeasure_\point\text{-almost surely}.
            \end{equation}
            \item $\multiplicativefunctional_\timepoint$ is $\sigma(\localprocess)$-measurable $\Rightarrow$ \begin{equation}\label{eq:memoryless-assumption}
                \begin{array}{rcl}
                    \expectation_\point\left[\multiplicativefunctional_{\prevtimepoint+\timepoint}\integrand(\localprocess_{\prevtimepoint+\timepoint})\expectationmid\filtration_\prevtimepoint\right]&=&\multiplicativefunctional_\prevtimepoint\expectation_\point\left[\multiplicativefunctional_\timepoint\integrand(\localprocess_\timepoint)\circ\shift_\prevtimepoint\expectationmid\filtration_\prevtimepoint\right]\\
                    &=&\multiplicativefunctional_\prevtimepoint\expectation_{\localprocess_\prevtimepoint}\left[\multiplicativefunctional_\timepoint\integrand(\localprocess_\timepoint)\right]\stackrel{\text{def}}=\multiplicativefunctional_\prevtimepoint(\killedlocalsemigroup_\timepoint\integrand)(\localprocess_\prevtimepoint).
                \end{array}
            \end{equation}
            \item Let $\tilde\filtrationset\in\tilde\filtration_\prevtimepoint$ $\Rightarrow$ $\exists\filtrationset\in\filtration_\prevtimepoint$ with \begin{equation}\label{eq:killed-process-is-markov-2}
                \tilde\filtrationset\cap\left\{\prevtimepoint<\tilde\lifetime\right\}=(\underbrace{\filtrationset\times[0,\infty]}_{\in\;\tilde\eventsystem_0})\cap\left\{\prevtimepoint<\tilde\lifetime\right\}.
            \end{equation}
            \item $\tilde\localprocess$ is $\tilde\eventsystem_0$-measurable $\Rightarrow$ \begin{equation}
                \begin{array}{rcl}
                    \tilde\expectation_\point\left[\integrand\left(\tilde\killedlocalprocess_{\prevtimepoint+\timepoint}\right);\tilde\filtrationset\right]&=&\tilde\expectation_\point\left[\integrand\left(\tilde\localprocess_{\prevtimepoint+\timepoint}\right);\tilde\filtrationset\cap\left\{\prevtimepoint+\timepoint<\tilde\lifetime\right\}\right]\\
                    &\stackrel{\eqref{eq:killed-process-is-markov-2}}=&\tilde\expectation_\point\left[\integrand\left(\tilde\localprocess_{\prevtimepoint+\timepoint}\right);(\filtrationset\times[0,\infty])\cap\left\{\prevtimepoint+\timepoint<\tilde\lifetime\right\}\right]\\
                    &=&\tilde\expectation_\point\underbraceinsidebracket{[}{}{\tilde\probabilitymeasure_\point\left[\prevtimepoint+\timepoint<\tilde\lifetime\expectationmid\tilde\eventsystem_0\right]}{\stackrel{\eqref{eq:killed-process-lifetime-distribution}}=\;\tilde\multiplicativefunctional_{\prevtimepoint+\timepoint}}{\integrand\left(\tilde\localprocess_{\prevtimepoint+\timepoint}\right);\filtrationset\times[0,\infty]}{]}\\
                    &\stackrel{\text{def}}=&\tilde\expectation_\point\left[\multiplicativefunctional_{\prevtimepoint+\timepoint}\integrand(\localprocess_{\prevtimepoint+\timepoint})1_\filtrationset\circ\probabilityprojection_1\right]\\
                    &\stackrel{\eqref{eq:killed-process-space-dependence-on-first-argument}}=&\expectation_\point\left[\multiplicativefunctional_{\prevtimepoint+\timepoint}\integrand(\localprocess_{\prevtimepoint+\timepoint});\filtrationset\right]\\
                    &\stackrel{\eqref{eq:memoryless-assumption}}=&\expectation_\point\left[\multiplicativefunctional_\prevtimepoint(\killedlocalsemigroup_\timepoint\integrand)(\localprocess_\prevtimepoint);\filtrationset\right]\\
                    &\stackrel{\eqref{eq:killed-process-space-dependence-on-first-argument}}=&\tilde\expectation_\point\left[\multiplicativefunctional_\prevtimepoint(\killedlocalsemigroup_\timepoint\integrand)(\localprocess_\prevtimepoint)1_\filtrationset\circ\probabilityprojection_1\right]\\
                    &\stackrel{\text{def}}=&\tilde\expectation_\point\left[\tilde\multiplicativefunctional_\prevtimepoint(\killedlocalsemigroup_\timepoint\integrand)(\tilde\localprocess_\prevtimepoint);\filtrationset\times[0,\infty]\right]\\
                    &=&\tilde\expectation_\point\underbraceinsidebracket{[}{}{\tilde\probabilitymeasure_\point\left[\prevtimepoint<\tilde\lifetime\expectationmid\tilde\eventsystem_0\right]}{\stackrel{\eqref{eq:killed-process-lifetime-distribution}}=\;\tilde\multiplicativefunctional_\timepoint}{(\killedlocalsemigroup_\timepoint\integrand)\left(\tilde\localprocess_\prevtimepoint\right);\filtrationset\times[0,\infty]}{]}\\
                    &=&\tilde\expectation_\point\left[(\killedlocalsemigroup_\timepoint\integrand)\left(\tilde\localprocess_\prevtimepoint\right);(\filtrationset\times[0,\infty])\cap\left\{\prevtimepoint<\tilde\lifetime\right\}\right]\\
                    &\stackrel{\eqref{eq:killed-process-is-markov-2}}=&\tilde\expectation_\point\left[(\killedlocalsemigroup_\timepoint\integrand)\left(\tilde\localprocess_\prevtimepoint\right);\tilde\filtrationset\cap\left\{\prevtimepoint<\tilde\lifetime\right\}\right]\\
                    &=&\tilde\expectation_\point\left[(\killedlocalsemigroup_\timepoint\integrand)\left(\tilde\killedlocalprocess_\prevtimepoint\right);\tilde\filtrationset\right].
                \end{array}
            \end{equation}
        \end{itemize}
    \end{proof}
\end{theorem}

Let $(\killedlocalsemigroup^\ast_\timepoint)_{\timepoint\ge0}$ denote the one-point extension of $(\killedlocalsemigroup_\timepoint)_{\timepoint\ge0}$ by $\deadstate$. For simplicity, assume now that $\measurablespace_\multiplicativefunctional=\measurablespace$, although the following conclusion can be adjusted for the general case.

\begin{corollary}[realization of killed semigroup]
    $\largebracket(\largebracket(\tilde\killedlocalprocess_\timepoint\largebracket)_{\timepoint\ge0},\largebracket(\tilde\probabilitymeasure_{\point^\ast}\largebracket)_{\point^\ast\in\measurablespace^\ast}\largebracket)$ is a realization of $(\killedlocalsemigroup^\ast_\timepoint)_{\timepoint\ge0}$ on $\largebracket(\tilde\probabilityspace,\tilde{\eventsystem},\largebracket(\tilde\filtration_\timepoint\largebracket)_{\timepoint\ge0},\largebracket(\tilde\shift_\timepoint\largebracket)_{\timepoint\ge0}\largebracket)$.
    \begin{proof}[Proof\textup:\nopunct]
        \leavevmode
        \begin{itemize}[$\circ$]
            \item\autoref{prop:killed-process-is-adapted} $\Rightarrow$ $\largebracket(\tilde \killedlocalprocess_\timepoint\largebracket)_{\timepoint\ge0}$ is an $(\measurablespace^\ast,\measurablesystem^\ast)$-valued $\largebracket(\tilde\filtration_\timepoint\largebracket)_{\timepoint\ge0}$-adapted process on $\largebracket(\tilde\probabilityspace,\tilde\eventsystem\largebracket)$.
            \item\autoref{thm:killed-process-is-markov} $\Rightarrow$ \eqref{eq:killed-process-is-markov-1} for all $(\point,\integrand)\in\measurablespace^\ast\times\measurablesystem_b$ $\Rightarrow$ \begin{equation}
                \tilde\expectation_\point\left[\integrand\left(\tilde \localprocess_{\prevtimepoint+\timepoint}\right)\expectationmid\tilde\filtration_\prevtimepoint\right]=(\killedlocalsemigroup^\ast_\timepoint\integrand)\left(\tilde \localprocess_\prevtimepoint\right)
            \end{equation} for all $(\point,\integrand)\in\measurablespace^\ast\times\measurablesystem^\ast_b$ by \autoref{rem:markov-property-extension}.
            \item \autoref{rem:kernel-induced-by-probability-measure} $\Rightarrow$ \begin{equation}
                \measurablespace^\ast\to[0,1]\;,\;\;\;\point^\ast\mapsto\tilde\probabilitymeasure_{\point^\ast}\left[\tilde\event\right]
            \end{equation} is $\measurablesystem^\ast$-measurable for all $\tilde\event\in\tilde\eventsystem$.
        \end{itemize}
    \end{proof}
\end{corollary}

We close this section by the identification of the generator of the killed process. To do so, let $\localgenerator$ denote the pointwise generator of $(\localsemigroup_\timepoint)_{\timepoint\ge0}$. Assume \begin{equation}\label{eq:multiplicative-functional-derivative-condition}
    \expectation_\point\left[\frac{\multiplicativefunctional_0-\multiplicativefunctional_\timepoint}{\timepoint}\integrand(\localprocess_\timepoint)\right]\xrightarrow{\timepoint\to0+}\killingrate(\point)\integrand(\point)\;\;\;\text{for all }(\point,\integrand)\in\measurablespace_\multiplicativefunctional\times\mathcal D(\localgenerator)
\end{equation} for some $\left.\measurablesystem\right|_{\measurablespace_\multiplicativefunctional}$-measurable $\killingrate:\measurablespace_\multiplicativefunctional\to[0,\infty)$.

\begin{example}\label{ex:multiplicative-functional-derivative-condition-1}
    If \begin{enumerate}[(i)]
    	\item\label{asm:multiplicative-functional-derivative-condition-1-i}$\multiplicativefunctional$ is (right-)differentiable at $0$ $\probabilitymeasure_\point$-almost surely and \begin{equation}\label{eq:domination-for-derivative-of-mf}
    		\begin{split}
    		    &\exists\timepoint_0\in(0,\infty]:\exists\dominatingvariable\in\mathcal L^1(\probabilitymeasure_\point):\\&\frac{\multiplicativefunctional_0-\multiplicativefunctional_\timepoint}\timepoint\le \dominatingvariable\;\;\;\text{for all }\timepoint\in(0,\timepoint_0)\;\probabilitymeasure_\point\text{-almost surely}
    		\end{split}
    	\end{equation} for all $\point\in\measurablespace_\multiplicativefunctional$; and
    	\item\label{asm:multiplicative-functional-derivative-condition-1-ii}$\integrand\circ\localprocess$ is (right-)continuous at $0$ $\probabilitymeasure_\point$-almost surely for all $(\point,\integrand)\in\measurablespace_\multiplicativefunctional\times\mathcal D\left(\localgenerator\right)$,
    \end{enumerate} then \eqref{eq:multiplicative-functional-derivative-condition} is satisfied with \begin{equation}
        \killingrate(\point)=-\expectation_\point\left[\left.\frac\dif{\dif\timepoint}\multiplicativefunctional_\timepoint\right|_{\timepoint=0+}\right]\;\;\;\text{for all }\point\in\measurablespace_\multiplicativefunctional.
    \end{equation}
    \begin{proof}[Proof\textup:\nopunct]
        \leavevmode
        \begin{itemize}[$\circ$]
            \item\ref{asm:multiplicative-functional-derivative-condition-1-i}$\Rightarrow$ \begin{equation}
                \frac{\multiplicativefunctional_0-\multiplicativefunctional_\timepoint}{\timepoint}\xrightarrow{\timepoint\to0+}-\left.\frac\dif{\dif\timepoint}\multiplicativefunctional_\timepoint\right|_{\timepoint=0+}\;\;\;\probabilitymeasure_\point\text{-almost surely}.
            \end{equation}
            \item\ref{asm:multiplicative-functional-derivative-condition-1-ii}$\Rightarrow$ \begin{equation}
                \integrand(\localprocess_\timepoint)\xrightarrow{\timepoint\to0+}\integrand(\localprocess_0)=\integrand(\point)\;\;\;\probabilitymeasure_\point\text{-almost surely}.
            \end{equation}
            \item \eqref{eq:domination-for-derivative-of-mf} and $\integrand$ is bounded $\Rightarrow$ Claim by Lebesgue's dominated convergence theorem.
        \end{itemize}
    \end{proof}
\end{example}

\begin{example}[killing at an exponential rate (cont.)]\label{ex:mf-derivative-for-af-induced-by-integral}
    Assume \autoref{set:additive-functional-induced-by-integral}. If $\killingrate$ is bounded and $\killingrate\circ\localprocess$ is (right-)continuous at $0$ $\probabilitymeasure_\point$-almost surely for all $\point\in\measurablespace_\multiplicativefunctional$, then \eqref{eq:domination-for-derivative-of-mf} is satisfied for $\timepoint_0=\infty$ and \begin{equation}
        \killingrate(\point)=-\left.\frac\dif{\dif\timepoint}\multiplicativefunctional_\timepoint(\outcome)\right|_{\timepoint=0+}\;\;\;\text{for }\probabilitymeasure_\point\text{-almost all }\outcome\in\probabilityspace\text{ for all }\point\in\measurablespace.
    \end{equation} That is, the definition of $\killingrate$ in \eqref{eq:multiplicative-functional-derivative-condition} is consistent with the usage of this symbol in \autoref{ex:mf-derivative-for-af-induced-by-integral}.
    \begin{flushright}
        $\square$
    \end{flushright}
\end{example}

\begin{example}\label{ex:multiplicative-functional-derivative-condition-2}
    If \begin{enumerate}[(i)]
        \item\label{ex:multiplicative-functional-derivative-condition-2-i} \begin{equation}
            [0,\infty)\ni\timepoint\mapsto\expectation_\point\left[\multiplicativefunctional_\timepoint\right]
        \end{equation} is (right-)differentiable at $0$; and
        \item\label{ex:multiplicative-functional-derivative-condition-2-ii} $\multiplicativefunctional$ and $\localprocess$ are $\probabilitymeasure_\point$-independent for all $\point\in\measurablespace_\multiplicativefunctional$,
    \end{enumerate} then \eqref{eq:multiplicative-functional-derivative-condition} is satisfied with \begin{equation}
        \killingrate(\point)=-\left.\frac\dif{\dif\timepoint}\expectation_\point\left[\multiplicativefunctional_\timepoint\right]\right|_{\timepoint=0+}\;\;\;\text{for all }\point\in\measurablespace_\multiplicativefunctional.
    \end{equation}
    \begin{flushright}
        $\square$
    \end{flushright}
\end{example}

Let $\killedlocalgenerator$ denote the pointwise generator of $(\killedlocalsemigroup_\timepoint)_{\timepoint\ge0}$.

\begin{theorem}[pointwise generator of killed process]\label{thm:pointwise-generator-of-process-killed-by-mf}
	\begin{equation}
		\frac{(\killedlocalsemigroup_\timepoint\integrand)(\point)-\integrand(\point)}{\timepoint}\xrightarrow{\timepoint\to0+}\left(\localgenerator\integrand\right)(\point)-\killingrate(\point)\integrand(\point)
	\end{equation} for all $(\point,\integrand)\in\measurablespace_\multiplicativefunctional\times\mathcal D\left(\localgenerator\right)$. In particular, if $\measurablespace=\measurablespace_\multiplicativefunctional$, then $\mathcal D\left(\localgenerator\right)\subseteq\mathcal D\left(\killedlocalgenerator\right)$ and \begin{equation}
		\killedlocalgenerator\integrand=\left(\localgenerator-\killingrate\right)\integrand\;\;\;\text{for all }\integrand\in\mathcal D\left(\localgenerator\right).
	\end{equation}
	\begin{proof}[Proof\textup:\nopunct]
		\leavevmode
		\begin{itemize}[$\circ$]
			\item Let $(\point,\integrand)\in\measurablespace_\multiplicativefunctional\times\measurablesystem_b$ $\Rightarrow$ \begin{equation}
				\begin{split}
					\frac{(\killedlocalsemigroup_\timepoint\integrand)(\point)-\integrand(\point)}\timepoint&=\frac{\expectation_\point\left[\multiplicativefunctional_\timepoint\integrand(\localprocess_\timepoint)-\integrand(\localprocess_0)\right]}\timepoint\\&=\frac{(\localsemigroup_\timepoint\integrand)(\point)-\integrand(\point)}\timepoint-\expectation_\point\left[\frac{\multiplicativefunctional_0-\multiplicativefunctional_\timepoint}\timepoint\integrand(\localprocess_\timepoint)\right]
				\end{split}
			\end{equation} for all $\timepoint>0$.
			\item\eqref{eq:multiplicative-functional-derivative-condition} $\Rightarrow$ Claim.
		\end{itemize}
	\end{proof}
\end{theorem}

\section{Concatenation of processes}\label{sec:concatenation}

We will consider the concatenation of (countably many) processes. Each process is executed until a given time, after which it will be terminated and replaced by the subsequent process. In order for the concatenation to make sense, each of the processes must have a finite lifetime. Otherwise execution of a process with infinite lifetime would never terminate and no subsequent process could be executed. We keep things as general as possible.

Let \begin{itemize}[$\circ$]
    \item $\indexset:=[1,n]\cap\mathbb N$ for some $n\in\mathbb N\uplus\{\infty\}$ and $\indexset_0:=\indexset_0$;
    \item $(\measurablespace_\indexsetelement,\measurablesystem_\indexsetelement)$ be a measurable space with \begin{equation}
        \{\point\}\in\measurablesystem_\indexsetelement\;\;\;\text{for all }\point\in\measurablespace_\indexsetelement;
    \end{equation}
    \item $\largebracket(\killedlocalsemigroup^{(\indexsetelement)}_\timepoint\largebracket)_{\timepoint\ge0}$ be a sub-Markov semigroup on $(\measurablespace_\indexsetelement,\measurablesystem_\indexsetelement)$;
    \item $\deadstate_\indexsetelement\in\measurablespace_\indexsetelement$ and $(\measurablespace_\indexsetelement^\ast,\measurablesystem_\indexsetelement^\ast)$ denote the one-point extension of $(\measurablespace_\indexsetelement,\measurablesystem_\indexsetelement)$ by $\deadstate_\indexsetelement$;
    \item $\largebracket(\killedlocalsemigroup^{(\indexsetelement)\ast}_\timepoint\largebracket)_{\timepoint\ge0}$ denote the one-point extension of $\largebracket(\localsemigroup^{(\indexsetelement)}_\timepoint\largebracket)_{\timepoint\ge0}$ by $\deadstate_\indexsetelement$;
    \item $\largebracket(\probabilityspace_\indexsetelement,\eventsystem_\indexsetelement,\largebracket(\filtration^{(\indexsetelement)}_\timepoint\largebracket)_{\timepoint\ge0},\largebracket(\shift^{(\indexsetelement)}_\timepoint\largebracket)_{\timepoint\ge0}\largebracket)$ be a stochastic base\footnote{see \autoref{def:stochastic-base}.};
    \item $\largebracket(\largebracket(\killedlocalprocess^{(\indexsetelement)}_\timepoint\largebracket)_{\timepoint\ge0},\largebracket(\probabilitymeasure^{(\indexsetelement)}_\point\largebracket)_{\point\in\measurablespace_\indexsetelement^\ast}\largebracket)$ be a realization\footnote{see \autoref{def:realization-of-markov-semigroup}.} of $\largebracket(\killedlocalsemigroup^{(\indexsetelement)\ast}_\timepoint\largebracket)_{\timepoint\ge0}$ on $\largebracket(\probabilityspace_\indexsetelement,\eventsystem_\indexsetelement,\largebracket(\filtration^{(\indexsetelement)}_\timepoint\largebracket)_{\timepoint\ge0},\largebracket(\shift^{(\indexsetelement)}_\timepoint\largebracket)_{\timepoint\ge0}\largebracket)$ and \begin{equation}
        \lifetime_\indexsetelement:=\inf\left\{\timepoint\ge0:\killedlocalprocess^{(\indexsetelement)}_\timepoint=\deadstate_\indexsetelement\right\}.
    \end{equation}
\end{itemize}


In \autoref{sec:killing-by-a-multiplicative-functional} we have shown one possibility how a Markov process with finite lifetime can be constructed. We will make use of this concrete construction for specific conclusions:

\begin{setting}[killing at an exponential rate]\label{set:concatenation-of-killing-at-exponential-rate}
    Let \begin{itemize}[$\circ$]
        \item $\largebracket(\localsemigroup^{(\indexsetelement)}_\timepoint\largebracket)_{\timepoint\ge0}$ be a Markov semigroup on $(\measurablespace_\indexsetelement,\measurablesystem_\indexsetelement)$;
        \item $\largebracket(\localprocess^{(\indexsetelement)}_\timepoint\largebracket)_{\timepoint\ge0}$ be an $(\measurablespace,\measurablesystem)$-valued process on $(\probabilityspace_\indexsetelement,\eventsystem_\indexsetelement)$;
        \item $\killingrate:\measurablespace_\indexsetelement\to[0,\infty)$ be $\measurablesystem_\indexsetelement$-measurable
    \end{itemize} for $\indexsetelement\in\indexset_0$. Assume $\largebracket(\largebracket(\localprocess^{(\indexsetelement)}_\timepoint\largebracket)_{\timepoint\ge0},\largebracket(\probabilitymeasure^{(\indexsetelement)}_\point\largebracket)_{\point\in\measurablespace_\indexsetelement}\largebracket)$ is a realization of $\largebracket(\localsemigroup^{(\indexsetelement)}_\timepoint\largebracket)_{\timepoint\ge0}$ on $\largebracket(\probabilityspace_\indexsetelement,\eventsystem_\indexsetelement,\largebracket(\filtration^{(\indexsetelement)}_\timepoint\largebracket)_{\timepoint\ge0},\largebracket(\shift^{(\indexsetelement)}_\timepoint\largebracket)_{\timepoint\ge0}\largebracket)$ and $\largebracket(\killedlocalsemigroup^{(\indexsetelement)}_\timepoint\largebracket)_{\timepoint\ge0}$ is the killing\footnote{see \autoref{def:killed-process}} of $\largebracket(\localprocess^{(\indexsetelement)}_\timepoint\largebracket)_{\timepoint\ge0}$ at rate $\killingrate$ for all $\indexsetelement\in\indexset_0$.
    \begin{flushright}
        $\square$
    \end{flushright}
\end{setting}

Once again, orthogonal to the above setting, mild regularity properties are required for specific conclusions:

\begin{setting}\label{set:left-regular-killedlocalprocess}
    Assume $\measurablesystem_\indexsetelement=\sigma(\topology_\indexsetelement)$ for some topology $\topology_\indexsetelement$ on $\measurablespace_\indexsetelement$ and $\largebracket(\killedlocalprocess^{(\indexsetelement)}_\timepoint\largebracket)_{\timepoint\ge0}$ is left-regular with respect to the one-point extension $\topology_\indexsetelement^\ast$ of $\topology_\indexsetelement$ by $\deadstate_\indexsetelement$ for all $\indexsetelement\in\indexset_0$.
    \begin{flushright}
        $\square$
    \end{flushright}
\end{setting}

The concatenated process will live on an augmented space:

\begin{definition}
    \begin{equation}
        \begin{split}
             \probabilityspace&:=\bigtimes_{\indexsetelement\in\indexset}\probabilityspace_\indexsetelement;\\
             \eventsystem&:=\bigvee_{\indexsetelement\in\indexset}\eventsystem_\indexsetelement;
        \end{split}
    \end{equation} and $\probabilityprojection_\indexsetelement$ denote the projection from $\tilde\probabilityspace$ onto the $\indexsetelement$th coordinate $\Rightarrow$ \begin{equation}
        \begin{split}
            \tilde\killedlocalprocess^{(\indexsetelement)}&:=\killedlocalprocess^{(\indexsetelement)}\circ\probabilityprojection_\indexsetelement;\\
            \tilde\lifetime_\indexsetelement&:=\lifetime_\indexsetelement\circ\probabilityprojection_\indexsetelement.
        \end{split}
    \end{equation}
    \begin{flushright}
        $\square$
    \end{flushright}
\end{definition}

Each object defined on $(\probabilityspace_\indexsetelement,\eventsystem_\indexsetelement)$ is naturally defined on $(\probabilityspace,\eventsystem)$ as well by composing it with $\probabilityprojection_\indexsetelement$. In the same spirit, the information available on $(\probabilityspace_\indexsetelement,\eventsystem_\indexsetelement)$ is naturally embedded into $(\probabilityspace,\eventsystem)$:

\begin{definition}
    Let $\timepoint\ge0$ $\Rightarrow$ \begin{equation}
        \begin{split}
            \tilde{\filtration}^{(\ge\indexsetelement)}_\timepoint&:=\filtration^{(\indexsetelement)}_\timepoint\times\bigtimes_{\substack{\secondindexsetelement\in\indexset\\\secondindexsetelement>\indexsetelement}}\probabilityspace_\secondindexsetelement;\\\tilde{\filtration}^{(\indexsetelement)}_\timepoint&:=\bigtimes_{\secondindexsetelement=1}^{\indexsetelement-1}\probabilityspace_\secondindexsetelement\times\tilde{\filtration}^{(\ge\indexsetelement)}_\timepoint
        \end{split}
    \end{equation} for $\timepoint\ge0$.
    \begin{flushright}
        $\square$
    \end{flushright}
\end{definition}

\begin{lemma}\label{lem:concatenated-process-measurability}
    Let $\indexsetelement\in\indexset$ $\Rightarrow$ \begin{equation}
        \left.\tilde\killedlocalprocess^{(\indexsetelement)}\right|_{\probabilityspace\times[0,\;\timepoint]}^{-1}(\measurablesystem)=\bigtimes_{\secondindexsetelement=1}^{\indexsetelement-1}\probabilityspace_\secondindexsetelement\times\left.\killedlocalprocess^{(\indexsetelement)}\right|_{\probabilityspace\times[0,\;\timepoint]}^{-1}(\measurablesystem_\indexsetelement)\times\bigtimes_{\substack{\secondindexsetelement\in\indexset\\\secondindexsetelement>\indexsetelement}}\probabilityspace_\secondindexsetelement\;\;\;\text{for all }\timepoint\ge0.
    \end{equation} In particular, if $\left(\killedlocalprocess^{(\indexsetelement)}_\timepoint\right)_{\timepoint\ge0}$ is $\left(\filtration^{(\indexsetelement)}_\timepoint\right)_{\timepoint\ge0}$-adaptive/progressive, then $\left(\tilde\killedlocalprocess^{(\indexsetelement)}_\timepoint\right)_{\timepoint\ge0}$ is $\left(\tilde\filtration^{(\indexsetelement)}_\timepoint\right)_{\timepoint\ge0}$-adaptive/progressive.
    \begin{proof}[Proof\textup:\nopunct]
        \leavevmode
        \begin{itemize}[$\circ$]
            \item Let $\timepoint\ge0$ $\Rightarrow$ \begin{equation}
                \begin{array}{rcl}
                    \left.\tilde\killedlocalprocess^{(\indexsetelement)}\right|_{\probabilityspace\times[0,\;\timepoint]}^{-1}(\measurablesystem)&\stackrel{\text{\autoref{prop:measurability-in-subspace}-\ref{prop:measurability-in-subspace-ii}}}=&\left.\tilde\killedlocalprocess^{(\indexsetelement)}\right|_{\probabilityspace\times[0,\;\timepoint]}^{-1}(\measurablesystem_\indexsetelement)\\
                    &\stackrel{\text{\autoref{lem:measurability-on-larger-space}}}=&\bigtimes_{\secondindexsetelement=1}^{\indexsetelement-1}\probabilityspace_\secondindexsetelement\times\left.\killedlocalprocess^{(\indexsetelement)}\right|_{\probabilityspace\times[0,\;\timepoint]}^{-1}(\measurablesystem_\indexsetelement)\times\bigtimes_{\substack{\secondindexsetelement\in\indexset\\\secondindexsetelement>\indexsetelement}}\probabilityspace_\secondindexsetelement.
                \end{array}
            \end{equation}
            \item If $\left.\killedlocalprocess^{(\indexsetelement)}\right|_{\probabilityspace_\indexsetelement\times[0,\;\timepoint]}$ is $\largebracket(\filtration^{(\indexsetelement)}_\timepoint\otimes\mathcal B([0,\timepoint]),\measurablesystem_\indexsetelement\largebracket)$-measurable, then \begin{equation}
                \bigtimes_{\secondindexsetelement=1}^{\indexsetelement-1}\probabilityspace_\secondindexsetelement\times\left.\killedlocalprocess^{(\indexsetelement)}\right|_{\probabilityspace\times[0,\;\timepoint]}^{-1}(\measurablesystem_\indexsetelement)\times\bigtimes_{\substack{\secondindexsetelement\in\indexset\\\secondindexsetelement>\indexsetelement}}\probabilityspace_\secondindexsetelement\subseteq\tilde\filtration^{(\indexsetelement)}_\timepoint\otimes\mathcal B([0,\timepoint]).
            \end{equation}
        \end{itemize}
    \end{proof}
\end{lemma}

We define the time after which the lifetimes of the first $\thirdindexsetelement\in\indexsetelement$ processes have elapsed:

\begin{definition}
    Let $\thirdindexsetelement\in\indexset\cup\left\{\sup\indexset\right\}$ $\Rightarrow$ \begin{equation}
        \lifetimesum_{\indexsetelement,\;\thirdindexsetelement}:=\sum_{\secondindexsetelement=\indexsetelement}^\thirdindexsetelement\lifetime_\secondindexsetelement\left(\outcome_\secondindexsetelement\right)\;\;\;\text{for }\outcome\in\bigtimes_{\secondindexsetelement=\indexsetelement}^\thirdindexsetelement\probabilityspace_\secondindexsetelement
    \end{equation} for all $\indexsetelement\in\indexset$ with $\indexsetelement\le\thirdindexsetelement$ and \begin{equation}
        \begin{array}{rcl}
            \lifetimesum_\thirdindexsetelement&:=&\lifetimesum_{1,\;\thirdindexsetelement};\\
            \tilde\lifetimesum_\thirdindexsetelement&:=&\displaystyle\sum_{\secondindexsetelement=1}^\thirdindexsetelement\tilde\lifetime_\secondindexsetelement.
        \end{array}      
    \end{equation}
    \begin{flushright}
        $\square$
    \end{flushright}
\end{definition}

Since our goal is to define a probability measure on $(\probabilityspace,\eventsystem)$ with respect to which the concatenated process is Markov, we need to fix the information with respect to which the Markov property has to be understood. To simplify the subsequent definition, we write\footnote{We slightly abuse notation to simplify the definition and the subsequent usage of it. If $\indexsetelement=1$, then \begin{equation}
    \event_{(\outcome_1,\;\ldots\;,\;\outcome_{\indexsetelement-1})}=\event.
\end{equation}} \begin{equation}
    \event_{(\outcome_1,\;\ldots\;,\;\outcome_{\indexsetelement-1})}:=\left\{\left(\outcome_\secondindexsetelement\right)_{\substack{\secondindexsetelement\in\indexset\\\secondindexsetelement\ge\indexsetelement}}\in\bigtimes_{\substack{\secondindexsetelement\in\indexset\\\secondindexsetelement\ge\indexsetelement}}\probabilityspace_\secondindexsetelement:(\outcome_\secondindexsetelement)_{\secondindexsetelement\in\indexset}\in\event\right\}
\end{equation} for $\event\subseteq\probabilityspace$ and $(\outcome_1,\ldots,\outcome_{\indexsetelement-1})\in\bigtimes_{\secondindexsetelement=1}^{\indexsetelement-1}\probabilityspace_\secondindexsetelement$.

\begin{remark}
    \begin{equation}
        \filtration_\timepoint:=\left\{\event\in\eventsystem\;\middle|\;\begin{aligned}&\left(\event\cap\left\{\timepoint<\tilde\lifetimesum_\indexsetelement\right\}\right)_{(\outcome_1,\;\ldots\;,\;\outcome_{\indexsetelement-1})}\in\tilde\filtration^{(\ge\indexsetelement)}_{\timepoint-\lifetimesum_{\indexsetelement-1}((\outcome_1,\;\ldots\;,\;\outcome_{\indexsetelement-1}))}\\&\text{for all }\outcome\in\left\{\lifetimesum_{\indexsetelement-1}\le\timepoint\right\}\text{ and }\indexsetelement\in\indexset\end{aligned}\right\}
    \end{equation} for $\timepoint\ge0$ is a filtration on $(\probabilityspace,\eventsystem)$.
    \begin{flushright}
        $\square$
    \end{flushright}
\end{remark}

Assume \begin{enumerate}[(i)]
    \item $\largebracket(\killedlocalprocess^{(\indexsetelement)}_\timepoint\largebracket)_{\timepoint\ge0}$ is $\largebracket(\filtration^{(\indexsetelement)}_\timepoint\largebracket)_{\timepoint\ge0}$-progressive for all $\indexsetelement\in\indexset\setminus\{1\}$;
    \item $\lifetime_\indexsetelement$ is an $\largebracket(\filtration^{(\indexsetelement)}_\timepoint\largebracket)_{\timepoint\ge0}$-stopping time for all $\indexsetelement\in\indexset$.
\end{enumerate}

\begin{definition}
    Let \begin{equation}
        \begin{split}
            \measurablespace&:=\bigcup_{\indexsetelement\in\indexsetelement}\measurablespace_\indexsetelement;\\\measurablesystem&:=\bigvee_{\indexsetelement\in\indexset}\measurablesystem_\indexsetelement;
        \end{split}
    \end{equation}
    \begin{flushright}
        $\square$
    \end{flushright}
\end{definition}

Let $\deadstate\not\in\measurablespace$ and $(\measurablespace^\ast,\measurablesystem^\ast)$ denote the one-point extension of $(\measurablespace,\measurablesystem)$ by $\deadstate$.

\begin{definition}[concatenated process]
    Let \begin{equation}
        \concatenatedprocess_\timepoint:=\left.\begin{cases}
            \tilde\killedlocalprocess^{(\indexsetelement}_{\timepoint-\tilde\lifetimesum_{\indexsetelement-1}}&\text{, if }\indexsetelement\in\indexset\text{ and }\timepoint\in\left[\tilde\lifetimesum_{\indexsetelement-1},\tilde\lifetimesum_\thirdindexsetelement\right)\\
            \deadstate&\text{, otherwise}
        \end{cases}\right\}\;\;\;\text{for }\timepoint\in[0,\infty].
    \end{equation}
\end{definition}

By construction, each process $\killedlocalprocess^{(\indexsetelement)}$ is executed until just immediately before the clock has reached the lifetime $\lifetime_\indexsetelement$, after which it will be terminated and replaced by the subsequent process $\killedlocalprocess^{(\indexsetelement+1)}$.

\begin{proposition}\label{prop:concatenated-process-is-adapted}
    $(\concatenatedprocess_\timepoint)_{\timepoint\ge0}$ is an $(\measurablespace^\ast,\measurablesystem^\ast)$-valued $(\filtration_\timepoint)_{\timepoint\ge0}$-adapted process on $(\probabilityspace,\eventsystem)$.
    \begin{proof}[Proof\textup:\nopunct]
        \leavevmode
        \begin{itemize}[$\circ$]
            \item Let $\timepoint\ge0$.
            \item \begin{listclaim}\label{prop:concatenated-process-is-adapted-claim1}
                Let $\measurableset\in\measurablesystem$ $\Rightarrow$ \begin{equation}
                    \underbrace{\left\{\concatenatedprocess\in\measurableset\right\}}_{=:\;\event}\in\filtration_\timepoint.
                \end{equation}
                \begin{proof}[Proof\textup:\nopunct]
                    \leavevmode
                    \begin{itemize}[$\circ$]
                        \item \autoref{lem:stopping-time-on-larger-space} and \autoref{lem:concatenated-process-measurability} $\Rightarrow$ \begin{equation}
                            A\cap\{t<\tilde\sigma_1\}=\underbrace{\left\{\tilde \killedlocalprocess^{(1)}_t\in B\right\}}_{\in\:\tilde{\mathcal F}^{(1)}_t}\cap\underbrace{\{t<\tilde\sigma_1\}}_{\in\:\tilde{\mathcal F}^{(1)}_t}\in\tilde{\mathcal F}^{(1)}_t.
                        \end{equation}
                        \item  Let $i\ge2$ and $(\omega_1,\ldots,\omega_{i-1})\in\{\sigma_{i-1}\le t\}$ $\Rightarrow$ \begin{equation}
                            \begin{split}
                                &(A\cap\{t<\tilde\sigma_i\})_{(\omega_1,\:\ldots\:,\:\omega_{i-1})}\\&\;\;\;\;=\left\{\left(\omega_j\right)_{\substack{j\in I\\j\ge i}}\in\bigtimes_{\substack{j\in I\\j\ge i}}\Omega_j:\killedlocalprocess^{(i)}_{t-\sigma_{i-1}(\omega_1,\:\ldots\:,\:\omega_{i-1})}(\omega_i)\in B\right\}\\&\;\;\;\;\;\;\;\;\;\cap\{t<\tilde\sigma_i\}_{(\omega_1,\:\ldots\:,\:\omega_{i-1})}\\&\;\;\;\;=\underbrace{\left(\underbrace{\left\{\omega_i\in\Omega_i:\killedlocalprocess^{(i)}_{t-\sigma_{i-1}(\omega_1,\:\ldots\:,\:\omega_{i-1})}(\omega_i)\in B\right\}}_{\in\:\mathcal F^{(i)}_{t-\sigma_{i-1}(\omega_1,\:\ldots\:,\:\omega_{i-1})}}\times\bigtimes_{\substack{j\in I\\j>i}}\Omega_j\right)}_{\in\:\tilde{\mathcal F}^{(\ge i)}_{t-\sigma_{i-1}(\omega_1,\:\ldots\:,\:\omega_{i-1})}}\\&\;\;\;\;\;\;\;\;\;\cap\underbrace{\{t<\tilde\sigma_i\}_{(\omega_1,\:\ldots\:,\:\omega_{i-1})}}_{\in\:\tilde{\mathcal F}^{(\ge i)}_{t-\sigma_{i-1}(\omega_1,\:\ldots\:,\:\omega_{i-1})}}\in\tilde{\mathcal F}^{(\ge i)}_{t-\sigma_{i-1}(\omega_1,\:\ldots\:,\:\omega_{i-1})}.
                            \end{split}
                        \end{equation}
                    \end{itemize}
                \end{proof}
            \end{listclaim}
            \item \begin{listclaim}\label{prop:concatenated-process-is-adapted-claim2}
                \begin{equation}
                    \underbrace{\left\{\concatenatedprocess=\deadstate\right\}}_{=:\;\event}\in\filtration_\timepoint.
                \end{equation}
                \begin{proof}[Proof\textup:\nopunct]
                    \leavevmode
                    \begin{itemize}[$\circ$]
                        \item \begin{equation}
                            A=\left\{t\ge\tilde\sigma_{\sup\indexset}\right\}
                        \end{equation} and hence \begin{equation}
                            \tilde\tau_i(\omega)\le\tilde\sigma_i(\omega)\le\tilde\sigma_{\sup\indexset}(\omega)\le t\;\;\;\text{for all }\omega\in A\text{ and }i\in\mathbb N.
                        \end{equation}
                        \item Thus, \begin{equation}
                            A\cap\{t<\tilde\tau_1\}=\emptyset\in\tilde{\mathcal F}^{(1)}_t
                        \end{equation} and \begin{equation}
                            (A\cap\{t<\tilde\sigma_i\})_{(\omega_1,\:\ldots\:,\:\omega_{i-1})}\in\tilde{\mathcal F}^{(\ge i)}_{t-\sigma_{i-1}(\omega_1,\:\ldots\:,\:\omega_{i-1})}
                        \end{equation} for all $i\ge2$ and $(\omega_1,\ldots,\omega_{i-1})\in\{\sigma_{i-1}\le t\}$.
                    \end{itemize}
                \end{proof}
            \end{listclaim}
            \item \autoref{prop:concatenated-process-is-adapted-claim1} and \autoref{prop:concatenated-process-is-adapted-claim2} $\Rightarrow$ $\concatenatedprocess_\timepoint$ is $(\filtration_\timepoint,\measurablesystem)$-measurable by \autoref{prop:measurability-of-map-into-one-point-extension}.
        \end{itemize}
    \end{proof}
\end{proposition}

Since we desire that the concatenated process $\concatenatedprocess$ has a Markov property as well, we need to control the \emph{transfer} from the space $\measurablespace_\indexsetelement$ to $\measurablespace_{\indexsetelement+1}$ in a way which does not destroy the Markov property of the concatenated processes $\killedlocalprocess^{(i)}$. In general, this can be done by fixing a suitable \emph{transfer kernel} from $\probabilityspace_\indexsetelement$ to $\measurablespace_{\indexsetelement+1}$. We need to choose such a transfer kernel for every $\indexsetelement\in\indexset_0$. If $\indexset$ is finite, then $\sup\indexset$ can be omitted, since there is no transfer after the last process has been executed.

Assuming enough path regularity, we can intuitively think of the transfer kernel as a description of the transfer from the \textit{exit point} \begin{equation}\label{eq:exit-point}
    \concatenatedprocess_{\tilde\lifetimesum_\indexsetelement-}=\tilde\killedlocalprocess^{(\indexsetelement)}_{\tilde\lifetime_\indexsetelement-}:=\lim_{\timepoint\to\tilde\lifetime_\indexsetelement-}\tilde\killedlocalprocess^{(\indexsetelement)}_\timepoint
\end{equation} to the \textit{revival point} \begin{equation}\label{eq:revival-point}
    \concatenatedprocess_{\tilde\lifetimesum_\indexsetelement}=\tilde\killedlocalprocess^{(\indexsetelement+1)}_0.
\end{equation} Formally, that description is given by prescribing a regular version $\regenerationdistribution_\indexsetelement$ of the conditional distribution of $\process_{\tilde\lifetimesum_\indexsetelement}$ given $\process_{\tilde\lifetimesum_\indexsetelement-}$; see \autoref{set:regeneration-distribution} and \autoref{ex:revival-point-distribution} below.

By fixation of the transfer kernel, we now got a complete probabilistic description of the evolution of $\concatenatedprocess$ at hand. It will be convenient to have a distinguished outcome for which the lifetime $\lifetime_\indexsetelement$ of the process $\killedlocalprocess^{(\indexsetelement)}$ has immediately ended.  We obtain this by adjoining a \emph{dead path} - i.e. a path for which the lifetime of the corresponding process has immediately ended - to each probability space $(\probabilityspace_\indexsetelement,\eventsystem_\indexsetelement)$.

Formally, we let \begin{enumerate}[$\circ$]
    \item $\left[\deadstate_\indexsetelement\right]\in\probabilityspace_\indexsetelement$ with \begin{equation}
        \killedlocalprocess^{(\indexsetelement)}_\timepoint\left(\left[\deadstate_\indexsetelement\right]\right)=\deadstate_\indexsetelement\;\;\;\text{for all }\timepoint\ge0;
    \end{equation}
    \item $\transferkernel_\indexsetelement$ be a Markov kernel with source $(\probabilityspace_\indexsetelement,\eventsystem_\indexsetelement)$ and target $(\measurablespace_{\indexsetelement+1},\measurablesystem_{\indexsetelement+1})$
\end{enumerate} for $\indexsetelement\in\indexset_0$.

\begin{setting}[revival distribution]\label{set:regeneration-distribution}
    Assume \autoref{set:left-regular-killedlocalprocess} and \begin{equation}
        \transferkernel_\indexsetelement(\outcome_\indexsetelement,\;\cdot\;)=\regenerationdistribution_\indexsetelement\left(\killedlocalprocess^{(\indexsetelement)}_{\lifetime_\indexsetelement-1}(\outcome_\indexsetelement)\;\cdot\;\right)\;\;\;\text{for all }\outcome_\indexsetelement\in\probabilityspace_\indexsetelement
    \end{equation} for some Markov kernel $\regenerationdistribution_\indexsetelement$ with source $(\measurablespace_\indexsetelement,\measurablesystem_\indexsetelement)$ and target $(\measurablespace_{\indexsetelement+1},\measurablesystem_{\indexsetelement+1})$ for all $\indexsetelement\in\indexset_0$.
    \begin{flushright}
        $\square$
    \end{flushright}
\end{setting}

For a purely technical reason, we need to assume that \begin{equation}
    \measurablespace_\indexsetelement\cap\measurablespace_\secondindexsetelement=\emptyset\;\;\;\text{for all }\indexsetelement,\secondindexsetelement\in\indexset\text{ with }\indexsetelement\ne\secondindexsetelement.
\end{equation} In practice, this is not a restriction. For example, if we actually want that all processes evolve on a single measurable space $(\measurablespace_0,\measurablesystem_0)$, we can always artificially distinguish them by adding an enumeration to it.

We can now define a probability measure on $(\probabilityspace,\eventsystem)$ with respect to which the concatenated process $\process$ will be again Markov:

\begin{definition}
    Let $\indexsetelement\in\indexset$ and $\point\in\measurablespace_\indexsetelement$ $\Rightarrow$ \begin{equation}
        \probabilitykernel_\indexsetelement(\point,\;\cdot\;):=\probabilitymeasure^{(\indexsetelement)}_\point.
    \end{equation} and \begin{equation}
        \probabilitymeasure_\point:=\probabilitykernel(\point,\;\cdot\;):=\bigotimes_{\secondindexsetelement=1}^{\indexsetelement-1}\delta_{\left[\deadstate_\secondindexsetelement\right]}\otimes\probabilitymeasure^{(\indexsetelement)}_\point\otimes\bigotimes_{\substack{\secondindexsetelement\in\indexset\\\secondindexsetelement>\indexsetelement}}\left(\transferkernel_{\secondindexsetelement-1}\probabilitykernel_\secondindexsetelement\right),
    \end{equation}
    \begin{flushright}
        $\square$
    \end{flushright}
\end{definition}


Let $\probabilityprojection_{\MakeUppercase\secondindexsetelement}$ denote the canonical projection of $\probabilityspace$ onto $\bigtimes_{\secondindexsetelement\in\MakeUppercase\secondindexsetelement}\probabilityspace_\secondindexsetelement$ for nonempty $\MakeUppercase\secondindexsetelement\subseteq\indexset$. 

\begin{lemma}\label{lem:concatenation-projection}
    Let $\thirdindexsetelement\in\indexset$ and $\projectionvariable\in\left(\bigotimes_{\secondindexsetelement=1}^\thirdindexsetelement\eventsystem_\secondindexsetelement\right)_b$ $\Rightarrow$ $\tilde\projectionvariable:=\projectionvariable\circ\probabilityprojection_{\{1,\;\ldots\;,\;\thirdindexsetelement\}}\in\eventsystem_b$ and \begin{equation}
        \expectation_\point\left[\tilde\projectionvariable\right]=\begin{cases}
            \begin{array}{l}
                \displaystyle\int\probabilitymeasure^{(\indexsetelement)}_\point\left[\dif\outcome_\indexsetelement\right]\int\left(\bigotimes_{\secondindexsetelement=\indexsetelement+1}^{\thirdindexsetelement-1}\left(\transferkernel_{\secondindexsetelement-1}\probabilitykernel_\secondindexsetelement\right)\right)\left(\outcome_{\indexsetelement+1},\ldots,\outcome_{\thirdindexsetelement-1}\right)\\
                \displaystyle\int\transferkernel_{\thirdindexsetelement-1}(\outcome_{\thirdindexsetelement-1},\dif\point_\thirdindexsetelement)\\
                \expectation^{(\thirdindexsetelement)}_{\point_\thirdindexsetelement}\left[\projectionvariable\left(\left[\deadstate_1\right],\ldots,\left[\deadstate_{\indexsetelement-1}\right],\outcome_\indexsetelement,\ldots,\outcome_{\thirdindexsetelement-1},\;\cdot\;\right)\right]
            \end{array}&\text{, if }\indexsetelement<\thirdindexsetelement;\\
            \expectation^{(\indexsetelement)}_\point\left[\projectionvariable\left(\left[\deadstate_1\right],\ldots,\left[\deadstate_{\indexsetelement-1}\right],\;\cdot\;\right)\right]&\text{, if }\indexsetelement=\thirdindexsetelement;\\
            \projectionvariable\left(\left[\deadstate_1\right],\ldots,\left[\deadstate_\thirdindexsetelement\right]\right)&\text{, otherwise}
        \end{cases}
    \end{equation} for all $\point\in\measurablespace_\indexsetelement$ for some $\indexsetelement\in\indexset$.
    \begin{proof}[Proof\textup:\nopunct]
        \leavevmode
        \begin{itemize}[$\circ$]
            \item Let $\indexsetelement\in\indexset$ and $\point\in\measurablespace_\indexsetelement$ $\Rightarrow$ \begin{equation}
                \begin{split}
                    \expectation_\point\left[\tilde\projectionvariable\right]=&\int\left(\bigotimes_{\secondindexsetelement=1}^{\indexsetelement-1}\delta_{\left[\deadstate_\secondindexsetelement\right]}\right)\left(\dif\integrand(\outcome_1,\ldots,\outcome_{\indexsetelement-1}\right)\int\probabilitymeasure^{(\indexsetelement)}_\point\left[\dif\outcome_\indexsetelement\right]\\&\int\left(\bigotimes_{\secondindexsetelement=\indexsetelement+1}^\thirdindexsetelement\left(\transferkernel_{\secondindexsetelement-1}\probabilitykernel_\secondindexsetelement\right)\right)\left(\outcome_{\indexsetelement+1},\ldots,\outcome_\thirdindexsetelement\right)\projectionvariable(\outcome_1,\ldots,\outcome_\thirdindexsetelement).
                \end{split}
            \end{equation}
            \item If $\indexsetelement<\thirdindexsetelement$, then \begin{equation}
                \begin{split}
                    \expectation_\point\left[\tilde\projectionvariable\right]&=\int\probabilitymeasure^{(\indexsetelement)}_\point\left[\dif\outcome_\indexsetelement\right]\int\left(\bigotimes_{\secondindexsetelement=\indexsetelement+1}^\thirdindexsetelement\left(\transferkernel_{\secondindexsetelement-1}\probabilitykernel_\secondindexsetelement\right)\right)\left(\outcome_{\indexsetelement+1},\ldots,\outcome_\thirdindexsetelement\right)\\&\;\hphantom=\;\projectionvariable\left(\left[\deadstate_1\right],\ldots,\left[\deadstate_{\indexsetelement-1}\right],\outcome_\indexsetelement,\ldots,\outcome_\thirdindexsetelement\right)\\
                    &=\int\probabilitymeasure^{(\indexsetelement)}_\point\left[\dif\outcome_\indexsetelement\right]\int\left(\bigotimes_{\secondindexsetelement=\indexsetelement+1}^{\thirdindexsetelement-1}\left(\transferkernel_{\secondindexsetelement-1}\probabilitykernel_\secondindexsetelement\right)\right)\left(\outcome_{\indexsetelement+1},\ldots,\outcome_{\thirdindexsetelement-1}\right)\\&\;\hphantom=\;\int\transferkernel_{\thirdindexsetelement-1}(\outcome_{\thirdindexsetelement-1},\dif\point_\thirdindexsetelement)\expectation^{(\thirdindexsetelement)}_{\point_\thirdindexsetelement}\left[\projectionvariable\left(\left[\deadstate_1\right],\ldots,\left[\deadstate_{\indexsetelement-1}\right],\outcome_\indexsetelement,\ldots,\outcome_{\thirdindexsetelement-1},\;\cdot\;\right)\right].
                \end{split}
            \end{equation}
            \item If $\indexsetelement=\thirdindexsetelement$, then \begin{equation}
                \expectation_\point\left[\tilde\projectionvariable\right]=\expectation^{(\indexsetelement)}_\point\left[\projectionvariable\left(\left[\deadstate_1\right],\ldots,\left[\deadstate_{\indexsetelement-1}\right],\;\cdot\;\right)\right].
            \end{equation}
            \item If $\indexsetelement>\thirdindexsetelement$, then \begin{equation}
                \expectation_\point\left[\tilde\projectionvariable\right]=\projectionvariable\left(\left[\deadstate_1\right],\ldots,\left[\deadstate_\thirdindexsetelement\right]\right).
            \end{equation}
        \end{itemize}
    \end{proof}
\end{lemma}

From the following corollary, we immediately conclude that if every transfer is state-independent, then $\largebracket(\tilde\killedlocalprocess^{(\indexsetelement)}\largebracket)_{\indexsetelement\in\indexset}$ is $\probabilitymeasure_{\regenerationdistribution_0}$-independent for all $\regenerationdistribution_0\in\mathcal M_1(\measurablespace_1,\measurablesystem_1)$:

\begin{corollary}\label{cor:concatenation-projection}
    Let $\indexsetelement\in\indexset$ and $\projectionvariable\in(\eventsystem_\indexsetelement)_b$ $\Rightarrow$ $\tilde\projectionvariable:=\projectionvariable\circ\probabilityprojection_\indexsetelement\in\eventsystem_b$ and \begin{equation}
        \expectation_{\point_\indexsetelement}\left[\tilde\projectionvariable\right]=\expectation^{(\indexsetelement)}_{\point_\indexsetelement}\left[\projectionvariable\right]\;\;\;\text{for all }\point_\indexsetelement\in\measurablespace_\indexsetelement.
    \end{equation} If $\indexsetelement\ge2$ and $\transferkernel_{\indexsetelement-1}$ is state-independent, then \begin{equation}
        \expectation_{\point_1}\left[\tilde\projectionvariable\right]=\expectation^{(\indexsetelement)}_{\transferkernel_{\indexsetelement-1}}\left[\projectionvariable\right]\;\;\;\text{for all }\point_1\in\measurablespace_1.
    \end{equation}
    \begin{proof}[Proof\textup:\nopunct]
        \autoref{lem:concatenation-projection} $\Rightarrow$ Claim.
    \end{proof}
\end{corollary}

As a first result, we show that the revival point \eqref{eq:revival-point} has actually the previously mentioned distribution:

\begin{theorem}\label{thm:revival-distribution}
    Let $\indexsetelement\in\indexset_0$ and $(\point_\indexsetelement,\integrand)\in\measurablespace_\indexsetelement\times\measurablesystem_b$ with $\probabilitymeasure_{\point_\indexsetelement}\left[\tilde\lifetime_{\indexsetelement+1}>0\right]=1$ $\Rightarrow$ \begin{equation}\label{eq:prop:revival-distribution}
        \expectation_{\point_\indexsetelement}\left[\integrand\left(\concatenatedprocess_{\tilde\lifetimesum_\indexsetelement}\right);\tilde\lifetimesum_\indexsetelement<\infty\expectationmid\filtration_{\tilde\lifetimesum_\indexsetelement-}\right]=1_{\left\{\;\tilde\lifetimesum_\indexsetelement\;<\;\infty\;\right\}}\left(\left.\transferkernel_\indexsetelement\integrand\right|_{\measurablespace_{\indexsetelement+1}}\circ\probabilityprojection_\indexsetelement\right).
    \end{equation}
    \textit{Remark}: \autoref{lem:concatenation-projection} $\Rightarrow$ \begin{equation}
        \probabilitymeasure_{\point_\indexsetelement}\left[\tilde\lifetime_{\indexsetelement+1}>0\right]=\probabilitymeasure^{(\indexsetelement+1)}_{\probabilitymeasure^{(\indexsetelement)}_{\point_\indexsetelement}\transferkernel_\indexsetelement}\left[\lifetime_{\indexsetelement+1}>0\right].
    \end{equation}
    \begin{proof}[Proof\textup:\nopunct]
        \leavevmode
        \begin{itemize}[$\circ$]
            \item Let $\timepoint\ge0$ and $\event\in\filtration_\timepoint$ $\Rightarrow$ $\exists\event_\indexsetelement\in\filtration^{(\indexsetelement)}_\timepoint$ with \begin{equation}\label{eq:revival-distribution-proof-eq1}
                \left(\event\cap\left\{\timepoint<\lifetimesum_\indexsetelement \right\}\right)_{(\outcome_1,\;\ldots\;,\;\outcome_{\indexsetelement+1})}=\event_\indexsetelement\times\bigtimes_{\substack{\secondindexsetelement\in\indexset\\\indexsetelement>\indexsetelement}}\Omega_\secondindexsetelement\;\;\;\text{for all }\outcome\in\left\{\lifetimesum_{\indexsetelement-1}\le\timepoint\right\}.
            \end{equation}
            \item $\tilde\lifetime_{\indexsetelement+1}>0$ $\probabilitymeasure_{\point_\indexsetelement}$-almost surely $\Rightarrow$ $\tilde\lifetimesum_\indexsetelement<\tilde\lifetimesum_{\indexsetelement+1}$ $\probabilitymeasure_{\point_\indexsetelement}$-almost surely $\Rightarrow$ \begin{equation}\label{eq:revival-distribution-proof-eq2}
                \concatenatedprocess_{\tilde\lifetimesum_\indexsetelement}=\tilde\killedlocalprocess^{(\indexsetelement+1)}_0
            \end{equation} and hence \begin{equation}
                \begin{array}{cl}
                    &\expectation_{\point_\indexsetelement}\left[1_{\left\{\;\tilde\lifetimesum_\indexsetelement\;<\;\infty\;\right\}}\integrand\left(\concatenatedprocess_{\tilde\lifetimesum_\indexsetelement}\right);\event\cap\left\{\timepoint<\tilde\lifetimesum_\indexsetelement\right\}\right]\\
                    \stackrel{\eqref{eq:revival-distribution-proof-eq1}+\eqref{eq:revival-distribution-proof-eq2}}=&\displaystyle\int\probabilitymeasure^{(\indexsetelement)}_{\point_\indexsetelement}\left[\dif\outcome_\indexsetelement\right]1_{\event_\indexsetelement}(\outcome_\indexsetelement)1_{\left\{\;\lifetime_\indexsetelement\;<\;\infty\;\right\}}(\outcome_\indexsetelement)\\
                    &\displaystyle\int\transferkernel_\indexsetelement(\outcome_\indexsetelement,\dif\point_{\indexsetelement+1})\expectation^{(\indexsetelement+1)}_{\point_{\indexsetelement+1}}\underbraceinsidebracket{[}{\integrand\Bigl(}{\killedlocalprocess^{(\indexsetelement+1)}_0}{=\;\point_{\indexsetelement+1}}{\Bigr)}{]}\\
                    =&\displaystyle\int\probabilitymeasure_{\point_\indexsetelement}\left[\dif\outcome_\indexsetelement\right]1_{\event_\indexsetelement}(\outcome_\indexsetelement)1_{\left\{\;\lifetime_\indexsetelement\;<\;\infty\;\right\}}(\outcome_\indexsetelement)\left(\left.\transferkernel_\indexsetelement\integrand\right|_{\measurablespace_{\indexsetelement+1}}\right)(\outcome_\indexsetelement)\\
                    \\
                    \stackrel{\eqref{eq:revival-distribution-proof-eq1}}=&\expectation_{\point_\indexsetelement}\left[1_{\left\{\;\tilde\lifetimesum_\indexsetelement\;<\;\infty\;\right\}}\left(\left.\transferkernel_\indexsetelement\integrand\right|_{\measurablespace_{\indexsetelement+1}}\circ\probabilityprojection_\indexsetelement\right);\event\cap\left\{\indexsetelement<\tilde\lifetimesum_\indexsetelement\right\}\right].
                \end{array}
            \end{equation}
        \end{itemize}
    \end{proof}
\end{theorem}

\begin{example}\label{ex:revival-point-distribution}
    Assume \autoref{set:regeneration-distribution}. Let $\indexsetelement\in\indexset_0$ and $\point_\indexsetelement\in\measurablespace_\indexsetelement$ with $\probabilitymeasure_{\point_\indexsetelement}\left[\tilde\lifetime_{\indexsetelement+1}>0\right]=1$ $\Rightarrow$ $\regenerationdistribution_\indexsetelement$ is a regular version of the conditional distribution of $\concatenatedprocess_{\tilde\lifetimesum_\indexsetelement}$ under $\probabilitymeasure_{\point_\indexsetelement}$ given $\concatenatedprocess_{\tilde\lifetimesum_\indexsetelement-}$ on $\left\{\tilde\lifetimesum_\indexsetelement<\infty\right\}$.
    \begin{proof}[Proof\textup:\nopunct]
        \autoref{thm:revival-distribution} $\Rightarrow$
        \begin{equation}
            \begin{split}
                \expectation_\point\left[1_{\left\{\;\tilde\lifetimesum_\indexsetelement\;<\;\infty\;\right\}}\integrand\left(\concatenatedprocess_{\tilde\lifetimesum_\indexsetelement}\right)\expectationmid\concatenatedprocess_{\tilde\lifetimesum_\indexsetelement-}\right]&=\expectation_\point\left[\expectation_\point\left[1_{\left\{\;\tilde\lifetimesum_\indexsetelement\;<\;\infty\;\right\}}\integrand\left(\concatenatedprocess_{\tilde\lifetimesum_\indexsetelement}\right)\expectationmid\filtration_{\tilde\lifetimesum_\indexsetelement-}\right]\expectationmid\concatenatedprocess_{\tilde\lifetimesum_\indexsetelement-}\right]\\
                &=\expectation_\point\left[1_{\left\{\;\tilde\lifetimesum_\indexsetelement\;<\;\infty\;\right\}}\left(\left.\transferkernel_\indexsetelement\integrand\right|_{\measurablespace_{\indexsetelement+1}}\circ\probabilityprojection_\indexsetelement\right)\expectationmid\concatenatedprocess_{\tilde\lifetimesum_\indexsetelement-}\right]\\
                &=\expectation_\point\left[1_{\left\{\;\tilde\lifetimesum_\indexsetelement\;<\;\infty\;\right\}}\left(\left.\regenerationdistribution_\indexsetelement\integrand\right|_{\measurablespace_{\indexsetelement+1}}\right)\left(\concatenatedprocess_{\tilde\lifetimesum_\indexsetelement-}\right)\expectationmid\concatenatedprocess_{\tilde\lifetimesum_\indexsetelement-}\right]\\
                &=1_{\left\{\;\tilde\lifetimesum_\indexsetelement\;<\;\infty\;\right\}}\left(\left.\regenerationdistribution_\indexsetelement\integrand\right|_{\measurablespace_{\indexsetelement+1}}\right)\left(\concatenatedprocess_{\tilde\lifetimesum_\indexsetelement-}\right).
            \end{split}
        \end{equation} for all $\integrand\in\measurablesystem_b$.
    \end{proof}
\end{example}

\begin{remark}
    Let $\indexsetelement,\thirdindexsetelement\in\indexset$ with $\indexsetelement\le\thirdindexsetelement$ $\Rightarrow$ \begin{equation}
        \probabilitykernel_{\indexsetelement,\;\thirdindexsetelement}:=\probabilitykernel_\indexsetelement\otimes\bigotimes_{\secondindexsetelement=\indexsetelement+1}^\thirdindexsetelement\left(\transferkernel_{\secondindexsetelement-1}\probabilitykernel_\secondindexsetelement\right)
    \end{equation} is a Markov kernel with source $(\measurablespace_\indexsetelement,\measurablesystem_\indexsetelement)$ and target $\left(\bigtimes_{\secondindexsetelement=\indexsetelement}^\thirdindexsetelement\probabilityspace_\secondindexsetelement,\bigoplus_{\secondindexsetelement=\indexsetelement}^\thirdindexsetelement\eventsystem_\secondindexsetelement\right)$ and hence \begin{equation}
        \probabilitymeasure^{(\indexsetelement,\;\thirdindexsetelement)}_{\point_\indexsetelement}:=\probabilitykernel_{\indexsetelement,\;\thirdindexsetelement}(\point_\indexsetelement,\;\cdot\;)
    \end{equation} is a probability measure on $\left(\bigtimes_{\secondindexsetelement=\indexsetelement}^\thirdindexsetelement\probabilityspace_\secondindexsetelement,\bigoplus_{\secondindexsetelement=\indexsetelement}^\thirdindexsetelement\eventsystem_\secondindexsetelement\right)$ for all $\point_\indexsetelement\in\measurablespace_\indexsetelement$.\footnote{If $\indexsetelement=\thirdindexsetelement$, then \begin{equation}
        \probabilitykernel_{\indexsetelement,\;\thirdindexsetelement}=\probabilitykernel_\indexsetelement
    \end{equation} and \begin{equation}
        \probabilitymeasure^{(\indexsetelement,\;\thirdindexsetelement)}_{\point_\indexsetelement}=\probabilitymeasure^{(\indexsetelement)}_{\point_\indexsetelement}\;\;\;\text{for all }\point_\indexsetelement\in\measurablespace_\indexsetelement.
    \end{equation}}
    \begin{flushright}
        $\square$
    \end{flushright}
\end{remark}

We now define the semigroup, for which we will show in \autoref{thm:concatenated-process-is-markov} that it is the transition semigroup of the concatenated process $\concatenatedprocess$:

\begin{definition}[transition semigroup of concatenated process]
    Let $\indexsetelement\in\indexset$ and $(\point_\indexsetelement,\integrand)\in\measurablespace_\indexsetelement\times\measurablesystem_b$ $\Rightarrow$ \begin{equation}
        (\concatenatedsemigroup_\timepoint\integrand)(\point_\indexsetelement):=\left(\killedlocalsemigroup^{(\indexsetelement)}_\timepoint\left.\integrand\right|_{\measurablespace_\indexsetelement}\right)(\point_\indexsetelement)+\sum_{\substack{\secondindexsetelement\in\indexset_0\\\secondindexsetelement\ge\indexsetelement}}\expectation^{(\indexsetelement,\;\secondindexsetelement)}_{\point_\indexsetelement}\left[\transferkernel_\secondindexsetelement\left(\killedlocalsemigroup^{(\secondindexsetelement+1)}_{\timepoint-\lifetimesum_{\indexsetelement,\;\secondindexsetelement}}\left.\integrand\right|_{\measurablespace_{\secondindexsetelement+1}}\right);\timepoint\ge\lifetimesum_{\indexsetelement,\;\secondindexsetelement}\right]
    \end{equation} for $\timepoint\ge0$.
    \begin{flushright}
        $\square$
    \end{flushright}
\end{definition}

\begin{example}[Markov process with restarts]
    Assume \begin{equation}
        \regenerationdistribution_\indexsetelement:=\transferkernel_\indexsetelement
    \end{equation} is state-independent - i.e. simply a probability measure on $(\measurablespace_{\indexsetelement+1},\measurablesystem_{\indexsetelement+1})$ - for all $\indexsetelement\in\indexset_0$ $\Rightarrow$ \begin{equation}
        \probabilitykernel_{\indexsetelement,\;\thirdindexsetelement}(\point_\indexsetelement,\;\cdot\;)=\probabilitymeasure^{(\indexsetelement)}_{\point_\indexsetelement}\otimes\bigotimes_{\secondindexsetelement=\indexsetelement+1}^\thirdindexsetelement\probabilitymeasure^{(\secondindexsetelement)}_{\regenerationdistribution_{\secondindexsetelement-1}}\;\;\;\text{for all }\point_\indexsetelement\in\measurablespace_\indexsetelement
    \end{equation} for all $\point_\indexsetelement\in\measurablespace_\indexsetelement$ and $\indexsetelement,\thirdindexsetelement\in\indexset$ with $\indexsetelement\le\thirdindexsetelement$. Let \begin{equation}
        \Sigma_{\point_\indexsetelement}^{(\indexsetelement,\;\thirdindexsetelement)}:=\probabilitykernel_{\indexsetelement,\;\thirdindexsetelement}(\point_\indexsetelement,\;\cdot\;)\circ\lifetimesum_{\indexsetelement,\;\thirdindexsetelement}^{-1}
    \end{equation} denote the distribution of $\lifetimesum_{\indexsetelement,\;\thirdindexsetelement}$ with respect to $\probabilitykernel_{\indexsetelement,\;\thirdindexsetelement}(\point_\indexsetelement,\;\cdot\;)$ for $\point_\indexsetelement\in\measurablespace_\indexsetelement$ and $\indexsetelement,\thirdindexsetelement\in\indexset$ with $\indexsetelement\le\thirdindexsetelement$ $\Rightarrow$ \begin{equation}
        \expectation^{(\indexsetelement,\;\secondindexsetelement)}_{\point_\indexsetelement}\left[\transferkernel_\secondindexsetelement\left(\killedlocalsemigroup^{(\secondindexsetelement+1)}_{\timepoint-\lifetimesum_{\indexsetelement,\;\secondindexsetelement}}\left.\integrand\right|_{\measurablespace_{\secondindexsetelement+1}}\right);\timepoint\ge\lifetimesum_{\indexsetelement,\;\secondindexsetelement}\right]=\int_0^\timepoint\Sigma_{\point_\indexsetelement}^{(\indexsetelement,\;\secondindexsetelement)}(\dif\prevtimepoint)\regenerationdistribution_\secondindexsetelement\killedlocalsemigroup^{(\secondindexsetelement)}_{\timepoint-\prevtimepoint}\left.\integrand\right|_{\measurablespace_{\secondindexsetelement+1}}
    \end{equation} for all  $\point_\indexsetelement\in\measurablespace_\indexsetelement$ and $\indexsetelement,\secondindexsetelement\in\indexset$ with $\indexsetelement\le\secondindexsetelement$. In particular, in \autoref{set:concatenation-of-killing-at-exponential-rate}, if $\killingrate_\indexsetelement=\killingrate$ for some $\killingrate>0$ for all $\indexsetelement\in\indexset$ and $\indexset=\mathbb N$, then \begin{equation}
        \Sigma_{\point_\indexsetelement}^{(\indexsetelement,\;\thirdindexsetelement)}=\operatorname{Gamma}(\thirdindexsetelement-\indexsetelement+1,\killingrate)
    \end{equation} for all $\point_\indexsetelement\in\measurablespace_\indexsetelement$ and $\indexsetelement\in\indexset$. If, additionally, $(\measurablespace_\indexsetelement,\measurablesystem_\indexsetelement)=(\measurablespace_0\times\{\indexsetelement\},\measurablesystem_0\times\{\indexsetelement\})$ for some measurable space $(\measurablespace_0,\measurablesystem_0)$, $\largebracket(\localprocess^{(\indexsetelement)}\largebracket)_{\indexsetelement\in\indexset}$ is identically distributed and $\regenerationdistribution_\indexsetelement=\regenerationdistribution$ for some $\regenerationdistribution\in\mathcal M_1(\measurablespace_0,\measurablesystem_0)$, then\footnote{$\largebracket(\localsemigroup^{(1)}_\timepoint\largebracket)_{\timepoint\ge0}$ can clearly be identified with a Markov semigroup on $(\measurablespace_0,\measurablesystem_0)$.} \begin{equation}
        \begin{split}
            (\concatenatedsemigroup_\timepoint\integrand)(\point)&=\left(\killedlocalsemigroup^{(1)}\integrand\right)(\point)+\killingrate\int_0^\timepoint\regenerationdistribution\killedlocalsemigroup^{(1)}_{\timepoint-\prevtimepoint}\integrand\\
            &=e^{-\killingrate\timepoint}\left(\localsemigroup^{(1)}_\timepoint\integrand\right)(\point)+\killingrate\int_0^\timepoint e^{-\killingrate\prevtimepoint}\regenerationdistribution\localsemigroup^{(1)}_\prevtimepoint\integrand
        \end{split}
    \end{equation} for all $(\point,\integrand)\in\measurablespace_0\times(\measurablesystem_0)_b$ and $\timepoint\ge0$. That is, we recover the result from \citep{avrachenkov2013restart}, where a Markov process is forced to restart from a given distribution $\regenerationdistribution$ at time moments generated by an independent Poisson process.
    \begin{flushright}
        $\square$
    \end{flushright}
\end{example}

While $\lifetime_\indexsetelement$ cannot "look in the future", it might have a certain "memory of the past". this memory would destroy the Markov property of the concatenated process. For that reason, we need to impose a "memorylessness" assumption. More precisely, we assume  \begin{enumerate}[(i)]
    \setcounter{enumi}{2}
    \item $\lifetime_\indexsetelement$ is a terminal time\footnote{see \autoref{def:terminal-time}.} on $\largebracket(\probabilityspace_\indexsetelement,\eventsystem_\indexsetelement\largebracket(\filtration^{(\indexsetelement)}_\timepoint\largebracket)_{\timepoint\ge0},\largebracket(\shift^{(\indexsetelement)}_\timepoint\largebracket)_{\timepoint\ge0},\largebracket(\probabilitymeasure^{(\indexsetelement)}_{\point_\indexsetelement}\largebracket)_{\point_\indexsetelement\in\measurablespace_\indexsetelement}\largebracket)$ and hence \begin{equation}\label{eq:concatenation-of-mp-terminal-time-assumption}
        \prevtimepoint+\lifetime_\indexsetelement\circ\shift^{(\indexsetelement)}_\prevtimepoint=\lifetime_\indexsetelement\;\;\;\probabilitymeasure^{(\indexsetelement)}_{\point_\indexsetelement}\text{-almost surely on }\{\;\prevtimepoint<\;\lifetime_\indexsetelement\;\};
    \end{equation}
    \item \begin{equation}\label{eq:transfer-kernel-memoryless}
        \transferkernel_\indexsetelement\secondintegrand\circ\shift^{(\indexsetelement)}_\prevtimepoint=\transferkernel_\indexsetelement\secondintegrand\;\;\;\text{for all }\secondintegrand\in(\measurablesystem_{\indexsetelement+1})_b\;\probabilitymeasure^{(\indexsetelement)}_{\point_\indexsetelement}\text{-almost surely on }\{\;\prevtimepoint<\lifetime_\indexsetelement\;\}
    \end{equation}
\end{enumerate} for all $\prevtimepoint\ge0$, $\point_\indexsetelement\in\measurablespace_\indexsetelement$ and $\indexsetelement\in\indexset$.

Before we proceed, we need a trivial but important, technical extension of \eqref{eq:concatenation-of-mp-terminal-time-assumption}:

\begin{lemma}
    Let $\indexsetelement\in\indexset_0$ $\Rightarrow$ \begin{equation}
        \left(\transferkernel_\indexsetelement\secondintegrand(\outcome,\;\cdot\;)\circ\shift^{(\indexsetelement)}_\prevtimepoint\right)(\outcome)=(\transferkernel_\indexsetelement \secondintegrand)(\outcome)
    \end{equation} for all $\secondintegrand\in(\eventsystem_\indexsetelement\otimes\measurablesystem_{\indexsetelement+1})_b$ and $\probabilitymeasure^{(\indexsetelement)}_{\point_\indexsetelement}$-almost all $\outcome\in\{\;\prevtimepoint<\lifetime_\indexsetelement\;\}$ for all $\prevtimepoint\ge0$ and $\point_\indexsetelement\in \measurablespace_\indexsetelement$.
    \begin{proof}[Proof\textup:\nopunct]
        Let $\prevtimepoint\ge0$ and $\point_\indexsetelement\in\measurablespace_\indexsetelement$ $\Rightarrow$ \begin{equation}
            \left(\transferkernel_\indexsetelement \secondintegrand\circ\shift^{(\indexsetelement)}_\prevtimepoint\right)(\outcome)=(\transferkernel_\indexsetelement \secondintegrand)(\outcome)\;\;\;\text{for all }\secondintegrand\in(\measurablesystem_{\secondindexsetelement})_b\text{ and }\outcome\in\{\;\prevtimepoint<\lifetime_\indexsetelement\;\}\setminus N
        \end{equation} for some $\probabilitymeasure^{(\indexsetelement)}_{\point_\indexsetelement}$-null set $N$ $\Rightarrow$ \begin{equation}
            \left(\transferkernel_\indexsetelement \secondintegrand(\outcome,\;\cdot\;)\circ\shift^{(\indexsetelement)}_\prevtimepoint\right)(\outcome)=\left(\transferkernel_\indexsetelement \secondintegrand(\outcome,\;\cdot\;)\right)(\outcome)\stackrel{\text{def}}=(\transferkernel_\indexsetelement \secondintegrand)(\outcome)
        \end{equation} for all $\secondintegrand\in(\eventsystem_\indexsetelement\otimes\measurablesystem_\secondindexsetelement)_b$ and $\outcome\in\{\;\prevtimepoint<\lifetime_\indexsetelement\;\}\setminus N$.
    \end{proof}
\end{lemma}

We are finally in the position to prove that the concatenation is actually a Markov process:

\begin{theorem}[concatenated process is Markov]\label{thm:concatenated-process-is-markov}
    Let $\prevtimepoint,\timepoint\ge0$ and $(\point,\integrand)\in\measurablespace\times\measurablesystem_b$ $\Rightarrow$ \begin{equation}
        \expectation_\point\left[\integrand\left(\concatenatedprocess_{\prevtimepoint+\timepoint}\right)\mid\filtration_\prevtimepoint\right]=\left(\concatenatedsemigroup_\timepoint\integrand\right)\left(\concatenatedprocess_\prevtimepoint\right).
    \end{equation}
    \begin{proof}[Proof\textup:\nopunct]
        \leavevmode
        \begin{itemize}[$\circ$]
            \item $\point\in\measurablespace$ $\Rightarrow$ $\exists\indexsetelement\in\indexset$ with $\point\in\indexset$.
            \item Let $\event\in\filtration_\prevtimepoint$ $\Rightarrow$\footnote{Note that \begin{equation}
                \prevtimepoint\ge\tilde\sigma_\secondindexsetelement\Rightarrow\prevtimepoint+\timepoint\ge\tilde\sigma_\secondindexsetelement\Rightarrow\prevtimepoint+\timepoint-\tilde\sigma_{\secondindexsetelement-1}\ge\tilde\sigma_\secondindexsetelement-\tilde\sigma_{\secondindexsetelement-1}=\tilde\tau_\secondindexsetelement\Rightarrow\integrand\left(\tilde\killedlocalprocess^{(\secondindexsetelement)}_{\prevtimepoint+\timepoint-\tilde\sigma_{\secondindexsetelement-1}}\right)=0
            \end{equation} and \begin{equation}
                \secondindexsetelement\le\indexsetelement-1\Rightarrow\probabilitymeasure_\point\left[\tilde\sigma_\secondindexsetelement=0\right]=1.
            \end{equation}} \begin{equation}
                \begin{split}
                    &\expectation\left[\integrand(\concatenatedprocess_{\prevtimepoint+\timepoint});\event\right]\\
                    &\;\;\;\;\;\;\;\;\;\;\;\;=\sum_{\secondindexsetelement\in\indexset}\expectation_\point\left[\integrand(\concatenatedprocess_{\prevtimepoint+\timepoint});\event\cap\{\;\tilde\lifetimesum_{\secondindexsetelement-1}\le\prevtimepoint+\timepoint<\tilde\lifetimesum_\secondindexsetelement\;\}\right]\\
                    &\;\;\;\;\;\;\;\;\;\;\;\;\;\;\;\;\;\;\;\;\;\;\;\;+\expectation_\point\left[\integrand(\concatenatedprocess_{\prevtimepoint+\timepoint});\event\cap\{\;\prevtimepoint+\timepoint\ge\tilde\lifetimesum_{\sup\indexset}\;\}\right]\\
                    &\;\;\;\;\;\;\;\;\;\;\;\;=\sum_{\secondindexsetelement\in\indexset}\expectation_\point\left[\integrand\left(\tilde\killedlocalprocess^{(\secondindexsetelement)}_{\prevtimepoint+\timepoint-\tilde\lifetimesum_{\secondindexsetelement-1}}\right);\event\cap\{\;\tilde\lifetimesum_{\secondindexsetelement-1}\le\prevtimepoint+\timepoint<\tilde\lifetimesum_\secondindexsetelement\;\}\right]\\
                    &\;\;\;\;\;\;\;\;\;\;\;\;\;\;\;\;\;\;\;\;\;\;\;\;+\underbrace{\integrand(\deadstate)}_{=\;0}\probabilitymeasure_\point\left[\event\cap\{\;\prevtimepoint+\timepoint\ge\tilde\lifetimesum_{\sup\indexset}\;\}\right]\\
                    &\;\;\;\;\;\;\;\;\;\;\;\;=\sum_{\secondindexsetelement\in\indexset}\expectation_\point\left[\integrand\left(\tilde\killedlocalprocess^{(\secondindexsetelement)}_{\prevtimepoint+\timepoint-\tilde\lifetimesum_{\secondindexsetelement-1}}\right);\event\cap\{\;\tilde\lifetimesum_{\secondindexsetelement-1}\le\prevtimepoint+\timepoint\;\}\cap\;\{\;\prevtimepoint<\tilde\lifetimesum_\secondindexsetelement\;\}\right].
                \end{split}
            \end{equation}
            \item Note that \begin{equation}
                T_{i,\:k}:=\bigotimes_{j=i+1}^k\left(T_{j-1}\probabilitykernel_j\right)\otimes T_k
            \end{equation} is a Markov kernel with source $(\Omega_i,\mathcal A_i)$ and target $\left(\bigtimes_{j=i+1}^k\Omega_j\times E_{k+1},\bigotimes_{j=i+1}^k\mathcal A_j\otimes\mathcal E_{k+1}\right)$ for all $i,k\in I_0$ with $i\le k$.
            \item \begin{listclaim}\label{thm:concatenation-of-mp-is-markov-proof-claim1}
                Let $k\in I$ with $k\ge i$ $\Rightarrow$ \begin{equation}
                    \begin{split}
                        &\operatorname E_x\left[f\left(\tilde\killedlocalprocess^{(k)}_{s+t-\tilde\sigma_{k-1}}\right);A\cap\{\;\tilde\sigma_{k-1}\le s+t\;\}\cap\{\;s<\tilde\sigma_k\;\}\right]\\&\;\;\;\;\;\;\;\;\;\;\;\;=\sum_{j=i}^{k-1}\operatorname E_x\left[\operatorname E^{(j)}_{X_s}\left[T_{j,\;k-1}\varphi\right];A\cap\{\;\tilde\sigma_{j-1}\le s<\tilde\sigma_j\;\}\right]\\&\;\;\;\;\;\;\;\;\;\;\;\;\;\;\;\;\;\;\;\;\;\;\;\;+\operatorname E_x\left[\left(Q^{(k)}_t\left.f\right|_{E_k}\right)(X_s);A\cap\{\;\tilde\sigma_{k-1}\le s<\tilde\sigma_k\;\}\right],
                    \end{split}
                \end{equation} where \begin{equation}
                    \begin{split}
                        &\varphi(\outcome_\indexsetelement,\ldots,\outcome_{\thirdindexsetelement-1},\point_\thirdindexsetelement):=\\
                        &\;\;\;\;\;\;\;\;\;\;\;\;1_{\left\{\;\timepoint\;\ge\;\lifetimesum_{\indexsetelement,\;\thirdindexsetelement-1}(\outcome_1,\;\ldots\;,\;\outcome_{\thirdindexsetelement-1})\;\right\}}\left(\killedlocalsemigroup^{(\thirdindexsetelement)}_{\timepoint-\lifetimesum_{\indexsetelement,\;\thirdindexsetelement-1}(\outcome_1,\;\ldots\;,\;\outcome_{\thirdindexsetelement-1})}\left.\integrand\right|_{\measurablespace_\thirdindexsetelement}\right)(\point_\thirdindexsetelement)
                    \end{split}
                \end{equation} for $(\outcome_\indexsetelement,\ldots,\outcome_{\thirdindexsetelement-1},\point_\thirdindexsetelement)\in\bigtimes_{\secondindexsetelement=\indexsetelement}^{\thirdindexsetelement-1}\probabilityspace_\secondindexsetelement\times\measurablespace_\thirdindexsetelement$.
                \begin{proof}[Proof\textup:\nopunct]
                    \leavevmode
                    \begin{itemize}[$\circ$]
                        \item $A\in\mathcal F_s$ $\Rightarrow$ \begin{equation}
                            (A\cap\{\;s<\tilde\sigma_k\;\})_{\omega_1,\;\ldots\;,\;\omega_{k-1}}=A^{\omega_1,\;\ldots\;,\;\omega_{k-1}}_k\times\bigtimes_{\substack{j\in I\\j>l}}\Omega_j
                        \end{equation} for some $A^{\omega_1,\;\ldots\;,\;\omega_{k-1}}_k\in\mathcal F^{(k)}_{s-\sigma_{k-1}(\omega_1,\;\ldots\;,\;\omega_{k-1})}$\newline for all $(\omega_1,\ldots,\omega_{k-1})\in\{\;\sigma_{k-1}\le s\;\}$ $\Rightarrow$\footnote{Let $(\omega_1,\ldots,\omega_k)\in\Omega_i\times\cdots\times\Omega_{k-1}$ with $s<\sigma_k([\Delta_1],\ldots,[\Delta_{i-1}],\omega_i,\ldots,\omega_k)$ $\Rightarrow$ $\exists!j\in\{1,\ldots,k\}$ with \begin{equation}
                            s\in\left[\sigma_{j-1}\left([\Delta_1],\ldots,[\Delta_{i-1}],\omega_i,\ldots,\omega_{j-1}\right),\sigma_j\left([\Delta_1],\ldots,[\Delta_{i-1}],\omega_i,\ldots,\omega_j\right)\right).
                        \end{equation}} \begin{equation}
                            \begin{split}
                                &\operatorname E_x\left[f\left(\tilde\killedlocalprocess^{(k)}_{s+t-\tilde\sigma_{k-1}}\right);A\cap\{\;\tilde\sigma_{k-1}\le s+t\;\}\cap\{\;s<\tilde\sigma_k\;\}\right]\\
                                &=\sum_{j=i}^k\operatorname E_x\left[f\left(\tilde\killedlocalprocess^{(k)}_{s+t-\tilde\sigma_{k-1}}\right);\begin{split}&A\cap\{\;\tilde\sigma_{j-1}\le s<\tilde\sigma_j\;\}\\&\hphantom A\cap\{\;\tilde\sigma_{k-1}\le s+t\;\}\end{split}\right]\\
                                &=\sum_{j=i}^k\int\left(\bigotimes_{j=1}^{i-1}\delta_{\left[\Delta_j\right]}\right)\left(\dif\left(\omega_1,\ldots,\omega_{i-1}\right)\right)\int\probabilitykernel_i(x,\dif\omega_i)\\&\;\;\;\;\;\;\;\;\int\left(\bigotimes_{j=i+1}^k\left(T_{j-1}\probabilitykernel_j\right)\right)\left(\omega_i,\dif\left(\omega_{i+1},\ldots,\omega_k\right)\right)\\&\;\;\;\;\;\;\;\;1_{\{\;\sigma_{j-1}\;\le\;s\;\}}(\omega_1,\ldots,\omega_{j-1})1_{\{\;\sigma_{k-1}\;\le\;s+t\;\}}(\omega_1,\ldots,\omega_{k-1})\\&\;\;\;\;\;\;\;\;f\left(\killedlocalprocess^{(k)}_{s+t-\sigma_{k-1}(\omega_1,\;\ldots\;,\;\omega_{k-1})}(\omega_k)\right)\\&\;\;\;\;\;\;\;\;\int\left(\bigotimes_{\substack{j\in I\\j>k}}\left(T_{j-1}\probabilitykernel_j\right)\right)\left(\omega_k,\dif(\omega_j)_{\substack{j\in I\\j>k}}\right)\underbrace{1_{A\;\cap\;\{\;s\;<\;\tilde\sigma_j\;\}}\left(\left(\omega_j\right)_{j\in I}\right)}_{=\;1_{A^{\omega_1,\;\ldots\;,\;\omega_{j-1}}_j}\left(\omega_j\right)}
                            \end{split}
                        \end{equation}
                        \item If $k=i$, then \begin{equation}
                            \begin{split}
                                &\operatorname E_x\left[f\left(\tilde\killedlocalprocess^{(k)}_{s+t-\tilde\sigma_{k-1}}\right);A\cap\{\;\tilde\sigma_{k-1}\le s+t\;\}\cap\{\;s<\tilde\sigma_k\;\}\right]\\
                                &=\underbrace{1_{\{\;\sigma_{i-1}\;\le\;s+t\;\}}([\Delta_1],\ldots,[\Delta_{i-1}])}_{=\;1}\\&\;\;\;\;\;\;\;\;\int\probabilitykernel_i(x,\dif\omega_i)f\left(\killedlocalprocess^{(i)}_{s+t-\underbrace{\sigma_{i-1}\left([\Delta_1],\;\ldots\;,\;[\Delta_{i-1}]\right)}_{=\;0}}(\omega_i)\right)\\&\;\;\;\;\;\;\;\;\int\left(\bigotimes_{j>i}\left(T_{j-1}\probabilitykernel_j\right)\right)\left(\omega_i,\dif\left(\omega_j\right)_{\substack{j\in I\\j>i}}\right)\\&\;\;\;\;\;\;\;\;\underbrace{1_{A\;\cap\;\{\;s\;<\;\tilde\sigma_i\;\}}([\Delta_1],\ldots,[\Delta_{i-1}],\left(\omega_j\right)_{\substack{j\in I\\j\ge i}}}_{=\;1_{A_i}(\omega_i}\\
                                &=\operatorname E^{(i)}_x\left[f\left(\killedlocalprocess^{(i)}_{s+t}\right);A_i\right]\\&=\operatorname E^{(i)}_x\left[\left(Q^{(i)}_t\left.f\right|_{E_i}\right)\left(\killedlocalprocess^{(i)}_s\right);A_i\right]\\
                                &=\underbrace{1_{\{\;\sigma_{i-1}\;\le\;s\;\}}([\Delta_1],\ldots,[\Delta_{i-1}])}_{=\;1}\\&\;\;\;\;\;\;\;\;\int\probabilitykernel_i(x,\dif\omega_i)\left(Q^{(i)}_t\left.f\right|_{E_i}\right)\left(\killedlocalprocess^{(i)}_{s-\underbrace{\sigma_{i-1}\left([\Delta_1],\;\ldots\;,\;[\Delta_{i-1}]\right)}_{=\;0}}(\omega_i)\right)\\&\;\;\;\;\;\;\;\;\int\left(\bigotimes_{j>i}\left(T_{j-1}\probabilitykernel_j\right)\right)\left(\omega_i,\dif\left(\omega_j\right)_{\substack{j\in I\\j>i}}\right)\\&\;\;\;\;\;\;\;\;\underbrace{1_{A\;\cap\;\{\;s\;<\;\tilde\sigma_i\;\}}([\Delta_1],\ldots,[\Delta_{i-1}],\left(\omega_j\right)_{\substack{j\in I\\j\ge i}}}_{=\;1_{A_i}(\omega_i}\\
                                &=\operatorname E_x\left[\left(Q^{(i)}_t\left.f\right|_{E_i}\right)\left(\tilde\killedlocalprocess^{(i)}_{s-\tilde\sigma_{i-1}}\right);A\cap\left\{\;\tilde\sigma_{i-1}\le s<\tilde\sigma_i\;\right\}\right].
                            \end{split}
                        \end{equation}
                        \item If $k>i$, then \begin{equation}
                            \begin{split}
                                &\operatorname E_x\left[f\left(\tilde\killedlocalprocess^{(k)}_{s+t-\tilde\sigma_{k-1}}\right);A\cap\{\;\tilde\sigma_{k-1}\le s+t\;\}\cap\{\;s<\tilde\sigma_k\;\}\right]\\
                                &=\sum_{a=i}^k\int\left(\bigotimes_{j=1}^{i-1}\delta_{\left[\Delta_j\right]}\right)\left(\dif\left(\omega_1,\ldots,\omega_{i-1}\right)\right)\int\probabilitykernel_i(x,\dif\omega_i)\\&\;\;\;\;\;\;\;\;\int\left(\bigotimes_{j=i+1}^{k-1}\left(T_{j-1}\probabilitykernel_j\right)\right)\left(\omega_i,\dif\left(\omega_{i+1},\ldots,\omega_{k-1}\right)\right)\\&\;\;\;\;\;\;\;\;1_{\{\;\sigma_{a-1}\;\le\;s\;\}}(\omega_1,\ldots,\omega_{a-1})1_{\{\;\sigma_{k-1}\;\le\;s+t\;\}}(\omega_1,\ldots,\omega_{k-1})\\&\;\;\;\;\;\;\;\;\int T_{k-1}(\omega_{k-1},\dif x_k)\int\probabilitykernel_k(x_k,\dif\omega_k)\\&\;\;\;\;\;\;\;\;1_{A^{\omega_1,\;\ldots\;,\;\omega_{a-1}}_a}(\omega_a)f\left(\killedlocalprocess^{(k)}_{s+t-\sigma_{k-1}(\omega_1,\;\ldots\;,\;\omega_{k-1})}(\omega_k)\right)\\
                                &=\sum_{a=i}^{k-1}\int\left(\bigotimes_{j=1}^{i-1}\delta_{\left[\Delta_j\right]}\right)\left(\dif\left(\omega_1,\ldots,\omega_{i-1}\right)\right)\int\probabilitykernel_i(x,\dif\omega_i)\\&\;\;\;\;\;\;\;\;\int\left(\bigotimes_{j=i+1}^{k-1}\left(T_{j-1}\probabilitykernel_j\right)\right)\left(\omega_i,\dif\left(\omega_{i+1},\ldots,\omega_{k-1}\right)\right)\\&\;\;\;\;\;\;\;\;1_{\{\;\sigma_{a-1}\;\le\;s\;\}}(\omega_1,\ldots,\omega_{a-1})1_{A^{\omega_1,\;\ldots\;,\;\omega_{a-1}}_a}(\omega_a)\\&\;\;\;\;\;\;\;\;1_{\{\;\sigma_{k-1}\;\le\;s+t\;\}}(\omega_1,\ldots,\omega_{k-1})\\&\;\;\;\;\;\;\;\;\int T_{k-1}(\omega_{k-1},\dif x_k)\underbrace{\operatorname E^{(k)}_{x_k}\left[f\left(\killedlocalprocess^{(k)}_{s+t-\sigma_{k-1}(\omega_1,\;\ldots\;,\;\omega_{k-1})}\right)\right]}_{=\;\left(Q^{(k)}_{s+t-\sigma_{k-1}(\omega_1,\;\ldots\;,\;\omega_{k-1}}\left.f\right|_{E_k}\right)(x_k)}\\
                                &+\int\left(\sum_{j=1}^{i-1}\delta_{\left[\Delta_j\right]}\right)\left(\dif\left(\omega_1,\ldots,\omega_{i-1}\right)\right)\int\probabilitykernel_i(x,\dif\omega_i)\\&\;\;\;\;\;\;\;\;\int\left(\bigotimes_{j=i+1}^{k-1}\left(T_{j-1}\probabilitykernel_j\right)\right)\left(\omega_i,\dif\left(\omega_{i+1},\ldots,\omega_{k-1}\right)\right)\\&\;\;\;\;\;\;\;\;1_{\{\;\sigma_{k-1}\;\le\;s\;\}}(\omega_1,\ldots,\omega_{k-1})\\&\;\;\;\;\;\;\;\;\int T_{k-1}(\omega_{k-1},\dif x_k)\\&\;\;\;\;\;\;\;\;\underbrace{\operatorname E^{(k)}_{x_k}\left[f\left(\killedlocalprocess^{(k)}_{s+t-\sigma_{k-1}(\omega_1,\;\ldots\;,\;\omega_{k-1})}\right);A^{\omega_1,\;\ldots\;,\;\omega_{k-1}}_k\right]}_{=\;\operatorname E^{(k)}_{x_k}\left[\left(Q^{(k)}_t\left.f\right|_{E_k}\right)\left(\killedlocalprocess^{(k)}_{s-\sigma_{k-1}(\omega_1,\;\ldots\;,\;\omega_{k-1}}\right);A^{\omega_1,\;\ldots\;,\;\omega_{k-1}}_k\right]}\\
                                &=:\operatorname I+\operatorname{II}
                            \end{split}
                        \end{equation}
                        \item \begin{listclaim}
                            \begin{equation}
                                \operatorname I=\sum_{a=i}^{k-1}\operatorname E_x\left[\operatorname E^{(a)}_{X_s}\left[T_{a,\;k-1}\varphi\right];A\cap\{\;\tilde\sigma_{a-1}\le s<\tilde\sigma_a\;\}\right].
                            \end{equation}
                            \begin{proof}[Proof\textup:\nopunct]
                                \leavevmode
                                \begin{itemize}[$\sq$]
                                    \item Let \begin{equation}
                                        \begin{split}
                                            &g_{\omega_1,\;\ldots\;\omega_{a-1}}(\omega_a,x_{a+1}):=\int\probabilitykernel_{a+1}(x_{a+1},\dif\omega_{a+1})\\&\;\;\;\;\;\;\;\;\int\left(\bigotimes_{j=a+2}^{k-1}\left(T_{j-1}\probabilitykernel_j\right)\right)\left(\omega_{a+1},\dif\left(\omega_{a+2},\ldots,\omega_{k-1}\right)\right)\\&\;\;\;\;\;\;\;\;1_{\{\;\sigma_{k-1}\;\le\;s+t\;\}}(\omega_1,\ldots,\omega_{k-1})\int T_{k-1}(\omega_{k-1},\dif x_k)\\&\;\;\;\;\;\;\;\;\left(Q^{(k)}_{s+t-\sigma_{k-1}(\omega_1,\;\ldots\;,\;\omega_{k-1})}\left.f\right|_{E_k}\right)(x_k)
                                        \end{split}
                                    \end{equation} for $(\omega_a,x_{a+1})\in\Omega_a\times E_{a+1}$ and $(\omega_1,\ldots,\omega_{a-1})\in\bigtimes_{j=1}^{a-1}\Omega_j$ and \begin{equation}
                                        \begin{split}
                                            &h(\omega_a,x_{a+1}):=\int\probabilitykernel_{a+1}(x_{a+1},\dif\omega_{a+1})\\&\;\;\;\;\;\;\;\;\int\left(\bigotimes_{j=a+2}^{k-1}\left(T_{j-1}\probabilitykernel_j\right)\right)\left(\omega_{a+1},\dif\left(\omega_{a+2},\ldots,\omega_{k-1}\right)\right)\\&\;\;\;\;\;\;\;\;1_{\left\{\;\tau_a\;\le\;t-\sum_{j=a+1}^{k-1}\tau_j\left(\omega_j\right)\;\right\}}(\omega_a)\int T_{k-1}(\omega_{k-1},\dif x_k)\\&\;\;\;\;\;\;\;\;\left(Q^{(k)}_{t-\sum_{j=a}^{k-1}\tau_j\left(\omega_j\right)}\left.f\right|_{E_k}\right)(x_k)
                                        \end{split}
                                    \end{equation} for $(\omega_a,x_{a+1})\in\Omega_a\times E_{a+1}$.
                                    \item \eqref{eq:concatenation-of-mp-terminal-time-assumption} $\Rightarrow$ \begin{equation}
                                        s+\tau_a\left(\theta^{(a)}_{s-\sigma_{a-1}(\omega_1,\;\ldots\;,\;\omega_{a-1}}(\omega_a)\right)=\sigma_a(\omega_1,\ldots,\omega_a)
                                    \end{equation} for $\operatorname P^{(a)}_{x_a}$-almost all $\omega_a\in\Omega_a$ with $s<\sigma_a(\omega_1,\ldots,\omega_a)$\newline for all $(\omega_1,\ldots,\omega_{a-1})\in\bigtimes_{j=1}^{a-1}\Omega_j$ with $\sigma_{a-1}(\omega_1,\ldots,\omega_{a-1})\le s$\newline and $x_a\in E_a$ $\Rightarrow$ \begin{equation}
                                        h\left(\theta^{(a)}_{s-\sigma_{a-1}(\omega_1,\;\ldots\;,\;\omega_{a-1}}(\omega_a),x_{a+1}\right)=g_{\omega_1,\;\ldots\;,\;\omega_{a-1}}(\omega_a,x_{a+1})
                                    \end{equation} for $\operatorname P^{(a)}_{x_a}$-almost all $\omega_a\in\Omega_a$ with $s<\sigma_a(\omega_1,\ldots,\omega_a)$\newline for all $(\omega_1,\ldots,\omega_{a-1})\in\bigtimes_{j=1}^{a-1}\Omega_j$ with $\sigma_{a-1}(\omega_1,\ldots,\omega_{a-1})\le s$\newline and $x_a\in E_a$.
                                    \item \eqref{eq:transfer-kernel-memoryless} $\Rightarrow$ \begin{equation}
                                        \begin{split}
                                            &\left(T_ag_{\omega_1,\;\ldots\;,\;\omega_{a-1}}\right)(\omega_a)\\
                                            &\;\;\;\;\;\;\;\;\;\;\;\;=\left(T_ag_{\omega_1,\;\ldots\;,\;\omega_{a-1}}(\omega_a,\;\cdot\;)\circ\theta^{(a)}_{s-\sigma_{a-1}(\omega_1,\;\ldots\;,\;\omega_{a-1}}\right)(\omega_a)\\
                                            &\;\;\;\;\;\;\;\;\;\;\;\;=(T_ah)\left(\theta^{(a)}_{s-\sigma_{a-1}(\omega_1,\;\ldots\;,\;\omega_{a-1}}(\omega_a)\right)
                                        \end{split}
                                    \end{equation} for $\operatorname P^{(a)}_{x_a}$-almost all $\omega_a\in\Omega_a$ with $s<\sigma_a(\omega_1,\ldots,\omega_a)$\newline for all $(\omega_1,\ldots,\omega_{a-1})\in\bigtimes_{j=1}^{a-1}\Omega_j$ with $\sigma_{a-1}(\omega_1,\ldots,\omega_{a-1})\le s$\newline and $x_a\in E_a$ $\Rightarrow$ \begin{equation}
                                        \begin{split}
                                            &\operatorname E^{(a)}_{x_a}\left[T_ag_{\omega_1,\;\ldots\;,\;\omega_{a-1}};A^{\omega_1,\;\ldots\;\omega_{a-1}}_a\right]\\
                                            &\;\;\;\;\;\;\;\;\;\;\;\;=\operatorname E^{(a)}_{x_a}\left[T_ah\circ\theta^{(a)}_{s-\sigma_{a-1}(\omega_1,\;\ldots\;,\;\omega_{a-1}};A^{\omega_1,\;\ldots\;\omega_{a-1}}_a\right]\\
                                            &\;\;\;\;\;\;\;\;\;\;\;\;=\operatorname E^{(a)}_{x_a}\left[\operatorname E^{(a)}_{\killedlocalprocess^{(a)}_{s-\sigma_{a-1}(\omega_1,\;\ldots\;,\;\omega_{a-1}}}\left[T_ah\right];A^{\omega_1,\;\ldots\;\omega_{a-1}}_a\right]
                                        \end{split}
                                    \end{equation} for all $(\omega_1,\ldots,\omega_{a-1})\in\bigtimes_{j=1}^{a-1}\Omega_j$ with $\sigma_{a-1}(\omega_1,\ldots,\omega_{a-1})\le s$ and $x_a\in E_a$.
                                    \item Let \begin{equation}
                                        \varphi(\omega_a,\ldots,\omega_{k-1},x_k):=1_{\left\{\;t\;\ge\;\sum_{j=a}^{k-1}\tau_j\left(\omega_j\right)\;\right\}}\left(Q^{(k)}_{t-\sum_{j=a}^{k-1}\tau_j\left(\omega_j\right)}\left.f\right|_{E_k}\right)(x_k)
                                    \end{equation} for $(\omega_a,\ldots,\omega_{k-1},x_k)\in\bigtimes_{j=a}^{k-1}\Omega_j\times E_k$ $\Rightarrow$ \begin{equation}
                                        T_{a,\;k-1}\varphi=T_ah.
                                    \end{equation}
                                    \item Thus, \begin{equation}
                                        \begin{split}
                                            &\operatorname I=\sum_{a=i}^{k-1}\int\left(\bigotimes_{j=1}^{i-1}\delta_{\left[\Delta_j\right]}\right)\left(\dif\left(\omega_1,\ldots,\omega_{i-1}\right)\right)\int\probabilitykernel_i(x,\dif\omega_i)\\&\;\;\;\;\;\;\;\;\int\left(\bigotimes_{j=i+1}^{a-1}\left(T_{j-1}\probabilitykernel_j\right)\right)\left(\omega_i,\dif\left(\omega_{i+1},\ldots,\omega_{a-1}\right)\right)\\&\;\;\;\;\;\;\;\;1_{\{\;\sigma_{a-1}\;\le\;s\;\}}(\omega_1,\ldots,\omega_{a-1})\\&\;\;\;\;\;\;\;\;\int T_{a-1}(\omega_{a-1},\dif x_a)\\&\;\;\;\;\;\;\;\;\underbrace{\operatorname E^{(a)}_{x_a}\left[T_ag_{\omega_1,\;\ldots\;,\;\omega_{a-1}};A^{\omega_1,\;\ldots\;,\;\omega_{a-1}}_a\right]}_{=\;\operatorname E^{(a)}_{x_a}\left[\operatorname E^{(a)}_{\killedlocalprocess^{(a)}_{s-\sigma_{a-1}(\omega_1,\;\ldots\;,\;\omega_{a-1})}}\left[T_ah\right];A^{\omega_1,\;\ldots\;,\;\omega_{a-1}}_a\right]}\\
                                            &\hphantom{\operatorname I}=\sum_{a=i}^{k-1}\int\left(\bigotimes_{j=1}^{i-1}\delta_{\left[\Delta_j\right]}\right)\left(\dif\left(\omega_1,\ldots,\omega_{i-1}\right)\right)\int\probabilitykernel_i(x,\dif\omega_i)\\&\;\;\;\;\;\;\;\;\int\left(\bigotimes_{j=i+1}^{a-1}\left(T_{j-1}\probabilitykernel_j\right)\right)\left(\omega_i,\dif\left(\omega_{i+1},\ldots,\omega_{a-1}\right)\right)\\&\;\;\;\;\;\;\;\;1_{\{\;\sigma_{a-1}\;\le\;s\;\}}(\omega_1,\ldots,\omega_{a-1})\\&\;\;\;\;\;\;\;\;\int T_{a-1}(\omega_{a-1},\dif x_a)\\&\;\;\;\;\;\;\;\;\operatorname E^{(a)}_{x_a}\left[\operatorname E^{(a)}_{\killedlocalprocess^{(a)}_{s-\sigma_{a-1}(\omega_1,\;\ldots\;,\;\omega_{a-1})}}\left[T_{a,\;k-1}\varphi\right];A^{\omega_1,\;\ldots\;,\;\omega_{a-1}}_a\right]\\
                                            &\hphantom{\operatorname I}=\sum_{a=i}^{k-1}\int\left(\bigotimes_{j=1}^{i-1}\delta_{\left[\Delta_j\right]}\right)\left(\dif\left(\omega_1,\ldots,\omega_{i-1}\right)\right)\int\probabilitykernel_i(x,\dif\omega_i)\\&\;\;\;\;\;\;\;\;\int\left(\bigotimes_{j=i+1}^{a-1}\left(T_{j-1}\probabilitykernel_j\right)\right)\left(\omega_i,\dif\left(\omega_{i+1},\ldots,\omega_{a-1}\right)\right)\\&\;\;\;\;\;\;\;\;1_{\{\;\sigma_{a-1}\;\le\;s\;\}}(\omega_1,\ldots,\omega_{a-1})\\&\;\;\;\;\;\;\;\;\int T_{a-1}(\omega_{a-1},\dif x_a)\int\probabilitykernel_a(x_a,\dif\omega_a)\\&\;\;\;\;\;\;\;\;\operatorname E^{(a)}_{\killedlocalprocess^{(a)}_{s-\sigma_{a-1}(\omega_1,\;\ldots\;,\;\omega_{a-1})}}\left[T_{a,\;k-1}\varphi\right]\\&\;\;\;\;\;\;\;\;\int\left(\bigotimes_{\substack{j\in I\\j>a}}\left(T_{j-1}\probabilitykernel_j\right)\right)\left(\omega_a,\dif\left(\omega_j\right)_{\substack{j\in I\\j>a}}\right)\\&\;\;\;\;\;\;\;\;\underbrace{1_{A\;\cap\;\{\;s\;<\;<\tilde\sigma_a\;\}}\left(\left(\omega_j\right)_{j\in I}\right)}_{=\;1_{A^{\omega_1,\;\ldots\;,\;\omega_{a-1}}_a}}\\
                                             &\hphantom{\operatorname I}=\sum_{a=i}^{k-1}\operatorname E_x\left[\operatorname E^{(a)}_{\tilde\killedlocalprocess^{(a)}_{s-\tilde\sigma_{a-1}}}\left[T_{a,\;k-1}\varphi\right];A\cap\left\{\;\tilde\sigma_{a-1}\le s<\tilde\sigma_a\;\right\}\right]
                                        \end{split}
                                    \end{equation} \begin{equation*}
                                        \begin{split}
                                             &\hphantom{\operatorname I}=\operatorname E_x\left[\sum_{a=i}^{k-1}1_{\left\{\;\tilde\sigma_{a-1}\;\le\;s\;<\;\tilde\sigma_a\;\right\}}\operatorname E^{(a)}_{\tilde\killedlocalprocess^{(a)}_{s-\tilde\sigma_{a-1}}}\left[T_{a,\;k-1}\varphi\right];A\right]\\
                                             &\hphantom{\operatorname I}=\operatorname E_x\left[\sum_{a=i}^{k-1}1_{\left\{\;\tilde\sigma_{a-1}\;\le\;s\;<\;\tilde\sigma_a\;\right\}}\operatorname E^{(a)}_{X_s}\left[T_{a,\;k-1}\varphi\right];A\right].
                                        \end{split}
                                    \end{equation*}
                                \end{itemize}
                            \end{proof}
                        \end{listclaim}
                        \begin{listclaim}
                            \begin{equation}
                                \operatorname{II}=\operatorname E_x\left[\left(Q^{(k)}_t\left.f\right|_{E_k}\right)(X_s);A\cap\left\{\;\tilde\sigma_{k-1}\le s<\tilde\sigma_k\;\right\}\right].
                            \end{equation}
                            \begin{proof}[Proof\textup:\nopunct]
                                \begin{equation}
                                    \begin{split}
                                        &\operatorname{II}=\int\left(\bigotimes_{j=1}^{i-1}\delta_{\left[\Delta_j\right]}\right)\left(\dif\left(\omega_1,\ldots,\omega_{i-1}\right)\right)\int\probabilitykernel_i(x,\dif\omega_i)\\&\;\;\;\;\;\;\;\;\int\left(\bigotimes_{j=i+1}^{k-1}\left(T_{j-1}\probabilitykernel_j\right)\right)\left(\omega_i,\dif\left(\omega_{i+1},\ldots,\omega_{k-1}\right)\right)\\&\;\;\;\;\;\;\;\;1_{\{\;\sigma_{k-1}\;\le\;s\;\}}(\omega_1,\ldots,\omega_{k-1})\int T_{k-1}(\omega_{k-1},\dif x_k)\\&\;\;\;\;\;\;\;\;\operatorname E^{(k)}_{x_k}\left[\left(Q^{(k)}_t\left.f\right|_{E_k}\right)\left(\killedlocalprocess^{(k)}_{s-\sigma_{k-1}(\omega_1,\;\ldots\;,\;\omega_{k-1})}\right);A^{\omega_1,\;\ldots\;,\;\omega_{k-1}}_k\right]\\
                                        &\hphantom{\operatorname{II}}=\int\left(\bigotimes_{j=1}^{i-1}\delta_{\left[\Delta_j\right]}\right)\left(\dif\left(\omega_1,\ldots,\omega_{i-1}\right)\right)\int\probabilitykernel_i(x,\dif\omega_i)\\&\;\;\;\;\;\;\;\;\int\left(\bigotimes_{j=i+1}^{k-1}\left(T_{j-1}\probabilitykernel_j\right)\right)\left(\omega_i,\dif\left(\omega_{i+1},\ldots,\omega_{k-1}\right)\right)\\&\;\;\;\;\;\;\;\;1_{\{\;\sigma_{k-1}\;\le\;s\;\}}(\omega_1,\ldots,\omega_{k-1})\int T_{k-1}(\omega_{k-1},\dif x_k)\\&\;\;\;\;\;\;\;\;\int\probabilitykernel_k(x_k,\dif\omega_k)\left(Q^{(k)}_t\left.f\right|_{E_k}\right)\left(\killedlocalprocess^{(k)}_{s-\sigma_{k-1}(\omega_1,\;\ldots\;,\;\omega_{k-1})}(\omega_k)\right)\\&\;\;\;\;\;\;\;\;\int\left(\bigotimes_{\substack{j\in I\\j>k}}\left(T_{j-1}\probabilitykernel_j\right)\right)\left(\omega_k,\dif\left(\omega_j\right)_{\substack{j\in I\\j>k}}\right)\underbrace{1_{A\;\cap\;\{\;s\;<\;\tilde\sigma_k\;\}}}_{=\;1_{A^{\omega_1,\;\ldots\;,\;\omega_{k-1}}_k(\omega_k)}}\\
                                        &\hphantom{\operatorname{II}}=\operatorname E_x\left[\left(Q^{(k)}_t\left.f\right|_{E_k}\right)\left(\tilde\killedlocalprocess^{(k)}_{s-\tilde\sigma_{k-1}}\right);A\cap\left\{\;\tilde\sigma_{k-1}\le s<\tilde\sigma_k\;\right\}\right]\\
                                        &\hphantom{\operatorname{II}}=\operatorname E_x\left[\left(Q^{(k)}_t\left.f\right|_{E_k}\right)(X_s);A\cap\left\{\;\tilde\sigma_{k-1}\le s<\tilde\sigma_k\;\right\}\right].
                                    \end{split}
                                \end{equation}
                            \end{proof}
                        \end{listclaim}
                    \end{itemize}
                \end{proof}
            \end{listclaim}
            \item \autoref{thm:concatenation-of-mp-is-markov-proof-claim1} $\Rightarrow$ \begin{equation}
                \begin{split}
                    &\expectation_\point\left[\integrand(\concatenatedprocess_{\prevtimepoint+\timepoint};\event\right]\\
                    &\;\;\;\;\;\;\;\;\;\;\;\;=\sum_{\substack{\thirdindexsetelement\in\indexset\\\thirdindexsetelement>\indexsetelement}}\sum_{\fourthindexsetelement=i}^{\thirdindexsetelement-1}\expectation_\point\left[\expectation^{(\fourthindexsetelement)}_{\concatenatedprocess_\prevtimepoint}\left[\transferkernel_{\fourthindexsetelement,\;\thirdindexsetelement-1}\varphi\right];\event\cap\left\{\;\tilde\lifetimesum_{\fourthindexsetelement-1}\le \prevtimepoint<\tilde\lifetimesum_\fourthindexsetelement\;\right\}\right]\\
                    &\;\;\;\;\;\;\;\;\;\;\;\;\;\;\;\;\;\;\;\;\;\;\;\;+\sum_{\substack{\thirdindexsetelement\in\indexset\\\thirdindexsetelement\ge\indexsetelement}}\expectation_\point\left[\left(\killedlocalsemigroup^{(\thirdindexsetelement)}_\timepoint\left.\integrand\right|_{\measurablespace_\thirdindexsetelement}\right)(\concatenatedprocess_\prevtimepoint);\event\cap\left\{\;\tilde\lifetimesum_{\thirdindexsetelement-1}\le\prevtimepoint<\tilde\lifetimesum_\thirdindexsetelement\;\right\}\right]\\
                    &\;\;\;\;\;\;\;\;\;\;\;\;=\sum_{\substack{\fourthindexsetelement\in\indexset\\\fourthindexsetelement\ge\indexsetelement}}\sum_{\substack{\thirdindexsetelement\in\indexset\\\thirdindexsetelement>\fourthindexsetelement}}\expectation_\point\left[\expectation^{(\fourthindexsetelement)}_{\concatenatedprocess_\prevtimepoint}\left[\transferkernel_{\fourthindexsetelement,\;\thirdindexsetelement-1}\varphi\right];\event\cap\left\{\;\tilde\lifetimesum_{\fourthindexsetelement-1}\le \prevtimepoint<\tilde\lifetimesum_\fourthindexsetelement\;\right\}\right]\\
                    &\;\;\;\;\;\;\;\;\;\;\;\;\;\;\;\;\;\;\;\;\;\;\;\;+\sum_{\substack{\thirdindexsetelement\in\indexset\\\thirdindexsetelement\ge\indexsetelement}}\expectation_\point\left[\left(\killedlocalsemigroup^{(\thirdindexsetelement)}_\timepoint\left.\integrand\right|_{\measurablespace_\thirdindexsetelement}\right)(\concatenatedprocess_\prevtimepoint);\event\cap\left\{\;\tilde\lifetimesum_{\thirdindexsetelement-1}\le \prevtimepoint<\tilde\lifetimesum_\thirdindexsetelement\;\right\}\right]\\
                    &\;\;\;\;\;\;\;\;\;\;\;\;=\sum_{\substack{\thirdindexsetelement\in\indexset\\\thirdindexsetelement\ge\indexsetelement}}\expectation_\point\left[(\concatenatedsemigroup_\timepoint\integrand)(\concatenatedprocess_\prevtimepoint);\event\cap\left\{\;\tilde\lifetimesum_{\thirdindexsetelement-1}\le \prevtimepoint<\tilde\lifetimesum_\thirdindexsetelement\;\right\}\right]\\
                    &\;\;\;\;\;\;\;\;\;\;\;\;=\sum_{\substack{\thirdindexsetelement\in\indexset\\\thirdindexsetelement\ge\indexsetelement}}\expectation_\point\left[(\concatenatedsemigroup_\timepoint\integrand)(\concatenatedprocess_\prevtimepoint);\event\cap\left\{\;\tilde\lifetimesum_{\thirdindexsetelement-1}\le \prevtimepoint<\tilde\lifetimesum_\thirdindexsetelement\;\right\}\right]\\
                    &\;\;\;\;\;\;\;\;\;\;\;\;\;\;\;\;\;\;\;\;\;\;\;\;+\underbrace{(\concatenatedsemigroup_\timepoint\integrand)(\deadstate)}_{=\;0}\probabilitymeasure_\point\left[\event\cap\{\;\prevtimepoint\ge\tilde\lifetimesum\;\}\right]\\
                    &\;\;\;\;\;\;\;\;\;\;\;\;=\expectation_\point\left[(\concatenatedsemigroup_\timepoint\integrand)(\concatenatedprocess_\prevtimepoint);\event\right].
                \end{split}
            \end{equation}
        \end{itemize}
    \end{proof}
\end{theorem}

We close our investigation of the concatenation of Markov process by the identification of its generator.

Let $\killedlocalgenerator_\indexsetelement$ denote the pointwise generator of $\largebracket(\killedlocalsemigroup^{(\indexsetelement)}_\timepoint\largebracket)_{\timepoint\ge0}$ for $\indexsetelement\in\indexset$.

\begin{remark}
    If $\sup\indexset\in\indexset$, then \begin{equation}
        \frac{\left(\concatenatedsemigroup_\timepoint\integrand\right)(\point)-\integrand(\point)}\timepoint\xrightarrow{\timepoint\to0+}\left(\killedlocalgenerator_{\sup\indexset}\left.\integrand\right|_{\measurablespace_{\sup\indexset}}\right)(\point)
    \end{equation} for all $\point\in\measurablespace_{\sup\indexset}$ and $\integrand\in\measurablesystem_b$ with $\left.\integrand\right|_{\measurablespace_{\sup\indexset}}\in\mathcal D\left(\killedlocalgenerator_{\sup\indexset}\right)$.
    \begin{flushright}
        $\square$
    \end{flushright}
\end{remark}

In order to obtain a satisfactory result, we need to impose several assumptions; mostly on the distribution of the exit points, but also on the shape of the transfer kernels. More precisely, we assume that \autoref{set:regeneration-distribution} and \begin{equation}
    \regenerationdistribution_\indexsetelement\integrand_{\indexsetelement+1}\in C_b(\measurablespace_\indexsetelement)\;\;\;\text{for all }\integrand_{\indexsetelement+1}\in C_b(\measurablespace_{\indexsetelement+1})
\end{equation} for all $\indexsetelement\in\indexset_0$. Let $\largebracket(\MakeUppercase\killingrate^{(\indexsetelement)}_\timepoint\largebracket)_{\timepoint\ge0}$ be an $[0,\infty)$-valued product measurable process on $(\probabilityspace_\indexsetelement,\eventsystem_\indexsetelement)$ with \begin{enumerate}[(i)]
    \setcounter{enumi}{4}
    \item\label{asm:concatenation-v} $\largebracket(\MakeUppercase\killingrate^{(\indexsetelement)}_\timepoint\largebracket)_{\timepoint\ge0}$ is (right-)continuous at $0$;
    \item\label{asm:concatenation-vi}  \begin{equation}
        \exists\timepoint_0>0:\sup_{(\outcome,\;\timepoint)\in\probabilityspace\times[0,\;\timepoint_0)}\MakeUppercase\killingrate^{(\indexsetelement)}_\timepoint(\outcome)<\infty;
    \end{equation}
    \item\label{asm:concatenation-vii}  \begin{equation}
        \probabilitymeasure^{(\indexsetelement)}_\point\left[\left(\lifetime_\indexsetelement,\killedlocalprocess^{(\indexsetelement)}_{\lifetime_\indexsetelement-}\right)\in\timedomain\times\measurableset\right]
        =\int_\timedomain\expectation^{(\indexsetelement)}_\point\left[\MakeUppercase\killingrate^{(\indexsetelement)}_\timepoint;\killedlocalprocess^{(\indexsetelement)}_{\timepoint-}\in\measurableset\right]\dif\timepoint
    \end{equation} for all $(\measurableset,\timedomain)\in\measurablesystem_\indexsetelement\times\mathcal B([0,\infty))$ and $\point\in\measurablespace_\indexsetelement$
\end{enumerate} and \begin{equation}\label{eq:concatenation-killing-rate-definition}
    \killingrate_\indexsetelement(\point):=\expectation^{(\indexsetelement)}_\point\left[\MakeUppercase\killingrate^{(\indexsetelement)}_0\right]\;\;\;\text{for }\point\in\measurablespace_\indexsetelement
\end{equation} for $\indexsetelement\in\indexset_0$.

\begin{example}
    Assume \autoref{set:concatenation-of-killing-at-exponential-rate} $\Rightarrow$ \begin{equation}
        \MakeUppercase\killingrate^{(\indexsetelement)}_\timepoint:=\killingrate_\indexsetelement\largebracket(\localprocess^{(\indexsetelement)}_\timepoint\largebracket)\multiplicativefunctional^{(\indexsetelement)}_\timepoint\;\;\;\text{for }\timepoint\ge0
    \end{equation} does satisfy \ref{asm:concatenation-v}-\ref{asm:concatenation-vii} and \eqref{eq:concatenation-killing-rate-definition} is consistent with the usage the symbol $\killingrate_\indexsetelement$ in \autoref{set:concatenation-of-killing-at-exponential-rate} for all $\indexsetelement\in\indexset_0$.
    \begin{proof}[Proof\textup:\nopunct]
        \autoref{ex:killed-process-lifefime-joint-distribution-in-set-additive-functional-induced-by-integral} $\Rightarrow$ Claim.
    \end{proof}
\end{example}

Let $\closedsubspace_\indexsetelement$ be a closed subspace of $(\measurablesystem_\indexsetelement)_b$ with\footnote{i.e. $\largebracket(\killedlocalsemigroup^{(\indexsetelement)}_\timepoint\largebracket)_{\timepoint\ge0}$ is strongly continuous on $\closedsubspace_\indexsetelement$.} \begin{equation}
    \left\|\killedlocalsemigroup^{(\indexsetelement)}_\timepoint\integrand-\integrand\right\|_\infty\xrightarrow{\timepoint\to0+}\;\;\;\text{for all }\integrand\in\closedsubspace_\indexsetelement
\end{equation} for $\indexsetelement\in\indexset_0\setminus\{1\}$.


We finally are able to state and prove the main result of this work:

\begin{theorem}[pointwise generator of concatenated process]\label{thm:generator-of-concatenated-process}
    \normalfont Let $\indexsetelement\in\indexset_0$ and $(\point,\integrand)\in\measurablespace_\indexsetelement\times\measurablesystem_b$ $\Rightarrow$ If $\left.\integrand\right|_{\measurablespace_\indexsetelement}\in C(\measurablespace_\indexsetelement)\cap\mathcal D\left(\killedlocalgenerator_\indexsetelement\right)$ and $\left.\integrand\right|_{\measurablespace_{\indexsetelement+1}}\in C(\measurablespace_{\indexsetelement+1})\cap\closedsubspace_{\indexsetelement+1}$, then \begin{equation}\label{eq:generator-of-concatenated-process}
        \frac{\left(\concatenatedsemigroup_\timepoint\integrand\right)(\point)-\integrand(\point)}\timepoint\xrightarrow{\timepoint\to0+}\left(\killedlocalgenerator_\indexsetelement\left.\integrand\right|_{\measurablespace_\indexsetelement}\right)(\point)+\killingrate_\indexsetelement(\point)\left(\left.\regenerationdistribution_\indexsetelement\integrand\right|_{\measurablespace_{\indexsetelement+1}}\right)(\point);
	\end{equation}
    \begin{proof}[Proof\textup:\nopunct]
        \leavevmode
        \begin{itemize}[$\circ$]
			\item For simplicity of notation, we assume $I=\{1,2\}$.
            \item \begin{listclaim}
                \normalfont Let $t\ge0$ $\Rightarrow$
                \begin{proof}[Proof\textup:\nopunct]
                    \leavevmode
                    \begin{itemize}[$\circ$]
                        \item Let \begin{equation}
                            \begin{array}{rcl}
                                \operatorname I_t&:=&\left(f(X_t)-f(x)\right)1_{\left\{\:t\:<\:\tilde\tau_1\:\right\}};\\\operatorname{II}_t&:=&\left(f(X_{\tilde\tau_1})-f\left(X_{\tilde\tau_1-}\right)\right)1_{\left\{\:t\:\ge\:\tilde\tau_1\:\right\}};\\\operatorname{III}_t&:=&\left(f(X_{\tilde\tau_1-})-f(x)+f(X_t)-f\left(X_{\tilde\tau_1}\right)\right)1_{\left\{\:t\:\ge\:\tilde\tau_1\:\right\}}.
                            \end{array}
                        \end{equation}
                        \item \begin{listclaim}
                            \normalfont\begin{equation}
                                \operatorname E_x\left[\operatorname I_t\right]=\begin{cases}\left(Q^{(1)}_t\left.f\right|_{E_1}\right)(x)-f(x)\underbrace{\operatorname P^{(1)}_x\left[t<\tau_1\right]}_{=\:Q^{(1)}_t(x,\:E_1)}&\text{, if }x\in E_1;\\0&\text{, if }x\in x_2.\end{cases}
                            \end{equation}
                            \begin{proof}[Proof\textup:\nopunct]
                                \leavevmode
                                \begin{itemize}[$\circ$]
                                    \item If $x\in E_1$, then \begin{equation}
                                        \begin{array}{rcl}
                                            \operatorname E_x\left[f(X_t);t<\tilde\tau_1\right]&=&\operatorname E_x\left[f\left(\tilde \killedlocalprocess^{(1)}_t\right);t<\tilde\tau_1\right]\\&=&\operatorname E^{(1)}_x\left[f\left(\killedlocalprocess^{(1)}_t\right);t<\tau_1\right]\\&=&\operatorname E^{(1)}_x\left[f\left(\killedlocalprocess^{(1)}_t\right)\right]=\left(Q^{(1)}_t\left.f\right|_{E_1}\right)(x).
                                        \end{array}
                                    \end{equation}
                                    \item If $x\in E_j$ for some $j\in I$ with $j>i$, then \begin{equation}
                                        \operatorname E_x\left[f(X_t);t<\tilde\tau_1\right]=f\left(\killedlocalprocess^{(1)}_t([\Delta_1])\right)\underbrace{1_{\{\:t\:<\:\tau_1\:\}}([\Delta_1])}_{=\:0}
                                    \end{equation}
                                    Thus, \begin{equation}
                                        \begin{split}
                                            &\operatorname E_x\left[\operatorname I_t\right]=\operatorname E_x\left[f(X_t)-f(x);t<\tilde\tau_1\right]\\&\;\;\;\;\;\;\;\;\;\;\;\;=\begin{cases}\underbrace{\operatorname E_x\left[f(X_t);\tilde\tau_1\right]}_{=\:\left(Q^{(1)}_t\left.f\right|_{E_1}\right)(x)}-f(x)\underbrace{\operatorname P_x\left[t<\tilde\tau_1\right]}_{=\:Q^{(1)}_t(x,\:E_1)}&\text{, if }x\in E_1;\\\underbrace{\operatorname E_x\left[f(X_t);t<\tilde\tau_1\right]}_{=\:0}-f(x)\underbrace{\operatorname P_x\left[t<\tilde\tau_1\right]}_{=\:0}&\text{, if }x\in E_2.\end{cases}
                                        \end{split}
                                    \end{equation}
                                \end{itemize}
                            \end{proof}
                        \end{listclaim}
                        \item \begin{listclaim}
                            \normalfont\begin{equation}
                                \operatorname E_x\left[\operatorname{II}\right]=\begin{cases}\operatorname E^{(1)}_x\left[T_1\left.f\right|_{E_2}-f\left(\killedlocalprocess^{(1)}_{\tau_1-}\right);t\ge\tau_1\right]&\text{, if }x\in E_1;\\f(x)&\text{, if }x\in E_2.\end{cases}
                            \end{equation}
                            \leavevmode
                            \begin{itemize}[$\sq$]
                                \item Note that \begin{equation}
                                    \begin{split}
                                        &\operatorname E_x\left[f\left(X_{\tilde\tau_1}\right);t\ge\tilde\tau_1\mid\mathcal F_{\tilde\tau_1-}\right]\\
                                        &\;\;\;\;\;\;\;\;\;\;\;\;\begin{array}{@{}c@{}l@{}}
                                            =&1_{\{\:t\:\ge\:\tilde\tau_1\:\}}\operatorname E_x\left[f\left(X_{\tilde\tau_1}\right);\tilde\tau_1<\infty\mid\mathcal F_{\tilde\tau_1-}\right]\\=&1_{\{\:t\:\ge\:\tilde\tau_1\:\}}\begin{cases}T_1\left.f\right|_{E_2}\circ p_1&\text{, if }x\in E_1;\\f(x)&\text{, if }x\in E_2.\end{cases}
                                        \end{array}
                                    \end{split}
                                \end{equation}
                                \item $1_{\{\:t\:\ge\:\tilde\tau_1\:\}}f\left(\tilde \killedlocalprocess^{(1)}_{\tilde\tau_1-}\right)$ is $\mathcal F_{\tilde\tau_1-}$-measurable $\Rightarrow$ \begin{equation*}
                                    \begin{split}
                                        &\operatorname E_x\left[\operatorname{II}\right]=\operatorname E_x\left[f\left(X_{\tilde\tau_1}\right)-f\left(\tilde \killedlocalprocess^{(1)}_{\tilde\tau_1-}\right);t\ge\tilde\tau_1\right]\\
                                        &\;\;\;\;=\operatorname E_x\left[\operatorname E_x\left[f\left(X_{\tilde\tau_1}\right)-f\left(\tilde \killedlocalprocess^{(1)}_{\tilde\tau_1-}\right);t\ge\tilde\tau_1\:\middle|\:\mathcal F_{\tilde\tau_1-}\right]\right]\\
                                        &\;\;\;\;=\begin{cases}\operatorname E_x\left[\left(T_1\left.f\right|_{E_2}-f\left(\killedlocalprocess^{(1)}_{\tau_1-}\right)\right)\circ p_1;t\ge\tau_1\circ p_1\right]&\text{, if }x\in E_1;\\\begin{split}&f(x)\operatorname P_x\left[t\ge\tau_1\circ p_1\right]\\&\;\;\;\;-\operatorname E_x\left[f\left(\killedlocalprocess^{(1)}_{\tau_1-}\right)\circ p_1;t\ge\tau_1\circ p_1\right]\end{split}&\text{, if }x\in E_2.\end{cases}\\
                                        &\;\;\;\;=\begin{cases}\operatorname E^{(1)}_x\left[T_1\left.f\right|_{E_2}-f\left(\killedlocalprocess^{(1)}_{\tau_1-}\right);t\ge\tau_1\right]&\text{, if }x\in E_1;\\\underbrace{1_{\{\:t\:\ge\:\tau_1\:\}}([\Delta_1])}_{=\:1}\left(f(x)-\underbrace{f\left(\underbrace{\killedlocalprocess^{(1)}_{\tau_1-}([Delta_1])}_{=\:\Delta_1}\right)}_{=\:0}\right)&\text{, if }x\in E_2.\end{cases}
                                    \end{split}
                                \end{equation*}
                            \end{itemize}
                        \end{listclaim}
                        \item \begin{listclaim}
                            \normalfont\begin{equation*}
                                \operatorname E_x\left[\operatorname{III}_t\right]=\begin{cases}\begin{array}{ll}&\operatorname E^{(1)}_x\left[f\left(\killedlocalprocess^{(1)}_{\tau_1-}\right)-f(x);t\ge\tau_1\right]\\&\;\;\;\;+\operatorname E^{(1)}_x\left[T_1\left(Q^{(2)}_{t-\tau_1}\left.f\right|_{E_2}-\left.f\right|_{E_2}\right);t\ge\tau_1\right]\end{array}&\text{, if }x\in E_1;\\\left(Q^{(2)}_t\left.f\right|_{E_2}\right)(x)-2f(x)&\text{, if }x\in E_2.\end{cases}
                            \end{equation*}
                            \begin{proof}[Proof\textup:\nopunct]
                                \leavevmode
                                \begin{itemize}[$\sq$]
                                    \item If $x\in E_1$, then \begin{equation}
                                        \begin{split}
                                            &\operatorname E_x\left[f(X_t);t\ge\tilde\tau_1\right]\\
                                            &\;\;\;\;=\int\pi_1(x,\dif\omega_1)1_{\{\:t\:\ge\:\tau_1\:\}}(\omega_1)\\&\;\;\;\;\;\;\;\;\int(T_1\pi_2)(\omega_1,\dif\omega_2)f\left(\killedlocalprocess^{(2)}_{t-\tau_1(\omega_1)}(\omega_2)\right)\\
                                            &\;\;\;\;=\int\pi_1(x,\dif\omega_1)1_{\{\:t\:\ge\:\tau_1\:\}}(\omega_1)\\&\;\;\;\;\;\;\;\;\int T_1(\omega_1,\dif x_2)f\underbrace{\operatorname E^{(2)}_{x_2}\left[\left(\killedlocalprocess^{(2)}_{t-\tau_1(\omega_1)}\right)\right]}_{=\:\left(Q^{(2)}_{t-\tau_1(\omega_1)}\left.f\right|_{E_2}\right)(x_2)}\\
                                            &\;\;\;\;=\int\pi_1(x,\dif\omega_1)1_{\{\:t\:\ge\:\tau_1\:\}}(\omega_1)\left(T_1\left(Q^{(2)}_{t-\tau_1}\left.f\right|_{E_2}\right)\right)(\omega_1)\\
                                            &\;\;\;\;=\operatorname E^{(1)}_x\left[T_1\left(Q^{(2)}_{t-\tau_1}\left.f\right|_{E_2}\right);t\ge\tau_1\right].
                                        \end{split}
                                    \end{equation}
                                    \item If $x\in E_2$, then \begin{equation}
                                        \begin{split}
                                            &\operatorname E_x\left[f(X_t);t\ge\tilde\tau_1\right]\\
                                            &\;\;\;\;\;\;\;\;\;\;\;\;\begin{split}
                                                &=\underbrace{1_{\{\:t\:\ge\:\tau_1\:\}}([\Delta_1])}_{=\:1}\int\pi_2(x,\dif\omega_2)f\left(\killedlocalprocess^{(2)}_{t-\underbrace{\tau_1([\Delta_1])}_{=\:0}}(\omega_2)\right)\\&=\operatorname E^{(2)}_x\left[f(\killedlocalprocess^{(2)}_t\right]=\left(Q^{(2)}_t\left.f\right|_{E_2}\right)(x).
                                            \end{split}
                                        \end{split}
                                    \end{equation}
                                    \item Now, \begin{equation}
                                        \begin{split}
                                            \operatorname E_x\left[\operatorname{III}_t\right]&\stackrel{\text{def}}=\operatorname E_x\left[f\left(\tilde \killedlocalprocess^{(1)}_{\tilde\tau_1-}\right)-f(x)+f(X_t)-f(X_{\tilde\tau_1})\right]\\
                                            &=\operatorname E_x\left[f\left(\killedlocalprocess^{(1)}_{\tau_1-}\right)\circ p_1;t\ge\tau_1\circ p_1\right]\\&\;\;\;\;\;\;\;\;\;\;\;\;-f(x)\operatorname P_x\left[t\ge\tilde\tau_1\right]\\&\;\;\;\;\;\;\;\;\;\;\;\;+\operatorname E_x\left[f(X_t);t\ge\tilde\tau_1\right]\\&\;\;\;\;\;\;\;\;\;\;\;\;-\operatorname E_x\left[f(X_{\tilde\tau_1});t\ge\tilde\tau_1\right].
                                        \end{split}
                                    \end{equation}
                                    \item If $x\in E_1$, then \begin{equation}
                                        \begin{split}
                                            \operatorname E_x\left[\operatorname{III}_t\right]&=\operatorname E^{(1)}_x\left[f\left(\killedlocalprocess^{(1)}_{\tau_1-}\right);t\ge\tau_1\right]\\&\;\;\;\;\;\;\;\;\;\;\;\;-f(x)\left(1-Q^{(1)}_t(x,E_1)\right)\\&\;\;\;\;\;\;\;\;\;\;\;\;+\operatorname E^{(1)}_x\left[T_1\left(Q^{(2)}_{t-\tau_1}\left.f\right|_{E_2}\right);t\ge\tau_1\right]\\&\;\;\;\;\;\;\;\;\;\;\;\;-\operatorname E^{(1)}_x\left[T_1\left.f\right|_{E_2};t\ge\tau_1\right].
                                        \end{split}
                                    \end{equation}
                                    \item If $x\in E_2$, then \begin{equation}
                                        \begin{split}
                                            \operatorname E_x\left[\operatorname{III}_t\right]&=\underbrace{f\left(\underbrace{\killedlocalprocess^{(1)}_{\tau_1-}([\Delta_1])}_{=\:\Delta_1}\right)}_{=\:0}\underbrace{1_{\{\:t\:\ge\:\tau_1\:\}}}_{=\:1}\\&\;\;\;\;\;\;\;\;\;\;\;\;-f(x)\\&\;\;\;\;\;\;\;\;\;\;\;\;+\left(Q^{(2)}_t\left.f\right|_{E_2}\right)(x)\\&\;\;\;\;\;\;\;\;\;\;\;\;-f(x).
                                        \end{split}
                                    \end{equation}
                                \end{itemize}
                            \end{proof}
                        \end{listclaim}
                    \end{itemize}
                \end{proof}
            \end{listclaim}
            \item \begin{listclaim}
                \normalfont If $x\in E_1$, then \begin{equation}
                    \frac{\operatorname E_x\left[\operatorname I_t+\operatorname{II}_t+\operatorname{III}_t\right]}t\xrightarrow{t\to0+}\left(A^{(1)}_{\textnormal p}\left.f\right|_{E_1}\right)(x)+\operatorname E^{(1)}_x\left[C^{(1)}_0\right]\left(\mu\left.f\right|_{E_2}\right)(x).
                \end{equation}
                \begin{proof}[Proof\textup:\nopunct]
                    \leavevmode
                    \begin{itemize}[$\circ$]
                        \item \begin{listclaim}
                            \begin{equation}
                                \frac{\operatorname E_x\left[\operatorname I_t\right]}t\xrightarrow{t\to0+}\left(A^{(1)}_{\textnormal p}\left.f\right|_{E_1}\right)(x)+f(x)\operatorname E^{(1)}_x\left[C^{(1)}_0\right].
                            \end{equation}
                            \begin{proof}[Proof\textup:\nopunct]
                                \leavevmode
                                \begin{itemize}[$\sq$]
                                    \item Note that \begin{equation}
                                        \operatorname P^{(1)}_x\left[t\ge\tau_1\right]=\operatorname E^{(1)}_x\left[\int_0^tC^{(1)}_s\dif s\right]
                                    \end{equation} and \begin{equation}
                                        \operatorname E_x\left[\operatorname I_t\right]=\left(Q^{(1)}_t\left.f\right|_{E_1}\right)(x)-f(x)+f(x)\operatorname P^{(1)}_x\left[t\ge\tau_1\right]
                                    \end{equation} for all $t\ge0$.
                                    \item $\left(C^{(1)}_t\right)_{t\ge0}$ is (right-)continuous at $0$ $\Rightarrow$ \begin{equation}
                                        \frac1t\int_0^tC^{(1)}_s\dif s\xrightarrow{t\to0+}C^{(1)}_0.
                                    \end{equation}
                                    \item By assumption,  \begin{equation}
                                        M:=\sup_{(\omega,\:t)\in\Omega\times[0,\:t_0)}C^{(1)}_t(\omega)<\infty
                                    \end{equation} for some $t_0>0$ and hence \begin{equation}
                                        0\le\frac1t\int_0^tC^{(1)}_s\dif s\le M\;\;\;\text{for all }t\in[0,t_0].
                                    \end{equation}
                                    \item Lebesgue’s dominated convergence theorem $\Rightarrow$ \begin{equation}
                                        \frac1t\operatorname E^{(1)}_x\left[\int_0^tC^{(1)}_s\dif s\right]\xrightarrow{t\to0+}\operatorname E^{(1)}_x\left[C^{(1)}_0\right].
                                    \end{equation}
                                \end{itemize}
                            \end{proof}
                        \end{listclaim}
                        \item \begin{listclaim}
                            \begin{equation}
                                \frac{\operatorname E_x\left[\operatorname{II}_t\right]}t\xrightarrow{t\to0+}\operatorname E^{(1)}_x\left[C^{(1)}_0\right]\left(\left.\mu f\right|_{E_2}-f\right)(x).
                            \end{equation}
                            \begin{proof}[Proof\textup:\nopunct]
                                \leavevmode
                                \begin{itemize}[$\sq$]
                                    \item Note that \begin{equation}
                                        \begin{array}{rcl}
                                            \operatorname E_x\left[\operatorname{II}_t\right]&=&\operatorname E^{(1)}_x\left[\left(\left.\mu f\right|_{E_2}-f\right)\left(\killedlocalprocess^{(1)}_{\tau_1-}\right);t\ge\tau_1\right]\\
                                            &=&\operatorname E^{(1)}_x\left[\underbrace{\int_0^t\underbrace{C^{(1)}_s\left(\left.\mu f\right|_{E_2}-f\right)\left(\killedlocalprocess^{(1)}_{s-}\right)}_{=:\:Y_s}\dif s}_{=:\:Z_t}\right]
                                        \end{array}
                                    \end{equation} for all $t\ge0$.
                                    \item $\left(C^{(1)}_t\right)_{t\ge0}$ and $[0,\infty)\ni t\mapsto \killedlocalprocess^{(1)}_{t-}$ are (right-)continuous at $0$ $\Rightarrow$ \begin{equation}
                                        \frac{Z_t}t\xrightarrow{t\to0+}Y_0\stackrel{\text{def}}=C^{(1)}_0\left(\left.\mu f\right|_{E_2}-f\right)\left(\killedlocalprocess^{(1)}_0\right).
                                    \end{equation}
                                    \item By assumption, \begin{equation}
                                        M:=\sup_{(\omega,\:t)\in\Omega\times[0,\:t_0)}C^{(1)}_t(\omega)<\infty
                                    \end{equation} for some $t_0>0$ and hence \begin{equation}
                                        \frac{|Z_t|}t\le\frac1t\int_0^t|Y_s|\dif s\le2\|f\|_\infty\frac1t\int_0^tC^{(1)}_s\dif s\le2\|f\|_\infty M
                                    \end{equation} for all $t\in[0,t_0]$.
                                    \item Lebesgue’s dominated convergence theorem $\Rightarrow$ \begin{equation}
                                        \begin{split}
                                            &\frac{\operatorname E^{(1)}_x\left[Z_t\right]}t\xrightarrow{t\to0+}\operatorname E^{(1)}_x\left[C^{(1)}_0\left(\left.\mu f\right|_{E_2}-f\right)\left(\killedlocalprocess^{(1)}_0\right)\right]\\&\;\;\;\;\;\;\;\;\;\;\;\;=\operatorname E^{(1)}_x\left[C^{(1)}_0\right]\left(\left.\mu f\right|_{E_2}-f\right)(x).
                                        \end{split}
                                    \end{equation}
                                \end{itemize}
                            \end{proof}
                        \end{listclaim}
                        \item \begin{listclaim}
                            \begin{equation}
                                \frac{\operatorname E_x\left[\operatorname{III}_t\right]}t\xrightarrow{t\to0+}0.
                            \end{equation}
                            \begin{proof}[Proof\textup:\nopunct]
                                \leavevmode
                                \begin{itemize}[$\sq$]
                                    \item \begin{listclaim}
                                        \begin{equation}
                                            \frac{\operatorname E^{(1)}_x\left[f\left(\killedlocalprocess^{(1)}_{\tau_1-}\right)-f(x);t\ge\tau_1\right]}t\xrightarrow{t\to0+}0.
                                        \end{equation}
                                        \begin{proof}[Proof\textup:\nopunct]
                                            \leavevmode
                                            \begin{itemize}[$\sq$]
                                                \item Note that \begin{equation}
                                                    \begin{split}
                                                        &\operatorname E^{(1)}_x\left[f\left(\killedlocalprocess^{(1)}_{\tau_1-}\right)-f(x);t\ge\tau_1\right]\\&\;\;\;\;\;\;\;\;\;\;\;\;=\operatorname E^{(1)}_x\left[\underbrace{\int_0^t\underbrace{C^{(1)}_s\left(f\left(\killedlocalprocess^{(1)}_{s-}\right)-f(x)\right)}_{=:\:Y_s}\dif s}_{=:\:Z_t}\right]
                                                    \end{split}
                                                \end{equation} for all $t\ge0$.
                                                \item $\left(C^{(1)}_t\right)_{t\ge0}$ and $[0,\infty)\ni t\mapsto \killedlocalprocess^{(1)}_{t-}$ are (right-)continuous at $0$ $\Rightarrow$ \begin{equation}
                                                    \frac{Z_t}t\xrightarrow{t\to0+}Y_0\stackrel{\text{def}}=C^{(1)}_0\left(f\left(\killedlocalprocess^{(1)}_0\right)-f(x)\right).
                                                \end{equation}
                                                \item Lebesgue’s dominated convergence theorem $\Rightarrow$ \begin{equation}
                                                    \frac{\operatorname E^{(1)}_x\left[Z_t\right]}t\xrightarrow{t\to0+}\operatorname E^{(1)}_x\left[\killedlocalprocess^{(1)}_0\left(f\left(\killedlocalprocess^{(1)}_0\right)-f(x)\right)\right]=0.
                                                \end{equation}
                                            \end{itemize}
                                        \end{proof}
                                    \end{listclaim}
                                    \item \begin{listclaim}
                                        \begin{equation}
                                            \frac{\operatorname E^{(1)}_x\left[T_1\left(\left.Q^{(2)}_{t-\tau_1}f\right|_{E_2}-\left.f\right|_{E_2}\right);t\ge\tau_1\right]}t\xrightarrow{t\to0+}0.
                                        \end{equation}
                                        \begin{proof}[Proof\textup:\nopunct]
                                            \leavevmode
                                            \begin{itemize}[$\sq$]
                                                \item Let \begin{equation*}
                                                    \begin{split}
                                                        &Y_{s,\:t}(\omega_1):=y_{\omega_1}(s,t)\\
                                                        &:=\begin{cases}
                                                            C^{(1)}_s(\omega_1)\left(\mu\left(\left.Q^{(2)}_{t-s}f\right|_{E_2}-\left.f\right|_{E_2}\right)\right)\left(\killedlocalprocess^{(1)}_{s-}(\omega_1)\right)&\text{, if }s\le t;\\0&\text{, otherwise}
                                                        \end{cases}
                                                    \end{split}
                                                \end{equation*} for $s,t\ge0$ and \begin{equation}
                                                    Z_t(\omega_1):=z_{\omega_1}(t):=\int_0^ty_{\omega_1}(s,t)\dif s\;\;\;\text{for }t\ge0
                                                \end{equation} for $\omega_1\in\Omega_1$.
                                                \item Note hat \begin{equation}
                                                    \operatorname E^{(1)}_x\left[T_1\left(\left.Q^{(2)}_{t-\tau_1}f\right|_{E_2}-\left.f\right|_{E_2}\right);t\ge\tau_1\right]=\operatorname E^{(1)}_x\left[Z_t\right]
                                                \end{equation} for all $t\ge0$.
                                                \item \begin{listclaim}\label{claim:2.3.2.1}
                                                    \normalfont Let $\omega_1\in\Omega_1$ $\Rightarrow$ $z_{\omega_1}$ is (right-)differentiable at $0$ with \begin{equation}
                                                        z_{\omega_1}'(0)=0.
                                                    \end{equation}
                                                    \begin{proof}[Proof\textup:\nopunct]
                                                        \leavevmode
                                                        \begin{itemize}[$\circ$]
                                                            \item \begin{listclaim}\label{claim:2.3.2.1.1}
                                                                \normalfont Let $\varepsilon>0$ $\Rightarrow$ $\exists$ $\delta>0$ with \begin{equation}
                                                                    |y_{\omega_1}(s,t)-\underbrace{y_{\omega_1}(0,0)}_{=\:0}|<\varepsilon
                                                                \end{equation} for all $0\le s\le t<\delta$.
                                                                \begin{proof}[Proof\textup:\nopunct]
                                                                    \leavevmode
                                                                    \begin{itemize}[$\sq$]
                                                                        \item $\left(C^{(1)}_t\right)_{t\ge0}$ is (right-)continuous at $0$ $\Rightarrow$ $\exists\delta_1>0$ with \begin{equation}
                                                                            \left|C^{(1)}_s(\omega_1)-C^{(1)}_0(\omega_1)\right|<1
                                                                        \end{equation} and hence \begin{equation}
                                                                            \left|C^{(1)}_s(\omega_1)\right|<1+\left|C^{(1)}_0(\omega_1)\right|=:M
                                                                        \end{equation} for all $s\in[0,\delta)$.
                                                                        \item Let $\varepsilon>0$ and \begin{equation}
                                                                            \tilde\varepsilon:=\frac\varepsilon M.
                                                                        \end{equation}
                                                                        \item $\left.f\right|_{E_2}\in C_2$ $\Rightarrow$ $\exists\delta>0$ with \begin{equation}
                                                                            \left\|\left.Q^{(2)}_tf\right|_{E_2}-\left.f\right|_{E_2}\right\|_\infty<\tilde\varepsilon\;\;\;\text{for all }t\in[0,\delta_2)
                                                                        \end{equation} and hence \begin{equation}
                                                                            \begin{split}
                                                                                &\left|y_{\omega_1}(s,t)-y_{\omega_1}(0,0)\right|\\&\;\;\;\;=\left|C^{(1)}_s(\omega_1)\left(\mu\left(\left.Q^{(2)}_{t-s}f\right|_{E_2}-\left.f\right|_{E_2}\right)\right)\left(\killedlocalprocess^{(1)}_{s-}(\omega_1)\right)\right|\\&\;\;\;\;\le\underbrace{\left|C^{(1)}_s(\omega_1)\right|}_{<\:M}\underbrace{\left\|\left.Q^{(2)}_{t-s}f\right|_{E_2}-\left.f\right|_{E_2}\right\|_\infty}_{<\:\tilde\varepsilon}<\varepsilon
                                                                            \end{split}
                                                                        \end{equation} for all $0\le s\le t<\delta_1\wedge\delta_2$.
                                                                    \end{itemize}
                                                                \end{proof}
                                                            \end{listclaim}
                                                            \item \autoref{claim:2.3.2.1.1} $\Rightarrow$ $z_{\omega_1}$ is (right-)differentiable at $0$ with \begin{equation}
                                                                z_{\omega_1}'(0)=y_{\omega_1}(0,0)=0.
                                                            \end{equation}
                                                        \end{itemize}
                                                    \end{proof}
                                                \end{listclaim}
                                                \item Lebesgue’s dominated convergence theorem $\Rightarrow$ \begin{equation}
                                                    \frac{\operatorname E^{(1)}_x\left[Z_t\right]}t\xrightarrow{t\to0+}0
                                                \end{equation} by \autoref{claim:2.3.2.1}.
                                            \end{itemize}
                                        \end{proof}
                                    \end{listclaim}
                                \end{itemize}
                            \end{proof}
                        \end{listclaim}
                    \end{itemize}
                \end{proof}
            \end{listclaim}
            \item \begin{listclaim}
                \normalfont If $x\in E_2$, then \begin{equation}
                    \frac{\operatorname E_x\left[\operatorname I_t+\operatorname{II}_t+\operatorname{III}_t\right]}t\xrightarrow{t\to0+}\left(\left.A^{(2)}_{\textnormal p}f\right|_{E_2}\right)(x).
                \end{equation}
                \begin{proof}[Proof\textup:\nopunct]
                    \leavevmode
                    \begin{itemize}[$\circ$]
                        \item Simply note that \begin{equation}
                            \operatorname E_x\left[\operatorname I_t+\operatorname{II}_t+\operatorname{III}_t\right]=\left(\left.Q^{(2)}_t\right|_{E_2}\right)(x)-f(x)
                        \end{equation} for all $t\ge0$.
                    \end{itemize}
                \end{proof}
            \end{listclaim}
        \end{itemize}
    \end{proof}
\end{theorem}

Let 
\begin{equation}
    \globalgenerator_\indexsetelement\integrand:=\regenerationdistribution_\indexsetelement\left.\integrand\right|_{\measurablespace_{\indexsetelement+1}}-\integrand\;\;\;\text{for }\integrand\in\measurablesystem_b
\end{equation} for $\indexsetelement\in\indexset_0$.

\begin{example}[killing at an exponential rate (cont.)]
    Assume \autoref{set:concatenation-of-killing-at-exponential-rate}. Let $\indexsetelement\in\indexset_0$ and $\localgenerator_\indexsetelement$ denote the pointwise generator of $\largebracket(\localsemigroup^{(\indexsetelement)}_\timepoint\largebracket)_{\timepoint\ge0}$. If $\killingrate_\indexsetelement$ is bounded and $\killingrate_\indexsetelement\circ\localprocess^{(\indexsetelement)}$ is (right-)continuous at $0$ $\probabilitymeasure^{(\indexsetelement)}_\point$-almost surely for all $\point\in\measurablespace_\indexsetelement$, then $\mathcal D\left(\localgenerator_\indexsetelement\right)\subseteq\mathcal D\left(\killedlocalgenerator_\indexsetelement\right)$ and \begin{equation}
		\killedlocalgenerator_\indexsetelement\integrand=\left(\localgenerator_\indexsetelement-\killingrate\right)\integrand\;\;\;\text{for all }\integrand\in\mathcal D\left(\localgenerator_\indexsetelement\right)
	\end{equation} by \autoref{thm:pointwise-generator-of-process-killed-by-mf} and hence the right-hand side of \eqref{eq:generator-of-concatenated-process} is equal to \begin{equation}
	    \left(\localgenerator_\indexsetelement\left.\integrand\right|_{\measurablespace_\indexsetelement}\right)(\point)+\killingrate_\indexsetelement(\point)(\globalgenerator_\indexsetelement\integrand)(\point)\;\;\;\text{for all }(\point,\integrand)\in\measurablespace_\indexsetelement\times\measurablesystem_b\text{ with }\left.\integrand\right|_{\measurablespace_\indexsetelement}\in\mathcal D(\localgenerator_\indexsetelement).
	\end{equation}
    \begin{flushright}
        $\square$
    \end{flushright}
\end{example}

\appendix
\section{Basic measure theory}

\subsection{Measurability in larger space}

Let \begin{itemize}[$\circ$]
    \item $I$ be a nonempty set;
    \item $(E_i,\mathcal E_i)$ be a measurable space for $i\in I$ with \begin{equation}
    E_i\cap E_j=\emptyset\;\;\;\text{for all }i,j\in I\text{ with }i\ne j
\end{equation} and \begin{align}
    E&:=\biguplus_{i\in I}E_i;\\\mathcal E&:=\bigvee_{i\in I}\mathcal E_i.
\end{align}
\end{itemize}

\begin{proposition}\label{prop:measurability-in-subspace}
    Let $i\in I$.
    \begin{enumerate}[(i)]
        \item $\left.\mathcal E\right|_{E_i}=\mathcal E_i$.
        \item\label{prop:measurability-in-subspace-ii} Let $(\Omega,\mathcal A)$ be a measurable space and $X:\Omega\to E_i$ $\Rightarrow$\newline $X$ is $(\mathcal A,\mathcal E_i)$ measurable iff $X$ is $(\mathcal A,\mathcal E)$-measurable.\newline In particular, \begin{equation}
            X^{-1}(\mathcal E_i)=X^{-1}(\mathcal E).
        \end{equation}
    \end{enumerate}
    \begin{proof}[Proof\textup:\nopunct]
        \leavevmode
        \begin{enumerate}[(i)]
            \item \begin{equation}
                \left.\bigcup_{j\in I}\mathcal E_j\right|_{E_i}\stackrel{E_i\:\in\:\mathcal E_i}=\left\{B\in\bigcup_{j\in I}\mathcal E_j:B\subseteq E_i\right\}\stackrel{E_i\:\cap\:E_j\:=\:\emptyset\:\text{for all }j\in I\text{ with }i\ne j}=\mathcal E_i
            \end{equation} and hence \begin{equation}
                \left.\mathcal E\right|_{E_i}\stackrel{\text{\citep[Lemma~14.20]{klenke2020probability}}}=\sigma\left(\left.\bigcup_{j\in I}\mathcal E_j\right|_{E_i}\right)=\mathcal E_i.
            \end{equation}
            \item \begin{equation}
                \mathcal E_i\stackrel{\text{(i)}}=\left.\mathcal E\right|_{E_i}\stackrel{\text{def}}=\left\{B\cap E_i:B\in\mathcal E\right\}
            \end{equation} and \begin{equation}
                \{X\in B\}\stackrel{X(\Omega)\:\subseteq\:E_i}=\{X\in B\cap E_i\}\;\;\;\text{for all }B\in\mathcal E.
            \end{equation}
        \end{enumerate}
    \end{proof}
\end{proposition}

\noindent Let \begin{itemize}[$\circ$]
    \item $(\Omega_1,\mathcal A_1)$ be a measurable space;
    \item $\Omega_2$ be a set and \begin{equation}
        \Omega:=\Omega_1\times\Omega_2.
    \end{equation}
\end{itemize}

\begin{lemma}\label{lem:embedding-of-sigma-algebra-into-larger-space}
    \begin{equation}
        \tilde A_1:=\mathcal A_1\times\Omega_2
    \end{equation} is a $\sigma$-algebra on $\Omega$.
    \begin{proof}[Proof\textup:\nopunct]
        \leavevmode
        \begin{itemize}[$\circ$]
            \item $\Omega_1\in\mathcal A_1$ $\Rightarrow$ \begin{equation}
                \Omega\in\tilde{\mathcal A}_1.
            \end{equation}
            \item \begin{listclaim}
                \normalfont Let $A\in\tilde{\mathcal A}_1$ $\Rightarrow$ \begin{equation}
                    A^c\in\tilde{\mathcal A}_1.
                \end{equation}
                \begin{proof}[Proof\textup:\nopunct]
                    \leavevmode
                    \begin{itemize}[$\circ$]
                        \item $A\in\tilde{\mathcal A}_1$ $\Rightarrow$ $\exists A_1\in\mathcal A_1$ with \begin{equation}
                            A=A_1\times\Omega_2.
                        \end{equation} $A_1^c\in\mathcal A_1$ $\Rightarrow$ \begin{equation}
                            A^c=A_1^c\times\Omega_2\in\mathcal A_1.
                        \end{equation}
                    \end{itemize}
                \end{proof}
            \end{listclaim}
            \item \begin{listclaim}
                \normalfont Let $(A^n)_{n\in\mathbb N}\subseteq\tilde{\mathcal A}_1$ $\Rightarrow$ \begin{equation}
                    \bigcup_{n\in\mathbb N}A^n\in\tilde{\mathcal A}_1.
                \end{equation}
                \begin{proof}[Proof\textup:\nopunct]
                    \leavevmode
                    \begin{itemize}[$\circ$]
                        \item $(A^n)_{n\in\mathbb N}\subseteq\tilde{\mathcal A}_1$ $\Rightarrow$ $\exists(A^n_1)_{n\in\mathbb N}\subseteq\mathcal A_1$ with \begin{equation}
                            A^n=A^n_1\times\Omega_2\;\;\;\text{for all }n\in\mathbb N.
                        \end{equation} $\bigcup_{n\in\mathbb N}A^n_1\in\mathcal A_1$ $\Rightarrow$ \begin{equation}
                            \bigcup_{n\in\mathbb N}A^n=\bigcup_{n\in\mathbb N}A^n_1\times\Omega_2\in\tilde{\mathcal A}_1.
                        \end{equation}
                    \end{itemize}
                \end{proof}
            \end{listclaim}
        \end{itemize}
    \end{proof}
\end{lemma}

\noindent Let \begin{itemize}[$\circ$]
    \item $(I,\le)$ be a partially ordered set;
    \item $(\mathcal F^1_t)_{t\in I}$ be aq filtration on $(\Omega_1,\mathcal A_1)$;
    \item $\mathcal A_2$ be a $\sigma$-algebra on $\Omega_2$ and \begin{equation}
        \mathcal A:=\mathcal A_1\otimes\mathcal A_2.
    \end{equation}
\end{itemize}

\begin{lemma}
    \begin{equation}
        \tilde{\mathcal F}^1_t:=\mathcal F^1_t\times\Omega_2\;\;\;\text{for }t\in I
    \end{equation} is a filtration on $(\Omega,\mathcal A)$.
    \begin{proof}[Proof\textup:\nopunct]
        \leavevmode
        \begin{itemize}[$\circ$]
            \item \autoref{lem:embedding-of-sigma-algebra-into-larger-space} $\Rightarrow$ $\tilde{\mathcal F}^1_t$ is a $\sigma$-algebra on $\Omega$ for all $t\in I$.
            \item $\Omega_2\in\mathcal A_2$ $\Rightarrow$ \begin{equation}
                \tilde{\mathcal F}^1_t\subseteq\mathcal A\;\;\;\text{for all }t\in I.
            \end{equation}
            \item $(\mathcal F^1_t)_{t\in I}$ is a filtration on $(\Omega_1,\mathcal A_1)$ $\Rightarrow$ \begin{equation}
                \tilde{\mathcal F}^1_s\subseteq\tilde{\mathcal F}^1_t\;\;\;\text{for all }s,t\in I\text{ with }s\le t.
            \end{equation}
        \end{itemize}
    \end{proof}
\end{lemma}

\noindent Let $\pi_i$ denote the projection from $\Omega$ onto the $i$th coordinate.

\begin{lemma}\label{lem:measurability-on-larger-space}
    Let $(E,\mathcal E)$ be a measurable space, $X_1:\Omega_1\to E$ and \begin{equation}
        \tilde X_1:=X_1\circ\pi_1.
    \end{equation} Then, \begin{equation}
        \underbrace{\tilde X_1^{-1}(\mathcal E)}_{\stackrel{\text{def}}=\:\sigma\left(\tilde X_1\right)}=\underbrace{X_1^{-1}(\mathcal E)}_{\stackrel{\text{def}}=\:\sigma(X_1)}\times\Omega_2.
    \end{equation}
    \begin{proof}[Proof\textup:\nopunct]
        Trivial.
    \end{proof}
\end{lemma}

\noindent An immediate consequence of \autoref{lem:measurability-on-larger-space} is that the canonical extension of an $({\mathcal F}^1_t)_{t\in I}$-stopping time on $(\Omega_1,\mathcal A_1)$ is an $\left(\tilde{\mathcal F}^1_t\right)_{t\in I}$-stopping time on $(\Omega,\mathcal A)$:

\begin{lemma}\label{lem:stopping-time-on-larger-space}
    Let $\sigma_1$ be an $(\mathcal F^1_t)_{t\in I}$-stopping time on $(\Omega_1,\mathcal A_1)$ $\Rightarrow$ \begin{equation}
        \tilde\sigma_1:=\sigma_1\circ\pi_1
    \end{equation} is an $\left(\tilde{\mathcal F}^1_t\right)_{t\in I}$-stopping time on $(\Omega,\mathcal A)$.
    \begin{proof}[Proof\textup:\nopunct]
        \leavevmode
        \begin{itemize}[$\circ$]
            \item $\sigma_1$ is an $(\mathcal F^1_t)_{t\in I}$-stopping time on $(\Omega,\mathcal A)$ $\Rightarrow$ \begin{equation}
                \left\{\tilde\sigma_1\le t\right\}=\underbrace{\{\sigma_1\le t\}}_{\in\:\mathcal F^1_t}\times\Omega_2\in\tilde{\mathcal F}^1_t\;\;\;\text{for all }t\in I.
            \end{equation}
        \end{itemize}
    \end{proof}
\end{lemma}

\section{Topological preliminaries}

\subsection{Path properties}

\subsubsection{Traps}

Let \begin{itemize}[$\circ$]
    \item $E$ be a nonempty set;
    \item $x:[0,\infty)\to E$.
\end{itemize}

\begin{definition}\label{def:trap}
    \leavevmode
    \begin{enumerate}[(i)]
        \item $B\subseteq E$ is called \textbf{absorbing} (or \textbf{trap}) \textbf{for} $\bm x$ $:\Leftrightarrow$ \begin{equation}
            x(s)\in B\Rightarrow\forall t\ge s:x(t)\in B
        \end{equation} for all $s\ge0$.
        \item $x_0$ is called \textbf{absorbing} (or \textbf{trap}) \textbf{for} $\bm x$ $:\Leftrightarrow$ $\{x_0\}$ is absorbing for $x$.
    \end{enumerate}
    \begin{flushright}
        $\square$
    \end{flushright}
\end{definition}

\subsubsection{Hitting time of a closed set}

Let \begin{itemize}[$\circ$]
    \item $E$ be a topological space and \begin{equation}
        \mathcal N(x):=\left\{N\subseteq E:N\text{ is a neighborhood of }x\right\}\;\;\;\text{for }x\in E;
    \end{equation}
    \item $x:[0,\infty)\to E$.
\end{itemize}

\begin{definition}
    \leavevmode
    \begin{enumerate}[(i)]
        \item Let $t\ge0$ $\Rightarrow$ $x$ is called \textbf{right-continuous at} $\bm t$ $:\Leftrightarrow$ \begin{equation}
            \forall N\in\mathcal N\left(x(t)\right):\exists\delta>0:x\left((t,t+\delta)\right)\subseteq N.
        \end{equation}
        \item $x$ is called \textbf{right-continuous} $:\Leftrightarrow$ $x$ is right-continuous at $t$ for all $t\ge0$.
    \end{enumerate}
    \begin{flushright}
        $\square$
    \end{flushright}
\end{definition}

\noindent Assume $x$ is right-continuous. Let \begin{itemize}[$\circ$]
    \item $B\subseteq E$ be closed and \begin{equation}
        I:=\left\{t\ge0:x(t)\in B\right\};
    \end{equation}
    \item $\tau:=\inf I$.
\end{itemize}

\begin{lemma}[hitting time of closed set]\label{lem:hitting-time-of-closed-set}~
    \bigbreak\noindent Either \begin{enumerate}[(i)]
        \item $\tau\in I$; or
        \item $I=\emptyset$ and hence $\tau=\infty$.
    \end{enumerate}
    \begin{proof}[Proof\textup:\nopunct]
        \leavevmode
        \begin{itemize}[$\circ$]
            \item If $I=\emptyset$, then \begin{equation}
                \tau=\inf\emptyset\stackrel{\text{def}}=\infty.
            \end{equation}
            \item Assume $I\ne\emptyset$ $\Rightarrow$ \begin{equation}
                \tau<\infty.
            \end{equation}
            \item Assume $\tau\notin I$ $\Rightarrow$ \begin{equation}
                x(\tau)\in B^c\in\tau.
            \end{equation}
            \item $x$ is right-continuous at $\tau$ $\Rightarrow$ $\exists\delta>0$ with \begin{equation}
                x\left([\tau,\tau+\delta)\right)\subseteq B^c
            \end{equation} and hence \begin{equation}
                [\tau,\tau+\delta)\cap I=\emptyset.
            \end{equation}
            \item $\tau\stackrel{\text{def}}=\inf I$ is a lower bound for $I$ $\Rightarrow$ \begin{equation}
                I\subseteq[\tau,\infty)
            \end{equation} and hence \begin{equation}
                \tau\ge\tau+\delta;
            \end{equation} which is impossible.
        \end{itemize}
    \end{proof}
\end{lemma}

\begin{lemma}
    Let $t\ge0$ $\Rightarrow$ \begin{equation}
        \tau\le t\Leftrightarrow\exists s\in I:s\le t.
    \end{equation}
    \begin{proof}[Proof\textup:\nopunct]
        \leavevmode
        \begin{itemize}[$\circ$]
            \item \begin{listclaim}
                "$\Rightarrow$"
                \begin{proof}[Proof\textup:\nopunct]
                    \leavevmode
                    \begin{itemize}[$\circ$]
                        \item Assume the contrary $\Rightarrow$ \begin{equation}
                            \forall s\in I:s>t
                        \end{equation} and hence $t$ is a lower bound for $I$ $\Rightarrow$ \begin{equation}
                            t\le\tau
                        \end{equation} by definition of the infimum $\Rightarrow$ \begin{equation}
                            \tau=t
                        \end{equation} and hence \begin{equation}
                            \tau\not\in I.
                        \end{equation}
                        \item \autoref{lem:hitting-time-of-closed-set} $\Rightarrow$ $I=\emptyset$ and hence $\tau=\infty$; in contradiction to $\tau\le t<\infty$.
                    \end{itemize}
                \end{proof}
            \end{listclaim}
            \item \begin{listclaim}
                "$\Leftarrow$"
                \begin{proof}[Proof\textup:\nopunct]
                    \leavevmode
                    \begin{itemize}[$\circ$]
                        \item Let $s\in I$ with $s\le t$.
                        \item $\tau$ is a lower bound for $I$ $\Rightarrow$ $\tau\le s\le t$.
                    \end{itemize}
                \end{proof}
            \end{listclaim}
        \end{itemize}
    \end{proof}
\end{lemma}

\subsection{One-point extension}
\subsubsection{Measurable one-point extension}\label{sec:one-point-extension-measurable-one-point-extension}

Let $(E,\mathcal E)$ be a measurable space. We add a point $\Delta\not\in E$ which does not already belong to $E$ and consider the smallest $\sigma$-algebra on \begin{equation}
    E^\ast:=E\uplus\{\Delta\}
\end{equation} which does contain $\mathcal E$. By definition, this $\sigma$-algebra is given by \begin{equation}
    \mathcal E^\ast:=\sigma_{E^\ast}(\mathcal E).
\end{equation} Each set in $\mathcal E^\ast$ is either in $\mathcal E$ or the (disjoint) union of a set in $\mathcal E$ and the singleton set $\{\Delta\}$:

\begin{proposition}[$\sigma$-algebra of one-point extension]\label{prop:sigma-algebra-of-one-point-extension}
    \begin{equation}\label{eq:sigma-algebra-of-one-point-extension}
        \mathcal E^\ast=\underbrace{\sigma(\{B\uplus\{\Delta\}:B\in\mathcal E\})}_{=:\:\mathcal B_1}=\underbrace{\mathcal E\uplus\{B\uplus\{\Delta\}:B\in\mathcal E\}}_{=:\:\mathcal B_2}
    \end{equation}
    \begin{proof}[Proof\textup:\nopunct]
        \leavevmode
        \begin{itemize}[$\circ$]
            \item \begin{listclaim}\label{prop:sigma-algebra-of-one-point-extension-claim-1}
                $\mathcal E^\ast\subseteq\mathcal B_1$
                \begin{proof}[Proof\textup:\nopunct]
                    \leavevmode
                    \begin{itemize}[$\circ$]
                        \item \begin{listclaim}\label{prop:sigma-algebra-of-one-point-extension-claim-1.1}
                            $\{\Delta\}\in\mathcal B_1$
                            \begin{proof}[Proof\textup:\nopunct]
                                \leavevmode
                                \begin{itemize}[$\sq$]
                                    \item $\emptyset\in\mathcal E$ $\Rightarrow$ \begin{equation}
                                    \{\Delta\}\in\{B\uplus\{\Delta\}:B\in\mathcal E\}\subseteq\sigma(\{B\uplus\{\Delta\}:B\in\mathcal E\})=\mathcal B_1.
                                \end{equation}
                                \end{itemize}
                            \end{proof}
                        \end{listclaim}
                        \item \begin{listclaim}\label{prop:sigma-algebra-of-one-point-extension-claim-1.2}
                            $\mathcal E\subseteq\mathcal B_1$
                            \begin{proof}[Proof\textup:\nopunct]
                                \leavevmode
                                \begin{itemize}[$\sq$]
                                    \item Let $B\in\mathcal E$ $\Rightarrow$ \begin{equation}
                                    B\uplus\{\Delta\}\in\mathcal B_1
                                \end{equation} and hence \begin{equation}
                                    B=(B\uplus\{\Delta\})\setminus\{\Delta\}\in\mathcal B_1
                                \end{equation} by \autoref{prop:sigma-algebra-of-one-point-extension-claim-1.1}.
                                \end{itemize}
                            \end{proof}
                        \end{listclaim}
                        \item \autoref{prop:sigma-algebra-of-one-point-extension-claim-1.2} $\Rightarrow$ \begin{equation}
                            \mathcal E^\ast\stackrel{\text{def}}=\sigma_{E^\ast}(\mathcal E)\subseteq\mathcal B_1
                        \end{equation} by monotonicity of $\sigma_{E^\ast}$.
                    \end{itemize}
                \end{proof}
            \end{listclaim}
            \item \begin{listclaim}\label{prop:sigma-algebra-of-one-point-extension-claim-2}
                $\left\{B\uplus\{\Delta\}:B\in\mathcal E\right\}\subseteq\mathcal E^\ast$
                \begin{proof}[Proof\textup:\nopunct]
                    \leavevmode
                    \begin{itemize}[$\circ$]
                        \item \begin{listclaim}\label{prop:sigma-algebra-of-one-point-extension-claim-2.1}
                            $\{\Delta\}\in\mathcal E^\ast$
                            \begin{proof}[Proof\textup:\nopunct]
                                \leavevmode
                                \begin{itemize}[$\sq$]
                                    \item $E\in\mathcal E\subseteq\sigma_{E^\ast}(\mathcal E)\stackrel{\text{def}}=\mathcal E^\ast$ and $E^\ast\in\mathcal E^\ast$ $\Rightarrow$ \begin{equation}
                                    \{\Delta\}=E^\ast\setminus E\in\mathcal E^\ast.
                                \end{equation}
                                \end{itemize}
                            \end{proof}
                        \end{listclaim}
                        \item Let $B\in\mathcal E$. $\mathcal E\subseteq\sigma_{E^\ast}(\mathcal E)\stackrel{\text{def}}=\mathcal E^\ast$ and \autoref{prop:sigma-algebra-of-one-point-extension-claim-2.1} $\Rightarrow$ \begin{equation}
                            B\uplus\{\Delta\}\in\mathcal E^\ast.
                        \end{equation}
                    \end{itemize}
                \end{proof}
            \end{listclaim}
            \item \autoref{prop:sigma-algebra-of-one-point-extension-claim-1} and \autoref{prop:sigma-algebra-of-one-point-extension-claim-2} $\Rightarrow$ \begin{equation}
                \mathcal E^\ast=\mathcal B_1.
            \end{equation}
            \item \begin{listclaim}\label{prop:sigma-algebra-of-one-point-extension-claim-3}
                $\mathcal B_2\subseteq\mathcal E^\ast$
                \begin{proof}[Proof\textup:\nopunct]
                    \leavevmode
                    \begin{itemize}[$\circ$]
                        \item $\mathcal E\subseteq\sigma_{E^\ast}(\mathcal E)\stackrel{\text{def}}=\mathcal E^\ast$.
                        \item Let $B\in\mathcal E$. $\mathcal E\subseteq\mathcal E^\ast$ and \autoref{prop:sigma-algebra-of-one-point-extension-claim-2.1} $\Rightarrow$ \begin{equation}
                            B\uplus\{\Delta\}\in\mathcal E^\ast
                        \end{equation}
                    \end{itemize}
                \end{proof}
            \end{listclaim}
            \item \begin{listclaim}\label{prop:sigma-algebra-of-one-point-extension-claim-4}
                $\mathcal B_2$ is a $\sigma$-algebra on $E^\ast$.
                \begin{proof}[Proof\textup:\nopunct]
                    \leavevmode
                    \begin{itemize}[$\circ$]
                        \item Trivial.
                    \end{itemize}
                \end{proof}
            \end{listclaim}
            \item \begin{listclaim}\label{prop:sigma-algebra-of-one-point-extension-claim-5}
                $\mathcal E^\ast\subseteq\mathcal B_2$.
                \begin{proof}[Proof\textup:\nopunct]
                    \leavevmode
                    \begin{itemize}[$\circ$]
                        \item \autoref{prop:sigma-algebra-of-one-point-extension-claim-1} and \autoref{prop:sigma-algebra-of-one-point-extension-claim-4} $\Rightarrow$ \begin{equation}
                            \mathcal E^\ast\subseteq\mathcal B_1\subseteq\sigma(\mathcal E\uplus\{B\uplus\{Delta\}:B\in\mathcal E\})=\mathcal B_2.
                        \end{equation}
                    \end{itemize}
                \end{proof}
            \end{listclaim}
            \item \autoref{prop:sigma-algebra-of-one-point-extension-claim-3} and \autoref{prop:sigma-algebra-of-one-point-extension-claim-5} $\Rightarrow$ \begin{equation}
                \mathcal E^\ast=\mathcal B_2.
            \end{equation}
        \end{itemize}
    \end{proof}
\end{proposition}

\begin{proposition}[trace of one-point extension is original $\sigma$-algebra]\label{prop:trace-of-one-point-extension-is-original-sigma-algebra}
    \begin{equation}
        \left.\mathcal E^\ast\right|_E=\mathcal E
    \end{equation}
    \begin{proof}[Proof\textup:\nopunct]
        \begin{equation}
            \begin{array}{rcl}
                \left.\mathcal E^\ast\right|_E&\stackrel{\eqref{eq:sigma-algebra-of-one-point-extension}}=&\left\{B^\ast\cap E:B^\ast\in\mathcal E\uplus\left\{B\uplus\left\{\Delta\right\}:B\in\mathcal E\right\}\right\}\\&=&\{\underbrace{B^\ast\cap E}_{=\:B^\ast}:B^\ast\in\mathcal E\}=\mathcal E.
            \end{array}
        \end{equation}        
    \end{proof}
\end{proposition}

\begin{definition}\label{def:one-point-extension-of-measurable-space}
    $(E^\ast,\mathcal E^\ast)$ is called \textbf{measurable one-point extension of} $\bm{(E,\mathcal E)}$\newline \textbf{by} $\bm\Delta$.
    \begin{flushright}
        $\square$
    \end{flushright}
\end{definition}

\begin{proposition}\label{prop:measurability-of-map-into-one-point-extension}
    Let $(\Omega,\mathcal A)$ be a measurable space and $X:\Omega\to E^\ast$ $\Rightarrow$\newline $X$ is $(\mathcal A,\mathcal E^\ast)$-measurable $\Leftrightarrow$ \begin{equation}
        X^{-1}(\mathcal E)\subseteq\mathcal A
    \end{equation} and \begin{equation}
        X^{-1}(\{\Delta\})\in\mathcal A.
    \end{equation}
    \begin{proof}[Proof\textup:\nopunct]
        \leavevmode
        \begin{itemize}[$\circ$]
            \item Note that \begin{equation}
                X^{-1}(B\uplus\{\Delta\})=X^{-1}(B)\uplus X^{-1}(\{\Delta\})\;\;\;\text{for all }B\in\mathcal E
            \end{equation} and hence \begin{equation}
                \begin{split}
                    X^{-1}(\mathcal E^\ast)&=\{X^{-1}(B):B\in\mathcal E\}\cup\{X^{-1}(B\uplus\{\Delta\}):B\in\mathcal E\}\\&=X^{-1}(\mathcal E)\cup\{X^{-1}(B)\uplus X^{-1}(\{\Delta\}):B\in\mathcal E).
                \end{split}
            \end{equation}
            \item \begin{listclaim}
                \normalfont "$\Rightarrow$"
                \begin{proof}[Proof\textup:\nopunct]
                    \leavevmode
                    \begin{itemize}[$\circ$]
                        \item $X^{-1}(\mathcal E^\ast)\subseteq\mathcal A$ $\Rightarrow$ \begin{equation}
                            X^{-1}(\mathcal E)\subseteq\mathcal A
                        \end{equation} and \begin{equation}
                            X^{-1}(\{\Delta\})=X^{-1}(\emptyset)\uplus X^{-1}(\{\Delta\})\in\mathcal A.
                        \end{equation}
                    \end{itemize}
                \end{proof}
            \end{listclaim}
            \item \begin{listclaim}
                \normalfont "$\Leftarrow$"
                \begin{proof}[Proof\textup:\nopunct]
                    \leavevmode
                    \begin{itemize}[$\circ$]
                        \item $X^{-1}(\mathcal E)\subseteq\mathcal A$ and $X^{-1}(\{\Delta\})\in\mathcal A$ $\Rightarrow$ \begin{equation}
                            \underbrace{X^{-1}(B)}_{\in\:\mathcal A}\uplus \underbrace{X^{-1}(\{\Delta\})}_{\in\:\mathcal A}\in\mathcal A\;\;\;\text{for all }B\in\mathcal E.
                        \end{equation}
                    \end{itemize}
                \end{proof}
            \end{listclaim}
        \end{itemize}
    \end{proof}
\end{proposition}

\begin{proposition}[function extension to one-point extension]\label{prop:function-extension-to-one-point-extension}~
    \bigbreak\noindent Let $f:E\to\mathbb R$ be $\mathcal E$-measurable $\Rightarrow$ \begin{equation}
        f^\ast(x):=\left.\begin{cases}f(x)&\text{, if }x\in E;\\0&\text{, if }x=\Delta\end{cases}\right\}\;\;\;\text{for }x\in E^\ast
    \end{equation} is $\mathcal E^\ast$-measurable.
    \begin{proof}[Proof\textup:\nopunct]
        \leavevmode
        \begin{itemize}[$\circ$]
            \item \autoref{prop:trace-of-one-point-extension-is-original-sigma-algebra} $\Rightarrow$ \begin{equation}
                \left.\mathcal E^\ast\right|_E=\mathcal E.
            \end{equation} and hence \begin{equation}
                \left.f^\ast\right|_E=f
            \end{equation} is $\left.\mathcal E^\ast\right|_E$-measurable.
            \item $\left.f^\ast\right|_{\{\Delta\}}=0$ is constant and hence $\left.\mathcal E^\ast\right|_{\{\Delta\}}$-measurable.
        \end{itemize}
    \end{proof}
\end{proposition}

\begin{proposition}\label{prop:restriction-of-function-on-one-point-extension}
    Let $f^\ast:E^\ast\to\mathbb R$ $\Rightarrow$ \begin{enumerate}[(i)]
        \item $f^\ast$ is $\mathcal E^\ast$-measurable;
        \item $\left.f^\ast\right|_E$ is $\mathcal E$-measurable
    \end{enumerate} are equivalent.
    \begin{proof}[Proof\textup:\nopunct]
        \leavevmode
        \begin{itemize}[$\circ$]
            \item Assume $\left.f^\ast\right|_E$ is $\mathcal E$-measurable.
            \item \autoref{prop:trace-of-one-point-extension-is-original-sigma-algebra} $\Rightarrow$ \begin{equation}
                \left.\mathcal E^\ast\right|_E=\mathcal E
            \end{equation} and hence $\left.f^\ast\right|_E$ is $\left.\mathcal E^\ast\right|_E$-measurable.
            \item $\left.f^\ast\right|_{\{\Delta\}}$ is constant and hence $\left.\mathcal E^\ast\right|_{\{\Delta\}}$-measurable.
        \end{itemize}
    \end{proof}
\end{proposition}

\subsubsection{Alexandroff one-point extension}

The Alexandroff one-point extension is a method to construct a compact topology from a given (not necessarily compact) topological space by adjoining a single point. The $\sigma$-algebra generated by this topology will coincide with the one from the measurable one-point extension. 

\bigbreak\noindent Let $(E,\tau)$ be a topological space. Once again, we add a point $\Delta\not\in E$ which does not already belong to $E$ and consider on \begin{equation}
    E^\ast:=E\uplus\{\Delta\}
\end{equation} the set system\footnote{Remember that if $\tau$ is Hausdorff, then every $\tau$-compact subset of $E$ is $\tau$-closed.} \begin{equation}
    \tau^\ast:=\tau\uplus\underbrace{\left\{E^\ast\setminus B:B\subseteq E\text{ is }\tau\text{-closed and }\tau\text{-compact}\right\}}_{=:\:\mathcal C}.
\end{equation}

\begin{remark}
    $\tau^\ast$ is a topology on $E^\ast$.
    \begin{flushright}
        $\square$
    \end{flushright}
\end{remark}

\begin{proposition}[subspace topology in one-point extension is original topology]\label{prop:subspace-topology-in-one-point-extension-is-original-topology}
    \begin{equation}
        \left.\tau^\ast\right|_E=\tau
    \end{equation}
    \begin{proof}[Proof\textup:\nopunct]
        \begin{equation}
            (E^\ast\setminus B)\cap E=\underbrace{(E^\ast\cap E)}_{=\:E}\setminus\underbrace{(B\cap E)}_{=\:B}=E\setminus B\;\;\;\text{for all }B\subseteq E
        \end{equation} and hence \begin{equation}
            \left.\tau^\ast\right|_E\stackrel{\text{def}}=\{B^\ast\cap E:B^\ast\in\tau\uplus\mathcal C\}=\{\underbrace{B^\ast\cap E}_{=\:B^\ast}:B^\ast\in\tau\}=\tau.
        \end{equation}
    \end{proof}
\end{proposition}

\begin{remark}[singleton formed by point at infinity is closed in one-point extension]~
    \bigbreak\noindent $E\in\tau\subseteq\tau^\ast$ $\Rightarrow$ \begin{equation}
        \{\Delta\}=E^\ast\setminus E
    \end{equation} is $\tau^\ast$-closed.
    \begin{flushright}
        $\square$
    \end{flushright}
\end{remark}

\begin{remark}[inclusion into one-point extension is topological embedding]\label{prop:inclusion-into-one-point-extension-is-topological-embedding}~
    \bigbreak\noindent Let $\iota$ denote the inclusion from $E$ into $E^\ast$ $\Rightarrow$ \begin{enumerate}[(i)]
        \item $\iota$ is $(\tau,\tau^\ast)$-continuous;
        \item $\iota$ is $(\tau,\tau^\ast)$-open.
    \end{enumerate}
    \begin{flushright}
        $\square$
    \end{flushright}
\end{remark}

\begin{lemma}\label{lem:compact-implies-compact-in-one-point-extension}
    Let $B\subseteq E$ be $\tau$-compact $\Rightarrow$ $B$ is $\tau^\ast$-compact.
    \begin{proof}[Proof\textup:\nopunct]
        \leavevmode
        \begin{itemize}[$\circ$]
            \item \autoref{prop:inclusion-into-one-point-extension-is-topological-embedding} $\Rightarrow$ $\left.\iota\right|_B$ is $(\left.\tau\right|_B,\tau^\ast)$-continuous;\newline since the restriction of a continuous function is continuous.
            \item $B$ is $\left.\tau\right|_B$-compact $\Rightarrow$ $\left.\iota\right|_B(B)=B$ is $\left.\tau^\ast\right|_B$-compact;\newline since the continuous image of a compact space is compact.
        \end{itemize}
    \end{proof}
\end{lemma}

\begin{proposition}[one-point extension is compact]~
    \bigbreak\noindent $E^\ast$ is $\tau^\ast$-compact.
    \begin{proof}[Proof\textup:\nopunct]
        \leavevmode
        \begin{itemize}[$\circ$]
            \item Let $B^\ast\subseteq2^{E^\ast}$ be an $\tau^\ast$-open cover of $E^\ast$ $\Rightarrow$\newline $\exists B^\ast\in\mathcal B^\ast$ with $\{\Delta\}\subseteq B^\ast$ $\Rightarrow$ $B^\ast\in\mathcal C$ $\Rightarrow$ \begin{equation}
                B^\ast=E^\ast\setminus B
            \end{equation} for some $\tau$-compact $B\subseteq E$.
            \item \autoref{lem:compact-implies-compact-in-one-point-extension} $B$ is $\tau^\ast$-compact $\Rightarrow$ $B$ has a finite cover $B\subseteq\mathcal B^\ast$.
            \item By construction, \begin{equation}
                E^\ast=B\cup\underbrace{(E^\ast\setminus B)}_{=\:B^\ast}\subseteq\underbrace{\bigcup\mathcal B}_{\supseteq\:B}\cup B^\ast
            \end{equation} and hence $\mathcal B\cup\{B^\ast\}\subseteq\mathcal B^\ast$ is a finite cover of $E^\ast$.
        \end{itemize}
    \end{proof}
\end{proposition}

\begin{proposition}[singleton formed by point at infinity is open in one-point extension of compact space]~
    \bigbreak\noindent If $E$ is $\tau$-compact, then $\{\Delta\}\in\tau^\ast$.
    \begin{proof}[Proof\textup:\nopunct]
        $E$ is $\tau$-closed and $\tau$-compact $\Rightarrow$ \begin{equation}
            \{\Delta\}=E^\ast\setminus E\in\sigma(\tau^\ast)\in\mathcal C\subseteq\tau^\ast.
        \end{equation}
    \end{proof}
\end{proposition}

\begin{proposition}[one-point extension of locally compact space is Hausdorff iff original space is]~
    \bigbreak\noindent If $E$ is locally $\tau$-compact, then $E$ is $\tau$-Hausdorff iff $E^\ast$ is $\tau^\ast$-Hausdorff.
    \begin{proof}[Proof\textup:\nopunct]
        \author[Paragraph~29]{munkres}
    \end{proof}
\end{proposition}

\begin{definition}\label{def:alexanroff-one-point-extension}
    $(E^\ast,\tau^\ast)$ is called \textbf{Alexandroff one-point extension of} $\bm{(E,\tau)}$\newline \textbf{by} $\bm\Delta$.
    \begin{flushright}
        $\square$
    \end{flushright}
\end{definition}

\bigbreak\noindent We now turn to the $\sigma$-algebra on $E^\ast$ generated by $\tau^\ast$:

\begin{proposition}[Borel $\sigma$-algebra of Alexandroff one-point extension]\label{prop:borel-sigma-algebra-of-one-point-extension}
    \begin{equation}
        \sigma(\tau^\ast)=\underbrace{\sigma\left(\left\{B\uplus\{\Delta\}:B\in\tau\right\}\right)}_{=:\:\mathcal B_1}=\underbrace{\sigma(\tau)\uplus\left\{B\uplus\{\Delta\}:B\in\sigma(\tau)\right\}}_{=:\:\mathcal B_2}.
    \end{equation}
    \begin{proof}[Proof\textup:\nopunct]
        \leavevmode
        \begin{itemize}[$\circ$]
            \item \begin{listclaim}\label{prop:borel-sigma-algebra-of-one-point-extension-claim-1}
                $\tau^\ast\subseteq\{B\uplus\{\Delta\}:B\in\tau\}$
                \begin{proof}[Proof\textup:\nopunct]
                    \leavevmode
                    \begin{itemize}[$\circ$]
                        \item Let $B^\ast\in\tau^\ast$.
                        \item Assume $B^\ast\not\in\tau$ $\Rightarrow$ $B^\ast\in\mathcal C$ $\Rightarrow$ \begin{equation}
                            B^\ast=E^\ast\setminus B=\underbrace{(E\setminus B)}_{\in\:\tau}\uplus\{\Delta\}.
                        \end{equation} for some $\tau$-closed $b\subseteq E$:
                    \end{itemize}
                \end{proof}
            \end{listclaim}
            \item \begin{listclaim}\label{prop:borel-sigma-algebra-of-one-point-extension-claim-2}
                $\{B\uplus\{\Delta\}:B\in\tau\}\subseteq\sigma(\tau^\ast)$
                \begin{proof}[Proof\textup:\nopunct]
                    \leavevmode
                    \begin{itemize}[$\circ$]
                        \item $E\in\tau\subseteq\tau^\ast\subseteq\sigma(\tau^\ast)$ and $E^\ast\in\sigma(\tau^\ast)$ $\Rightarrow$ \begin{equation}\label{eq:borel-sigma-algebra-of-one-point-extension-2.1}
                            \{\Delta\}=E^\ast\setminus E\in\sigma(\tau^\ast).
                        \end{equation}
                        \item Let $B\in\tau$.
                        \item $B\in\tau\subseteq\tau^\ast\subseteq\sigma(\tau^\ast)$ and \eqref{eq:borel-sigma-algebra-of-one-point-extension-2.1} $\Rightarrow$ \begin{equation}
                            B\uplus\{\Delta\}\in\sigma(\tau^\ast.
                        \end{equation}
                    \end{itemize}
                \end{proof}
            \end{listclaim}
            \item \autoref{prop:borel-sigma-algebra-of-one-point-extension-claim-1} and \autoref{prop:borel-sigma-algebra-of-one-point-extension-claim-2} $\Rightarrow$ \begin{equation}
                \sigma(\tau^\ast)=\mathcal B_1.
            \end{equation}
            \item \begin{listclaim}\label{prop:borel-sigma-algebra-of-one-point-extension-claim-3}
                $\mathcal B_2\subseteq\sigma(\tau^\ast)$
                \begin{proof}[Proof\textup:\nopunct]
                    \leavevmode
                    \begin{itemize}[$\circ$]
                        \item $\tau\subseteq\tau^\ast$ $\Rightarrow$ $\sigma(\tau)\subseteq\sigma(\tau^\ast)$.
                        \item Let $B\in\sigma(\tau)$ $\Rightarrow$ \begin{equation}
                            B\uplus\{\Delta\}\in\sigma(\tau^\ast)
                        \end{equation} by \eqref{eq:borel-sigma-algebra-of-one-point-extension-2.1}.
                    \end{itemize}
                \end{proof}
            \end{listclaim}
            \item \begin{listclaim}\label{prop:borel-sigma-algebra-of-one-point-extension-claim-4}
                $\mathcal B_2$ is a $\sigma$-algebra on $E^\ast$
                \begin{proof}[Proof\textup:\nopunct]
                    \leavevmode
                    \begin{itemize}[$\circ$]
                        \item Trivial.
                    \end{itemize}
                \end{proof}
            \end{listclaim}
            \item \begin{listclaim}\label{prop:borel-sigma-algebra-of-one-point-extension-claim-5}
                $\sigma(\tau^\ast)\subseteq\mathcal B_2$
                \begin{proof}[Proof\textup:\nopunct]
                    \leavevmode
                    \begin{itemize}[$\circ$]
                        \item \autoref{prop:borel-sigma-algebra-of-one-point-extension-claim-1} $\Rightarrow$ \begin{equation}
                            \tau^\ast\subseteq\{B\uplus\{\Delta\}:B\in\tau\}\subseteq\{B\uplus\{\Delta\}:B\in\sigma(\tau)\}\subseteq\mathcal B_2
                        \end{equation} and hence \begin{equation}
                            \sigma(\tau^\ast)\subseteq\sigma(\mathcal B_2)=\mathcal B_2
                        \end{equation} by \autoref{prop:borel-sigma-algebra-of-one-point-extension-claim-4}.
                    \end{itemize}
                \end{proof}
            \end{listclaim}
            \item \autoref{prop:borel-sigma-algebra-of-one-point-extension-claim-3} and \autoref{prop:borel-sigma-algebra-of-one-point-extension-claim-5} $\Rightarrow$ \begin{equation}
                \sigma(\tau^\ast)=\mathcal B_2.
            \end{equation}
        \end{itemize}
    \end{proof}
\end{proposition}

\bigbreak\noindent Similar to \autoref{prop:function-extension-to-one-point-extension}, the canonical extension of a continuous function on $(E,\tau)$ to $(E^\ast,\tau^\ast)$ is still continuous:

\begin{proposition}[function extension to Alexandroff one-point extension]\label{prop:function-extension-to-alexandroff-one-point-extension}~
    \bigbreak\noindent Let $f:E\to\mathbb R$ be $\tau$-continuous $\Rightarrow$ \begin{equation}
        f^\ast(x):=\left.\begin{cases}f(x)&\text{, if }x\in E;\\0&\text{, if }x=\Delta\end{cases}\right\}\;\;\;\text{for }x\in E^\ast
    \end{equation} is $\tau^\ast$-continuous.
    \begin{proof}[Proof\textup:\nopunct]
        \leavevmode
        \begin{itemize}[$\circ$]
            \item \autoref{prop:subspace-topology-in-one-point-extension-is-original-topology} $\Rightarrow$ \begin{equation}
                \left.\tau^\ast\right|_E=\tau
            \end{equation} and hence \begin{equation}
                \left.f^\ast\right|_E=f
            \end{equation} is $\left.\tau^\ast\right|_E$-continuous.
            \item $\left.f^\ast\right|_{\{\Delta\}}=0$ is contant and hence $\left.\tau^\ast\right|_{\{\Delta\}}$-continuous.
        \end{itemize}
    \end{proof}
\end{proposition}

\bigbreak\noindent We have a result, analogous to \autoref{prop:restriction-of-function-on-one-point-extension}, for the restriction of a continuous function on $(E^\ast,\tau^\ast)$ to $(E,\tau)$ as well:

\begin{proposition}\label{prop:restriction-of-function-on-alexandroff-one-point-extension}
    Let $f^\ast:E^\ast\to\mathbb R$ $\Rightarrow$ \begin{enumerate}[(i)]
        \item $f^\ast$ is $\tau^\ast$-continuous;
        \item $\left.f^\ast\right|_E$ is $\tau$-continuous
    \end{enumerate} are equivalent.
    \begin{proof}[Proof\textup:\nopunct]
        \leavevmode
        \begin{itemize}[$\circ$]
            \item Assume $\left.f^\ast\right|_E$ is $\tau$-continuous.
            \item \autoref{prop:subspace-topology-in-one-point-extension-is-original-topology} $\Rightarrow$ \begin{equation}
                \left.\tau^\ast\right|_E=\tau
            \end{equation} and hence $\left.f^\ast\right|_E$ is $\left.|\tau^\ast\right|_E$-continuous.
            \item $\left.f^\ast\right|_{\{\Delta\}}$ is constant and hence $\left.\tau^\ast\right|_{\{\Delta\}}$-continuous.
        \end{itemize}        
    \end{proof}
\end{proposition}

\section{lifetime formalism}\label{chap:lifetime-formalism}

\subsection{One-point extension of a Markov semigroup}

Let \begin{itemize}[$\circ$]
    \item $(E,\mathcal E)$ be a measurable space and $\Delta\not\in E$;
    \item $(E^\ast,\mathcal E^\ast)$ denote the one-point extension\footnote{see \autoref{def:one-point-extension-of-measurable-space}} of $(E,\mathcal E)$ by $\Delta$.
\end{itemize}

\begin{proposition}
    Let $\kappa^\ast$ be a Markov kernel on $(E^\ast,\mathcal E^\ast)$ $\Rightarrow$ \begin{equation}
        \kappa(x,B):=\kappa^\ast(x,B)\;\;\;\text{for }(x,B)\in(E,\mathcal E)
    \end{equation} is a sub-Markov kernel on $(E,\mathcal E)$ with \begin{equation}
        \kappa^\ast(x^\ast,B^\ast)=\kappa(x^\ast,\underbrace{B^\ast\cap E}_{\in\:\mathcal E})+1_{B^\ast}(\delta)\left(1-\kappa(x^\ast, E)\right)
    \end{equation} for all $(x^\ast,B^\ast)\in(E^\ast,\mathcal E^\ast)$.
    \begin{proof}[Proof\textup:\nopunct]
        Let $(x^\ast,B^\ast)\in(E^\ast,\mathcal E^\ast)$ $\Rightarrow$ \begin{equation}
            \begin{split}
                \kappa^\ast(x^\ast,B^\ast)&=\kappa^\ast(x^\ast,B^\ast\cap E)+\kappa^\ast(x^\ast,B^\ast\cap\{\Delta\})\\&=\begin{cases}\kappa^\ast(x^\ast,B^\ast\cap E)+\kappa^\ast(x^\ast,\{\Delta\})&\text{, if }\Delta\in B^\ast;\\\kappa^\ast(x^\ast,B^\ast\cap E)+\kappa^\ast(x^\ast,\emptyset)&\text{, otherwise}\end{cases}\\&=\kappa(x^\ast,B^\ast\cap E)+1_{B^\ast}(\Delta)\left(1-\kappa(x^\ast,E)\right).
            \end{split}
        \end{equation}
    \end{proof}
\end{proposition}

\begin{proposition}\label{prop:markov-kernel-on-one-point-extension-induced-by-sub-markov-kernel-on-original-space}
    Let $\kappa$ be a sub-Markov kernel on $(E,\mathcal E)$ $\Rightarrow$ \begin{equation}
        \kappa^\ast(x^\ast,B^\ast):=\kappa(x^\ast,\underbrace{B^\ast\cap E}_{\in\:\mathcal E})+1_{B^\ast}(\Delta)\left(1-\kappa(x^\ast,E)\right)
    \end{equation} for $(x^\ast,B^\ast)\in(E^\ast,\mathcal E^\ast)$ is a Markov kernel on $(E^\ast,\mathcal E^\ast)$ with \begin{equation}\label{eq:markov-kernel-on-one-point-extension-induced-by-sub-markov-kernel-on-original-space-eq1}
        \kappa^\ast f^\ast=\left.\kappa f\right|_E+f^\ast(\Delta)\left(1-\kappa(\;\cdot\;,E)\right)=\kappa\left(\left.f^\ast\right|_E-f^\ast(\Delta)\right)+f^\ast(\Delta)
    \end{equation} for all $f^\ast\in\mathcal E^\ast_b$.
    \bigbreak\noindent\textit{Remark}:\begin{enumerate}[(i)]
        \item \autoref{prop:trace-of-one-point-extension-is-original-sigma-algebra} $\Rightarrow$ \begin{equation}
            B^\ast\cap E\in\left.\mathcal E\right|_E=\mathcal E\;\;\;\text{for all }B^\ast\in\mathcal E^\ast.
        \end{equation}
        \item \autoref{prop:function-extension-to-one-point-extension} $\Rightarrow$ \begin{equation}
            \kappa(\Delta,B)\stackrel{\text{def}}=(\kappa1_B)(\Delta)\stackrel{\text{def}}=0\;\;\;\text{for all }B\in\mathcal E
        \end{equation} and hence \begin{equation}
            \kappa^\ast(\Delta,B^\ast)=1_{B^\ast}(\Delta)\;\;\;\text{for all }B^\ast\in\mathcal E^\ast.
        \end{equation}
        \item If $\kappa$ is a Markov kernel on $(E,\mathcal E)$, then \begin{equation}
            \kappa^\ast(x^\ast,B^\ast)=\begin{cases}\kappa(x^\ast,B^\ast\cap E)&\text{, if }x^\ast\in E\\1_{B^\ast}(\Delta)&\text{, otherwise}\end{cases}
        \end{equation} for all $(x^\ast,B^\ast)\in(E^\ast,\mathcal E^\ast)$.
        \item Let $x^\ast\in E^\ast$ $\Rightarrow$ \begin{equation}\label{prop:markov-kernel-on-one-point-extension-induced-by-sub-markov-kernel-on-original-space-remark-iv}
            \kappa^\ast(x^\ast,\{\Delta\})=1-\underbrace{\kappa(x^\ast,E)}_{=\:0,\text{ if }x^\ast\:=\:\Delta}.
        \end{equation}
    \end{enumerate}
    \begin{proof}[Proof\textup:\nopunct]
        \leavevmode
        \begin{itemize}[$\circ$]
            \item \begin{listclaim}
                $\kappa^\ast$ is a Markov kernel on $(E^\ast,\mathcal E^\ast)$.
                \begin{proof}[Proof\textup:\nopunct]
                    \begin{equation}
                        \kappa^\ast(x^\ast,E^\ast)\stackrel{\text{def}}=\kappa(x^\ast,\underbrace{E^\ast\cap E}_{=\:E})+\underbrace{1_{E^\ast}(\Delta)}_{=\:1}\left(1-\kappa(x^\ast,E)\right)=1.
                    \end{equation}
                \end{proof}
            \end{listclaim}
            \item \begin{listclaim}
                Let $f^\ast\in\mathcal E^\ast_b$ $\Rightarrow$
                \begin{proof}[Proof\textup:\nopunct]
                    Assume $k:=|f^\ast(E^\ast)|\in\mathbb N$ \begin{equation}
                        f^\ast=\sum_{i=1}^kf^\ast_i1_{B^\ast_i}
                    \end{equation} for some $f^\ast_1,\ldots,f^\ast_k\in\mathbb R$ and disjoint $B^\ast_1,\ldots,B^\ast_k\in\mathcal E^\ast$ $\Rightarrow$ \begin{equation}
                        \begin{split}
                            (\kappa^\ast f^\ast)(x^\ast)&=\sum_{i=1}^kf^\ast_i\kappa^\ast(x^\ast,B^\ast_i)\\&\stackrel{\text{def}}=\sum_{i=1}^kf^\ast_i\left(\kappa(x^\ast,B^\ast_i\cap E)+1_{B^\ast_i}(\Delta)\left(1-\kappa(x^\ast,E)\right)\right)\\&=\int\kappa(x^\ast,\dif y)\underbrace{\sum_{i=1}^kf^\ast_i\underbrace{1_{B^\ast_i\cap E}(y)}_{=\:1_{B^\ast_i}(y)}}_{=\:\left.f^\ast\right|_E(y)}+\underbrace{\sum_{i=1}^kf^\ast_i1_{B^\ast_i}(\Delta)}_{=\:f^\ast(\Delta)}\left(1-\kappa(x^\ast,E)\right)\\&=(\left.\kappa f^\ast\right|_E)(x^\ast)+f^\ast(\Delta)\left(1-\kappa(x^\ast,E)\right)\\&=\left(\kappa\left(\left.f^\ast\right|_E-f^\ast(\Delta)\right)\right)(x^\ast)+f^\ast(\Delta)
                        \end{split}
                    \end{equation} for all $x^\ast\in E^\ast$.
                \end{proof}
            \end{listclaim}
        \end{itemize}
    \end{proof}
\end{proposition}

\begin{definition}\label{def:one-point-extension-of-markov-kernel}
    $\kappa^\ast$ is called \textbf{one-point extension of} $\bm\kappa$ \textbf{by} $\bm\Delta$.
    \begin{flushright}
        $\square$
    \end{flushright}
\end{definition}

\bigbreak\noindent Now let $(\kappa_t)_{t\ge0}$ be a sub-Markov semigroup on $(E,\mathcal E)$.

\begin{corollary}
    \begin{equation}
        \kappa_t^\ast(x^\ast,B^\ast):=\kappa_t(x^\ast,B^\ast\cap E)+1_{B^\ast}(\Delta)\left(1-\kappa_t(x^\ast,E)\right)
    \end{equation} for $(x^\ast,B^\ast)\in(E^\ast,\mathcal E^\ast)$ and $t\ge0$ is a Markov semigroup on $(E^\ast,\mathcal E^\ast)$.
    \begin{proof}[Proof\textup:\nopunct]
        \leavevmode
        \begin{itemize}[$\circ$]
            \item \autoref{prop:markov-kernel-on-one-point-extension-induced-by-sub-markov-kernel-on-original-space} $\Rightarrow$ $\kappa_t^\ast$ is a Markov kernel on $(E^\ast,\mathcal E^\ast)$ for all $t\ge0.$
            \item The semigroup property of $(\kappa^\ast)_{t\ge0}$ is easily verified.
        \end{itemize}
    \end{proof}
\end{corollary}

\begin{definition}\label{def:one-point-extension-of-markov-semigroup}
    $(\kappa_t^\ast)_{t\ge0}$ is called \textbf{one-point extension of} $\bm{(\kappa_t)_{t\ge0}}$ \textbf{by} $\bm\Delta$.
    \begin{flushright}
        $\square$
    \end{flushright}
\end{definition}

\bigbreak\noindent Let \begin{itemize}[$\circ$]
    \item $(\Omega,\mathcal A)$ be a measurable space;
    \item $(\mathcal F_t)_{t\ge0}$ be a filtration on $(\Omega,\mathcal A)$;
    \item $(X_t)_{t\ge0}$ be an $(E^\ast,\mathcal E^\ast)$-valued $(\mathcal F_t)_{t\ge0}$-adapted process on $(\Omega,\mathcal A)$.
\end{itemize}

\begin{proposition}\label{prop:markov-property-implies-markov-property-on-one-point-extension}
    Let $\operatorname P\in\mathcal M_1(\Omega,\mathcal A)$ and $s,t\ge0$ $\Rightarrow$ If \begin{equation}\label{eq:markov-property-implies-markov-property-on-one-point-extension-eq1}
        \operatorname E\left[f(X_{s+t})\mid\mathcal F_s\right]=(\kappa_tf)(X_s)
    \end{equation} for all $f\in\mathcal E_b$\footnote{see \autoref{prop:function-extension-to-one-point-extension}.}, then \begin{equation}
        \operatorname E\left[f^\ast(X_{s+t})\mid\mathcal F_s\right]=(\kappa_t^\ast f^\ast)(X_s)
    \end{equation} for all $f^\ast\in\mathcal E^\ast_b$.
    \begin{proof}[Proof\textup:\nopunct]
        \leavevmode
        \begin{enumerate}[$\circ$]
            \item Let $f^\ast\in\mathcal E^\ast_b$.
            \item \autoref{prop:restriction-of-function-on-one-point-extension} $\Rightarrow$ $\left.f^\ast\right|_E$ is $\mathcal E$-measurable $\Rightarrow$ \begin{equation}
                \operatorname E\left[\left.f^\ast\right|_E(X_{s+t})\mid\mathcal F_s\right]\stackrel{\eqref{eq:markov-property-implies-markov-property-on-one-point-extension-eq1}}=\left(\left.\kappa_tf^\ast\right|_E\right)(X_s)
            \end{equation}
            \item Note that \begin{equation}\label{eq:markov-property-implies-markov-property-on-one-point-extension-eq2}
                \operatorname P\left[X_{s+t}=\Delta\mid\mathcal F_s\right]=1-\underbrace{\operatorname P\left[X_{s+t}\in E\mid\mathcal F_s\right]}_{\stackrel{\eqref{eq:markov-property-implies-markov-property-on-one-point-extension-eq1}}=\:\kappa_t(X_s,\:E)}=1-\kappa_t(X_s,E)
            \end{equation} and hence\footnote{Remember that \begin{equation}
                \left.f^\ast\right|_E(\Delta)\stackrel{\text{def}}=0;
            \end{equation} see \autoref{prop:function-extension-to-one-point-extension}.} \begin{equation}
                \begin{array}{@{}r@{}c@{}l@{}}
                    \operatorname E\left[f^\ast(X_{s+t})\mid\mathcal F_s\right]&=&\operatorname E\left[f^\ast(X_{s+t});X_{s+t}\in E\mid\mathcal F_s\right]\\&\hphantom=&\;\;\;\;\;\;\;\;\;\;\;\;+\operatorname E\left[f^\ast(X_{s+t});X_{s+t}=\Delta\mid\mathcal F_s\right]\\&=&\operatorname E\left[\left.f^\ast\right|_E(X_{s+t})\mid\mathcal F_s\right]+f^\ast(\Delta)\operatorname P\left[X_{s+t}=\Delta\mid\mathcal F_s\right]\\&=&(\left.\kappa_tf^\ast\right|_E)(X_s)+f^\ast(\Delta)\left(1-\kappa_t(X_s,E)\right)\\&\stackrel{\eqref{eq:markov-kernel-on-one-point-extension-induced-by-sub-markov-kernel-on-original-space-eq1}}=&(\kappa_t^\ast f^\ast)(X_s).
                \end{array}
            \end{equation}
        \end{enumerate}
    \end{proof}
\end{proposition}

\subsection{Lifetime of a path}

Let \begin{itemize}[$\circ$]
    \item $E$ be a set;
    \item $\Delta\not\in E$ and $E^\ast:=E\uplus\{\Delta\}$;
    \item $x:^\ast:[0,\infty)\to E^\ast$.
\end{itemize}

\begin{definition}
    \begin{equation}
        \zeta:=\inf\underbrace{\left\{t\ge0:x^\ast(t)=\Delta\right\}}_{=:\:I}
    \end{equation} is called \textbf{lifetime of} $\bm{x^\ast}$.
    \begin{flushright}
        $\square$
    \end{flushright}
\end{definition}

\begin{lemma}\label{lem:path-before-end-of-lifetime}
    \begin{equation}
        x(t)\in E\;\;\;\text{for all }t\in[0,\zeta).
    \end{equation}
    \begin{proof}[Proof\textup:\nopunct]
        \leavevmode
        \begin{itemize}[$\circ$]
            \item Assume the contrary $\Rightarrow$ $\exists t\in[0,\zeta)\cap I$.
            \item $\zeta$ is a lower bound for $I$ $\Rightarrow$ \begin{equation}
                \zeta\le t;
            \end{equation} in contradiction to $t<\zeta$.
        \end{itemize}
    \end{proof}
\end{lemma}

Let \begin{itemize}[$\circ$]
    \item $\tau$ be a topology on $E$;
    \item $\tau^\ast$ denote the Alexandroff one-point extension\footnote{see \autoref{def:alexanroff-one-point-extension}.} of $\tau$ by $\Delta$.
\end{itemize}

\begin{lemma}
    If $x^\ast$ is $\tau^\ast$-right-continuous, then \begin{equation}\label{eq:weak-right-continuity}
        \forall t\in[0,\infty)\setminus I:\exists\delta>0:x^\ast\left((t,t+\delta)\right)\subseteq E.
    \end{equation}
    \begin{proof}[Proof\textup:\nopunct]
        \leavevmode
        \begin{itemize}[$\circ$]
            \item Let $t\in[0,\infty)\setminus I$ $\Rightarrow$ \begin{equation}
                x^\ast(t)\in E\in\tau\subseteq\tau^\ast
            \end{equation} and hence $N$ is a $\tau^\ast$-neighborhood of $x^\ast(t)$.
            \item $x^\ast$ is $\tau^\ast$-right-continuous at $t$ $\Rightarrow$ \begin{equation}
                x^\ast\left((t,t+\delta)\right)\subseteq E.
            \end{equation}
        \end{itemize}
    \end{proof}
\end{lemma}

\noindent Assume \eqref{eq:weak-right-continuity} is satisfied. The following result is similar to the corresponding claim about the hitting time of a closed set; see \autoref{lem:hitting-time-of-closed-set}.

\begin{lemma}\label{lem:lifetime-finiteness}
    Either \begin{enumerate}[(i)]
        \item $\zeta\in I$; or
        \item $I=\emptyset$ and hence $\zeta=\infty$.
    \end{enumerate}
    \begin{proof}[Proof\textup:\nopunct]
        \leavevmode
        \begin{itemize}[$\circ$]
            \item If $I=\emptyset$, then \begin{equation}
                \zeta=\inf\emptyset=\infty.
            \end{equation}
            \item Assume $I\ne\emptyset$ $\Rightarrow$ \begin{equation}
                \zeta<\infty.
            \end{equation}
            \item Assume $\zeta\not\in I$ $\Rightarrow$ $\exists\delta>0$ with \begin{equation}
                x^\ast\left([\zeta,\zeta+\delta)\right)\subseteq E
            \end{equation} and hence \begin{equation}
                [\zeta,\zeta+\delta)\cap I=\emptyset.
            \end{equation}
            \item $\zeta=\inf I$ is a lower bound for $I$ $\Rightarrow$ \begin{equation}
                I\subseteq[\zeta,\infty)
            \end{equation} and hence \begin{equation}
                \zeta\ge\zeta+\delta;
            \end{equation} which is impossible.
        \end{itemize}
    \end{proof}
\end{lemma}

\begin{lemma}\label{lem:lifetime-smaller-than-t}
    Let $t\ge0$ $\Rightarrow$ \begin{equation}
        \zeta\le t\Leftrightarrow\exists s\in[0,t]:x^\ast(s)=\Delta
    \end{equation} and hence \begin{equation}\label{eq:lifetime-smaller-than-t-2}
        \zeta>t\Leftrightarrow\forall s\in[0,t]:x^\ast(s)\in E.
    \end{equation}
    \begin{proof}[Proof\textup:\nopunct]
        \leavevmode
        \begin{itemize}[$\circ$]
            \item\begin{listclaim}
                "$\Rightarrow$"
                \begin{proof}[Proof\textup:\nopunct]
                    \leavevmode
                    \begin{itemize}[$\circ$]
                        \item Assume the contrary $\Rightarrow$ \begin{equation}
                            \forall s\in I:s>t
                        \end{equation} and hence $t$ is a lower bound for $I$.
                        \item $\zeta=\inf I$ $\Rightarrow$ \begin{equation}
                            t\le\zeta.
                        \end{equation}
                        \item $\zeta\le t$ $\Rightarrow$ \begin{equation}
                            \zeta=t.
                        \end{equation} and hence \begin{equation}
                            \zeta\not\in I.
                        \end{equation}
                        \item \autoref{lem:lifetime-finiteness} $\Rightarrow$ $I=\emptyset$ and hence $\zeta=\infty$; in contradiction to $\zeta\le t<\infty$.
                    \end{itemize}
                \end{proof}
            \end{listclaim}
            \item\begin{listclaim}
                "$\Leftarrow$"
                \begin{proof}[Proof\textup:\nopunct]
                    \leavevmode
                    \begin{itemize}[$\circ$]
                        \item Let $s\in I$ with $s\le t$.
                        \item $\zeta$ is a lower bound for $I$ $\Rightarrow$ $\zeta\le s\le t$.
                    \end{itemize}
                \end{proof}
            \end{listclaim}
        \end{itemize}
    \end{proof}
\end{lemma}

\begin{lemma}\label{lem:lifetime-when-trap-is-absorbing}
    If $\Delta$ is absorbing\footnote{see \autoref{def:trap}.} for $x^\ast$, then \begin{enumerate}[(i)]
        \item \begin{equation}\label{eq:trap-at-lifetime}
            x^\ast(t)=\Delta\;\;\;\text{for all }t\in[\zeta,\infty).
        \end{equation}
        \textit{Remark}: If \begin{equation}
            x^\ast(\infty):=\Delta,
        \end{equation} then \begin{equation}
            x^\ast(t)=x^\ast(\zeta)=\Delta\;\;\;\text{for all }t\in[\zeta,\infty].
        \end{equation}
        \item Let $t\ge0$ $\Rightarrow$ \begin{equation}\label{eq:lifetime-when-trap-is-absorbing-eq2}
            \zeta>t\Leftrightarrow x^\ast(t)\in E.
        \end{equation}
    \end{enumerate}
    \begin{proof}[Proof\textup:\nopunct]
        \leavevmode
        \begin{enumerate}[(i)]
            \item\leavevmode\begin{itemize}[$\circ$]
                \item \autoref{lem:lifetime-finiteness} $\Rightarrow$ Either $\zeta\in I$ or $I=\emptyset$.
                \item If $\zeta\in I$, then \begin{equation}
                    x^\ast(\zeta)=\Delta
                \end{equation} and hence \eqref{eq:trap-at-lifetime}; since $\Delta$ is absorbing for $x^\ast$.
                \item If $I=\emptyset$, then $\zeta=\infty$ and hence \eqref{eq:trap-at-lifetime} is trivial.
            \end{itemize}
            \item\leavevmode\begin{itemize}[$\circ$]
                \item If $\zeta>t$, then\footnote{Assume $x^\ast(t)=\Delta$ $\Rightarrow$ $t\in I$. $\zeta$ is a lower bound for $I$ $\Rightarrow$ $\zeta\le t$.} \begin{equation}
                    x^\ast(t)\in E.
                \end{equation}
                \item Assume $x^\ast(t)\in E$ and $\zeta\le t$.
                \item $\zeta\le t$ and \autoref{lem:lifetime-smaller-than-t} $\Rightarrow$ $\exists s\in[0,t]$ with \begin{equation}
                    x^\ast(s)=\Delta.
                \end{equation}
                \item $\Delta$ is absorbing for $x^\ast$ and $t\ge s$ $\Rightarrow$ \begin{equation}
                    x^\ast(t)=\Delta\not\in E;
                \end{equation} in contradiction to $x^\ast(t)\in E$.
            \end{itemize}
        \end{enumerate}
    \end{proof}
\end{lemma}

\subsection{Markov processes with restricted lifetime}

Let \begin{itemize}[$\circ$]
    \item $(E,\mathcal E)$ be a measurable space;
    \item $(\kappa_t)_{t\ge0}$ be a sub-Markov semigroup on $(E,\mathcal E)$;
    \item $\Delta\not\in E$ and $(E^\ast,\mathcal E^\ast)$ denote the one-point extension\footnote{see \autoref{def:one-point-extension-of-measurable-space}.} of $(E,\mathcal E)$ by $\Delta$;
    \item $(\kappa^\ast_t)_{t\ge0}$ denote the one-point extension\footnote{see \autoref{def:one-point-extension-of-markov-semigroup}.} of $(\kappa_t)_{t\ge0}$ by $\Delta$;
    $(\Omega,\mathcal A,\operatorname P)$ be a probability space;
    \item $(\mathcal F_t)_{t\ge0}$ be a filtration on $(\Omega,\mathcal A)$;
    \item $(X^\ast_t)_{t\ge0}$ be an $(E^\ast,\mathcal E^\ast)$-valued $(\mathcal F_t)_{t\ge0}$-adapted process on $(\Omega,\mathcal A,\operatorname P)$ with \begin{equation}\label{eq:restricted-lifetime-markov-property}
        \operatorname E\left[f(X^\ast_{s+t})\mid\mathcal F_s\right]=(\kappa_tf)(X^\ast_s)\;\;\;\text{or all }f\in\mathcal E_b\text{ and }s,t\ge0.
    \end{equation}
\end{itemize}

\begin{remark}\label{rem:markov-property-extension}
    \eqref{eq:restricted-lifetime-markov-property} and \autoref{prop:markov-property-implies-markov-property-on-one-point-extension} $\Rightarrow$ \begin{equation}\label{eq:restricted-lifetime-extended-markov-property}
        \operatorname E\left[f^\ast(X^\ast_{s+t})\mid\mathcal F_s\right]=(\kappa_t^\ast f^\ast)(X^\ast_s)\;\;\;\text{or all }f^\ast\in\mathcal E^\ast_b\text{ and }s,t\ge0.
    \end{equation}
    \begin{flushright}
        $\square$
    \end{flushright}
\end{remark}

\begin{proposition}\label{prop:markov-process-end-of-lifetime-prop1}
    \leavevmode
    \begin{enumerate}[(i)]
        \item Let $t\ge0$ $\Rightarrow$ \begin{equation}\label{eq:markov-process-end-of-lifetime-prop1-eq1}
            \operatorname P\left[X^\ast_t=\Delta\right]=1-\operatorname E\left[\kappa_t(X^\ast_0,E)\right]
        \end{equation} and hence \begin{equation}\label{eq:markov-process-end-of-lifetime-prop1-eq2}
            \operatorname P\left[X^\ast_t\in E\right]=\operatorname E\left[\kappa_t(X^\ast_0,E)\right].
        \end{equation}
        \textit{Remark}: If $X^\ast_0=\Delta$ almost surely, then \begin{equation}\label{eq:markov-process-end-of-lifetime-prop1-eq3}
            \operatorname E\left[\kappa_t(X^\ast_0,E)\right]=\kappa_t(\Delta,E)=0.
        \end{equation}
        \item If $X^\ast_0=\Delta$ almost surely, then \begin{equation}\label{eq:markov-process-end-of-lifetime-prop1-eq4}
            \operatorname P\left[X^\ast_t=\Delta\text{ for all }t\ge0\right]=1.
        \end{equation}
    \end{enumerate}
    \begin{proof}[Proof\textup:\nopunct]
        \leavevmode
        \begin{enumerate}[(i)]
            \item \begin{equation}
                \operatorname P\left[X^\ast_t\in E\right]\stackrel{\eqref{eq:restricted-lifetime-extended-markov-property}}=\operatorname E\left[\kappa_t^\ast(X^\ast_0,\{\Delta\})\right]\stackrel{\eqref{prop:markov-kernel-on-one-point-extension-induced-by-sub-markov-kernel-on-original-space-remark-iv}}=1-\operatorname E\left[\kappa_t(X^\ast_0,E)\right].
            \end{equation}
            \item If $X^\ast_0=\Delta$ almost surely, then \begin{equation}
                \begin{split}
                    &\operatorname P\left[\left(X^\ast_{t_1},\ldots,X^\ast_{t_n}\right)\in\bigtimes_{i=1}^nB_i\right]\\&\;\;\;\;\;\;\;\;\;\;\;\;=\int\underbrace{\operatorname P\left[X^\ast_0\in\dif x^\ast_0\right]}_{=\:\delta_\Delta(\dif x^\ast_0)}\int_{B_1}\kappa_{t_1}^\ast(x^\ast_0,\dif x^\ast_1)\\&\;\;\;\;\;\;\;\;\;\;\;\;\;\;\;\;\;\;\;\;\;\;\;\;\int_{B_2}\kappa^\ast_{t_2-t_1}(x^\ast_1,\dif x^\ast_2)\cdots\int_{B_n}\kappa^\ast_{t_n-t_{n-1}}(x^\ast_{n-1},\dif x^\ast_n)\\&\;\;\;\;\;\;\;\;\;\;\;\;=\int_{B_1}\underbrace{\kappa_{t_1}^\ast(\Delta,\dif x^\ast_1)}_{=\:\delta_\Delta(\dif x^\ast_1)}\int_{B_2}\\&\;\;\;\;\;\;\;\;\;\;\;\;\;\;\;\;\;\;\;\;\;\;\;\;\kappa^\ast_{t_2-t_1}(x^\ast_1,\dif x^\ast_2)\cdots\int_{B_n}\kappa^\ast_{t_n-t_{n-1}}(x^\ast_{n-1},\dif x^\ast_n)\\&\;\;\;\;\;\;\;\;\;\;\;\;=1_{B_1}(\Delta)\int_{B_2}\underbrace{\kappa^\ast_{t_2-t_1}(\Delta,\dif x^\ast_2)}_{=\:\delta_\Delta(\dif x^\ast_1)}\cdots\int_{B_n}\kappa^\ast_{t_n-t_{n-1}}(x^\ast_{n-1},\dif x^\ast_n)\\&\;\;\;\;\;\;\;\;\;\;\;\;=\prod_{i=1}^n1_{B_i}(\Delta)
                \end{split}
            \end{equation} by \eqref{eq:restricted-lifetime-extended-markov-property}
            for all $n\in\mathbb N$, $0\le t_1<\cdots<t_n$ and $B_1,\ldots,B_n\in\mathcal E$ $\Rightarrow$ Claim.
        \end{enumerate}
    \end{proof}
\end{proposition}

\noindent Assume \begin{equation}
    \{x\}\in\mathcal E\;\;\;\text{for all }x\in E.
\end{equation}

\noindent Let $(\operatorname P_{x^\ast})_{x^\ast\in E^\ast}\subseteq\mathcal M_1(\Omega,\mathcal A)$ with \begin{equation}
    \operatorname P_{x^\ast}\left[X^\ast_0=x^\ast\right]=1
\end{equation} and \eqref{eq:restricted-lifetime-markov-property} for $\operatorname P$ replaced by $\operatorname P_{x^\ast}$ for all $x^\ast\in E^\ast$. Assume \begin{equation}
    E\to[0,1]\;,\;\;\;x\mapsto\operatorname P_x[A]
\end{equation} is $\mathcal E$-measurable for all $A\in\mathcal A$.

\begin{remark}\label{rem:kernel-induced-by-probability-measure}
    \autoref{prop:restriction-of-function-on-one-point-extension} $\Rightarrow$ \begin{equation}
        E^\ast\to[0,1]\;,\;\;\;x^\ast\mapsto\operatorname P_{x^\ast}[A]
    \end{equation} is $\mathcal E^\ast$-measurable for all $A\in\mathcal A$.
    \begin{flushright}
        $\square$
    \end{flushright}
\end{remark}

\begin{proposition}\label{prop:markov-process-end-of-lifetime-prop2}
    \leavevmode
    \begin{enumerate}[(i)]
        \item\label{prop:markov-process-end-of-lifetime-prop2-i} Let $s\ge0$ $\Rightarrow$ \begin{equation}
            \operatorname P\left[X^\ast_t=\Delta\text{ for all }t\ge s\right]=\operatorname P\left[X^\ast_s=\Delta\right]
        \end{equation} and hence \begin{equation}
            \operatorname P\left[X^\ast_s=\Delta\Rightarrow\forall t\ge s:X^\ast_t=\Delta\right]=1;
        \end{equation} i.e. $\Delta$ is almost surely absorbing for $(X^\ast_t)_{t\ge0}$.
        \smallbreak\textit{Remark}: If $X^\ast_s=\Delta$ almost surely, then \begin{equation}
            \operatorname P\left[X^\ast_t=\Delta\text{ for all }t\ge s\right]=1.
        \end{equation}
        \item\label{prop:markov-process-end-of-lifetime-prop2-ii} Let $\tau$ be an $(\mathcal F_t)_{t\ge0}$-stopping time on $(\Omega,\mathcal A)$ $\Rightarrow$ If $X^\ast_\tau=\Delta$ almost surely on $\{\:\tau<\infty\:\}$ and $(X^\ast_t)_{t\ge0}$ is strongly $(\mathcal F_t)_{t\ge0}$-Markov at $\tau$, then \begin{equation}
            \operatorname P\left[X^\ast_t=\Delta\text{ for all }t\in[\tau,\infty)\right]=1.
        \end{equation}
    \end{enumerate}
    \begin{proof}[Proof\textup:\nopunct]
        \leavevmode
        \begin{enumerate}[(i)]
            \item\leavevmode\begin{itemize}[$\circ$]
                \item Let \begin{equation}\label{eq:markov-process-end-of-lifetime-prop2-proof-eq0}
                    B:=\left\{x^\ast\in\left(E^\ast\right)^{[0,\:\infty)}:x^\ast(t)=\Delta\text{ for all }t\ge0\right\}.
                \end{equation}
                \item \eqref{eq:restricted-lifetime-extended-markov-property} $\Rightarrow$ \begin{equation}\label{eq:markov-process-end-of-lifetime-prop2-proof-eq1}
                    \operatorname P\left[(X^\ast_{s+t})_{t\ge0}\in B\mid\mathcal F_s\right]=\operatorname P_{X^\ast_s}\left[X^\ast\in B\right]
                \end{equation} and hence \begin{equation}\label{eq:markov-process-end-of-lifetime-prop2-proof-eq2}
                    \begin{split}
                        &\operatorname P\left[X^\ast_t=\Delta\text{ for all }t\ge s\right]\\&\;\;\;\;\;\;\;\;\;\;\;\;\begin{array}{cl}
                            =&\operatorname P\left[X^\ast_t=\Delta\text{ for all }t\ge s;X^\ast_s=\Delta\right]\\=&\operatorname P\left[(X^\ast_{s+t})_{t\ge0}\in B;X^\ast_s=\Delta\right]\\=&\operatorname E\left[\operatorname P\left[(X^\ast_{s+t})_{t\ge0}\in B\mid\mathcal F_s\right];X^\ast_s=\Delta\right]\\\stackrel{\eqref{eq:markov-process-end-of-lifetime-prop2-proof-eq1}}=&\operatorname E\left[\operatorname P_{X^\ast_s}\left[X^\ast\in B\right];X^\ast_s=\Delta\right]\\=&\underbrace{\operatorname P_\Delta\left[X^\ast\in B\right]}_{\stackrel{\eqref{eq:markov-process-end-of-lifetime-prop1-eq4}}=\:1}\operatorname P\left[X^\ast_s=\Delta\right].
                        \end{array}
                    \end{split}
                \end{equation}
                \item \eqref{eq:markov-process-end-of-lifetime-prop2-proof-eq2} $\Rightarrow$ \begin{equation}
                    \begin{split}
                        &\operatorname P\left[X^\ast_s=\Delta\Rightarrow\forall t\ge s:X^\ast_t=\Delta\right]\\&\;\;\;\;\;\;\;\;\;\;\;\;=\operatorname P\left[\{X^\ast_s\in E\}\uplus\{\forall t\ge s:X^\ast_t=\Delta\}\right]\\&\;\;\;\;\;\;\;\;\;\;\;\;=\operatorname P\left[X^\ast_s\in E\right]+\operatorname P\left[\forall t\ge s:X^\ast_t=\Delta\right]\\&\;\;\;\;\;\;\;\;\;\;\;\;=\operatorname P\left[X^\ast_s\in E\right]+\operatorname P\left[X^\ast_s=\Delta\right]=1.
                    \end{split}
                \end{equation}
                \item If $X^\ast_s=\Delta$ almost surely, then \begin{equation}
                    \operatorname P\left[(X^\ast_{s+t})_{t\ge0}\in B\right]=\operatorname E\left[\operatorname P_{X^\ast_s}\left[X^\ast\in B\right]\right]=\operatorname P\left[X^\ast\in B\right]\stackrel{\eqref{eq:markov-process-end-of-lifetime-prop1-eq4}}=1.
                \end{equation}
            \end{itemize}
            \item\leavevmode\begin{itemize}[$\circ$]
                \item Define $B$ as in \eqref{eq:markov-process-end-of-lifetime-prop2-proof-eq0}.
                \item $(X^\ast_t)_{t\ge0}$ is strongly $(\mathcal F_t)_{t\ge0}$-Markov at $\tau$ $\Rightarrow$ \begin{equation}
                    \operatorname P\left[(X^\ast_{\tau+t})_{t\ge0}\in B;\tau<\infty\mid\mathcal F_\tau\right]=1_{\{\:\tau\:<\:\infty\:\}}\operatorname P_{X^\ast_\tau}\left[X^\ast\in B\right].
                \end{equation}
                \item $X^\ast_\tau=\Delta$ almost surely on $\{\:\tau<\infty\:\}$ $\Rightarrow$ \begin{equation}
                    \begin{split}
                        &\operatorname P\left[X^\ast_\tau=\Delta\text{ for all }t\in[\tau,\infty);\tau<\infty\right]\\&\;\;\;\;\;\;\;\;\;\;\;\;=\operatorname P\left[(X^\ast_{\tau+t})_{t\ge0}\in B;\tau<\infty\right]\\&\;\;\;\;\;\;\;\;\;\;\;\;=\operatorname P\left[1_{\{\:\tau\:<\:\infty\:\}}\operatorname P_{X^\ast_\tau}\left[X^\ast\in B\right]\right]\\&\;\;\;\;\;\;\;\;\;\;\;\;=\operatorname P\left[\tau<\infty\right]\underbrace{\operatorname P_\Delta\left[X^\ast\in B\right]}_{\stackrel{\eqref{eq:markov-process-end-of-lifetime-prop1-eq4}}=\:1}
                    \end{split}
                \end{equation} and hence \begin{equation}
                    \begin{split}
                    &\operatorname P\left[X^\ast_t=\Delta\text{ for all }t\in[\tau,\infty)\right]\\&\;\;\;\;\;\;\;\;\;\;\;\;=\operatorname P\left[X^\ast_t=\Delta\text{ for all }t\in[\tau,\infty);\tau<\infty\right]+\operatorname P\left[\tau=\infty\right]=1.
                    \end{split}
                \end{equation}
            \end{itemize}
        \end{enumerate}
    \end{proof}
\end{proposition}

\noindent Let $\zeta$ denote the lifetime of $(X^\ast_t)_{t\ge0}$\footnote{i.e. \begin{equation}
    \zeta:=\inf\left\{t\ge0:X^\ast_t=\Delta\right\}.
\end{equation}}. Assume \begin{itemize}[$\circ$]
    \item $E$ is a topological space and $\mathcal E=\mathcal B(E)$;
    \item \eqref{eq:weak-right-continuity} is satisfied for $x^\ast$ replaced by $X^\ast(\omega)$ for all $\omega\in\Omega$.
\end{itemize}

\begin{proposition}[lifetime of Markov process is stopping time]\label{prop:lifetime-of-Markov-process-is-stopping-time}~
    If $(\mathcal F_t)_{t\ge0}$ is $\operatorname P$-complete, then $\zeta$ is an $(\mathcal F_t)_{t\ge0}$-stopping time.
    \begin{proof}[Proof\textup:\nopunct]
        \leavevmode
        \begin{itemize}[$\circ$]
            \item \autoref{prop:markov-process-end-of-lifetime-prop2}-\ref{prop:markov-process-end-of-lifetime-prop2-i} $\Rightarrow$ $\Delta$ is absorbing for $X^\ast(\omega)$ for all $\omega\in\Omega\setminus N$ for some $\operatorname P$-null set $N$
            \item $N$ is a $\operatorname P$-null set $\Rightarrow$ $N\in\mathcal F_0$ $\Rightarrow$ \begin{equation}
                \Omega\setminus N\in\mathcal F_0.
            \end{equation}
            \item Let $t\ge0$ $\Rightarrow$ $\{\:\zeta\le t\:\}\cap N$ is a $\operatorname P$-null set $\Rightarrow$ \begin{equation}
                \{\:\zeta\le t\:\}\cap N\in\mathcal F_0.
            \end{equation}
            \item $\Delta$ is absorbing for $X^\ast(\omega)$ for all $\omega\in\Omega\setminus N$ $\Rightarrow$ \begin{equation}
                \{\:\zeta\le t\:\}\cap(\Omega\setminus N)\stackrel{\eqref{eq:lifetime-smaller-than-t-2}}=\underbrace{\{\:X^\ast_t=\Delta\:\}}_{\in\:\mathcal F_t}\cap\underbrace{(\Omega\setminus N)}_{\in\:\mathcal F_0}\in\mathcal F_t
            \end{equation} and hence \begin{equation}
                \{\:\zeta\le t\:\}=\underbrace{\{\:\zeta\le t\:\}\cap N}_{\in\:\mathcal F_0}\uplus\underbrace{\{\:\zeta\le t\:\}\cap(\Omega\setminus N)}_{\in\:\mathcal F_t}\in\mathcal F_t.
            \end{equation}
        \end{itemize}
    \end{proof}
\end{proposition}

\begin{proposition}\label{prop:markov-process-dead-at-lifetime}
    \leavevmode
    \begin{enumerate}[(i)]
        \item\label{prop:markov-process-dead-at-lifetime-i} \begin{equation}
            X^\ast_\zeta=\Delta\;\;\;\text{on }\{\:\zeta<\infty\:\}
        \end{equation} and \begin{equation}
            X^\ast_t\in E\;\;\;\text{for all }t\ge0\text{ on }\{\:\zeta=\infty\:\}.
        \end{equation}
        \item Let $t\ge0$ $\Rightarrow$ \begin{equation}\label{eq:markov-process-dead-at-lifetime-ii}
            \operatorname P\left[\zeta>t\right]=\operatorname E\left[\kappa_t(X^\ast_0,E)\right].
        \end{equation}
        \item If $\zeta$ is an $(\mathcal F_t)_{t\ge0}$-stopping time and $(X^\ast_t)_{t\ge0}$ is strongly $(\mathcal F_t)_{t\ge0}$-Markov at $\zeta$, then \begin{equation}
            \operatorname P\left[X^\ast_t=\Delta\text{ for all }t\in[\zeta,\infty)\right]=1.
        \end{equation}
    \end{enumerate}
    \begin{proof}[Proof\textup:\nopunct]
        \leavevmode
        \begin{enumerate}[(i)]
            \item \autoref{lem:lifetime-finiteness} $\Rightarrow$ Claim.
            \item \autoref{prop:markov-process-end-of-lifetime-prop2}-\ref{prop:markov-process-end-of-lifetime-prop2-i} $\Rightarrow$ $\Delta$ is almost surely absorbing for $(X^\ast_t)_{t\ge0}$ $\Rightarrow$ \begin{equation}
                \operatorname P\left[\zeta>t\right]\stackrel{\eqref{eq:lifetime-when-trap-is-absorbing-eq2}}=\operatorname P\left[X^\ast_t\in E\right]\stackrel{\eqref{eq:restricted-lifetime-markov-property}}=\operatorname E\left[\kappa_t(X^\ast_0,E)\right].
            \end{equation}
            \item \ref{prop:markov-process-dead-at-lifetime-i} $\Rightarrow$ $X^\ast_\zeta=\Delta$ on $\{\:\zeta<\infty\:\}$ $\Rightarrow$ Claim by \autoref{prop:markov-process-end-of-lifetime-prop2}-\ref{prop:markov-process-end-of-lifetime-prop2-ii}.
        \end{enumerate}
    \end{proof}
\end{proposition}


\newpage
\bibliographystyle{ACM-Reference-Format}
\bibliography{bibliography}

\end{document}